\pdfoutput=1
\RequirePackage{ifpdf}
\ifpdf 
\documentclass[pdftex]{sigma}
\else
\documentclass{sigma}
\fi

\numberwithin{equation}{section}

\newtheorem{Theorem}{Theorem}[section]
\newtheorem*{Theorem*}{Theorem}

\newtheorem{Lemma}[Theorem]{Lemma}
\newtheorem{Conjecture}[Theorem]{Conjecture}
\newtheorem{Proposition}[Theorem]{Proposition}
 { \theoremstyle{definition}
\newtheorem{Definition}[Theorem]{Definition}

\newtheorem{Example}[Theorem]{Example}
\newtheorem{Remark}[Theorem]{Remark}
\newtheorem{Algorithm}[Theorem]{Algorithm} }

\usepackage{tikz-cd}
\usepackage{enumerate}
\usepackage{stmaryrd}

\newcommand{\Hom}{\textup{Hom}}
\newcommand{\id}{\textup{id}}
\newcommand{\GL}{\textup{GL}}
\newcommand{\SL}{\textup{SL}}
\newcommand{\AGL}{\textup{AGL}}
\newcommand{\Sp}{\textup{Sp}}
\newcommand{\SO}{\textup{SO}}
\newcommand{\U}{\mathbb{U}}
\newcommand{\Stab}{\textup{Stab}}
\newcommand{\Span}{\textup{Span}}
\newcommand{\op}{\textup{op}}
\newcommand{\CC}{\mathbb{C}}
\newcommand{\ZZ}{\mathbb{Z}}
\newcommand{\FF}{\mathbb{F}}
\newcommand{\GG}{\mathbb{G}}
\renewcommand{\AA}{\mathbb{A}}
\newcommand{\K}{\textup{K}}
\newcommand{\Var}{\textup{\bf Var}}
\newcommand{\Mod}{\textup{\bf Mod}}
\newcommand{\Grpd}{\textbf{Grpd}}
\newcommand{\Bd}{\textbf{Bd}}
\newcommand{\Bdp}{\textbf{Bdp}}
\newcommand{\Tb}{\textbf{Tb}}
\DeclareMathOperator{\Spec}{Spec}

\newcommand{\bdscale}{0.5}

\def\bdcupleft {
\begin{tikzpicture}[semithick, scale=\bdscale, baseline=-0.5ex]
\begin{scope}
 \draw (0,0) ellipse (0.2cm and 0.4cm);
 \draw (0,-0.4) arc (-90:90:0.75cm and 0.4cm);
\end{scope}
\end{tikzpicture}
}

\def\bdcupright {
\begin{tikzpicture}[semithick, scale=\bdscale, baseline=-0.5ex]
\begin{scope}
 \draw (0,0) ellipse (0.2cm and 0.4cm);
 \draw (0,0.4) arc (90:270:0.75cm and 0.4cm);
\end{scope}
\end{tikzpicture}
}

\def\bdgenus {
\begin{tikzpicture}[semithick, scale=\bdscale, baseline=-0.5ex]
\begin{scope}
 \draw (-1,0) ellipse (0.2cm and 0.4cm);
 \draw (-1,0.4) .. controls (-0.5,0.6) and (0.5,0.6) .. (1,0.4);
 \draw (-1,-0.4) .. controls (-0.5,-0.6) and (0.5,-0.6) .. (1,-0.4);
 \draw (-0.5,0.1) .. controls (-0.5,-0.125) and (0.5,-0.125) .. (0.5,0.1);
 \draw (-0.4,0.0) .. controls (-0.4,0.0625) and (0.4,0.0625) .. (0.4,0.0);
 \draw (1,0) ellipse (0.2cm and 0.4cm);
\end{scope}
\end{tikzpicture}
}

\def\bdgenerictube#1 {
\begin{tikzpicture}[semithick, scale=\bdscale, baseline=-0.5ex]
\begin{scope}
 \draw (-1,0) ellipse (0.2cm and 0.4cm);
 \draw (-1,0.4) cos (-0.875, 0.5) sin (-0.75,0.6) cos (-0.625,0.5) sin (-0.5,0.4) cos (-0.375,0.5) sin (-0.25, 0.6) cos (-0.125,0.5) sin (0, 0.4) cos (0.125,0.5) sin (0.25,0.6) cos (0.375,0.5) sin (0.5,0.4) cos (0.625,0.5) sin (0.75,0.6) cos (0.875,0.5) sin (1,0.4);
 \draw (-1,-0.4) cos (-0.875, -0.5) sin (-0.75,-0.6) cos (-0.625,-0.5) sin (-0.5,-0.4) cos (-0.375,-0.5) sin (-0.25, -0.6) cos (-0.125,-0.5) sin (0,-0.4) cos (0.125,-0.5) sin (0.25,-0.6) cos (0.375,-0.5) sin (0.5,-0.4) cos (0.625,-0.5) sin (0.75,-0.6) cos (0.875,-0.5) sin (1,-0.4);
 \draw (1,0) ellipse (0.2cm and 0.4cm);
 \draw (0,0) node {$#1$};
\end{scope}
\end{tikzpicture}
}

\def\bdpgenus {
\begin{tikzpicture}[semithick, scale=2.0*\bdscale, baseline=-0.5ex]
\begin{scope}
 \draw (-1,0) ellipse (0.2cm and 0.4cm);
 \draw (-1,0.4) .. controls (-0.5,0.45) and (0.5,0.45) .. (1,0.4);
 \draw (-1,-0.4) .. controls (-0.5,-0.45) and (0.5,-0.45) .. (1,-0.4);
 \draw (-0.5,0.1) .. controls (-0.5,-0.125) and (0.5,-0.125) .. (0.5,0.1);
 \draw (-0.4,0.0) .. controls (-0.4,0.0625) and (0.4,0.0625) .. (0.4,0.0);
 \draw (1,0.4) arc (90:-90:0.2cm and 0.4cm);
 \draw[dashed] (1,0.4) arc (90:270:0.2cm and 0.4cm);
 \draw[black,fill=black] (-0.8,0) circle (.4ex);
 \draw[black,fill=black] (1.2,0) circle (.4ex);
\end{scope}
\end{tikzpicture}
}

\def\bdpcupleft {
\begin{tikzpicture}[semithick, scale=2.0*\bdscale, baseline=-0.5ex]
\begin{scope}
 \draw (0,0) ellipse (0.2cm and 0.4cm);
 \draw (0,-0.4) arc (-90:90:0.75cm and 0.4cm);
 \draw[black,fill=black] (0.2,0) circle (.4ex);
\end{scope}
\end{tikzpicture}
}

\def\bdpcupright {
\begin{tikzpicture}[semithick, scale=2.0*\bdscale, baseline=-0.5ex]
\begin{scope}
 \draw (0,0.4) arc (90:-90:0.2cm and 0.4cm);
 \draw (0,0.4) arc (90:270:0.75cm and 0.4cm);
 \draw[dashed] (0,0.4) arc (90:270:0.2cm and 0.4cm);
 \draw[black,fill=black] (0.2,0) circle (.4ex);
\end{scope}
\end{tikzpicture}
}

\def\bdpparabolic {
\begin{tikzpicture}[semithick, scale=2.0*\bdscale, baseline=-0.5ex]
\begin{scope}
 \draw (-1,0) ellipse (0.2cm and 0.4cm);
 \draw (-1,0.4) -- (1,0.4);
 \draw (-1,-0.4) -- (1,-0.4);
 \draw (1,0.4) arc (90:-90:0.2cm and 0.4cm);
 \draw[dashed] (1,0.4) arc (90:270:0.2cm and 0.4cm);
 \draw[black,fill=black] (-0.8,0) circle (.4ex);
 \draw[black,fill=black] (1.2,0) circle (.4ex);
 \draw[black,fill=black] (0,0.15) circle (.4ex);
\end{scope}
\end{tikzpicture}
}

\begin{document}
\allowdisplaybreaks

\newcommand{\arXivNumber}{2008.06679}

\renewcommand{\PaperNumber}{095}

\FirstPageHeading

\ShortArticleName{Virtual Classes of Representation Varieties of Upper Triangular Matrices via TQFTs}

\ArticleName{Virtual Classes of Representation Varieties\\ of Upper Triangular Matrices\\ via Topological Quantum Field Theories}

\Author{M\'arton HABLICSEK and Jesse VOGEL}

\AuthorNameForHeading{M.~Hablicsek and J.~Vogel}

\Address{Mathematical Institute, Niels Bohrweg 1, 2333 CA Leiden, The Netherlands}
\Email{\href{mailto:m.hablicsek@math.leidenuniv.nl}{m.hablicsek@math.leidenuniv.nl}, \href{mailto:j.t.vogel@math.leidenuniv.nl}{j.t.vogel@math.leidenuniv.nl}}

\ArticleDates{Received February 28, 2022, in final form November 28, 2022; Published online December 06, 2022}

\Abstract{In this paper, we use a geometric technique developed by Gonz\'alez-Prieto, Logares, Mu\~noz, and Newstead to study the $G$-representation variety of surface groups $\mathfrak{X}_G(\Sigma_g)$ of arbitrary genus for $G$ being the group of upper triangular matrices of fixed rank. Explicitly, we compute the virtual classes in the Grothendieck ring of varieties of the $G$-representation variety and the moduli space of $G$-representations of surface groups for $G$ being the group of complex upper triangular matrices of rank~$2$, $3$, and~$4$ via constructing a~topological quantum field theory. Furthermore, we show that in the case of upper triangular matrices the character map from the moduli space of $G$-representations to the $G$-character variety is not an isomorphism.}

\Keywords{representation variety; character variety; topological quantum field theory; Gro\-thendieck ring of varieties}

\Classification{14D23; 14D21; 14C30; 14D20; 14D07; 57R56}

\section{Introduction}
Let $X$ be a closed connected manifold with finitely generated fundamental group $\pi_1(X)$, and $G$ an algebraic group over a field $k$. The set of group representations $\rho\colon \pi_1(X) \to G$,
\[ \mathfrak{X}_G(X) = \Hom(\pi_1(X), G) , \]
carries a natural structure of an algebraic variety, and is called the \textit{$G$-representation variety} of~$X$. Indeed, given a set of generators $\gamma_1, \dots, \gamma_n$ of $\pi_1(X)$, the morphism
\[ \mathfrak{X}_G(X) \to G^n, \qquad \rho \mapsto (\rho(\gamma_1), \dots, \rho(\gamma_n)) \]
identifies the $G$-representation variety $\mathfrak{X}_G(X)$ with a subvariety of $G^n$, and this structure can be shown to be independent of the chosen generators. When $X = \Sigma_g$ is a closed oriented surface of genus $g$, the $G$-representation variety is the closed subvariety of $G^{2g}$ given by
\begin{equation} \label{eq:explicit_expression_representation_variety}
 \mathfrak{X}_G(\Sigma_g) = \left\{ (A_1, B_1, \dots, A_g, B_g) \in G^{2g} \,\bigg|\, \prod_{i = 1}^{g} [A_i, B_i] = 1 \right\} .
\end{equation}
The algebraic group $G$ acts by conjugation on the variety $\mathfrak{X}_G(X)$, so one can look at the categorical quotient (the affine GIT quotient)
\[ \mathcal{M}_G(X) = \mathfrak{X}_G(X) \sslash G , \]
which is known as the \emph{moduli space of $G$-representations}.

When the group $G$ is a linear algebraic group ($G\leq \GL_r(k)$ for some $r>0$), there exists another natural variety that parametrizes $G$-representations up to conjugation, which is defined as follows. To a representation $\rho\colon \pi_1(X)\to G$ we can associate its character
\[ \chi_\rho \colon \ \pi_1(X) \to k, \qquad \gamma \mapsto \operatorname{tr}(\rho(\gamma)) . \]
The image of the so-called \emph{character map}
\[ \chi \colon \ \mathfrak{X}_G(X) \to k^\Gamma, \qquad \rho \mapsto \chi_\rho \]
 is called the \emph{$G$-character variety} and is denoted by $\chi_G(X)$. By results from \cite{culler_and_shalen}, there exists a finite set of elements $\gamma_1, \dots, \gamma_a \in \pi_1(X)$ such that $\chi_\rho$ is determined by the images $(\chi_\rho(\gamma_1), \dots , \chi_\rho(\gamma_a))$ for any $\rho$. This way, $\chi_G(X)$ can be identified with the image of the map $\mathfrak{X}_G(X) \to k^a$ given by $\rho \mapsto (\chi_\rho(\gamma_1), \dots, \chi_\rho(\gamma_a))$, which indeed provides $\chi_G(X)$ a natural structure of a variety. Again, this structure is independent of the chosen $\gamma_i$.

Note that the character map $\chi$ is a $G$-invariant morphism: indeed the trace map is invariant under conjugation. Since the affine GIT quotient is a categorical quotient, we obtain a natural morphism
\begin{equation}\label{eq:charmap}
 \overline{\chi} \colon \ \mathcal{M}_{G}(X) \to \chi_G(X) .
\end{equation}
The map $\overline{\chi}$ is an isomorphism, for instance, in the case of $G = \SL_n(\CC)$, $\Sp_{2n}
(\CC)$ or $\SO_{2n+1}(\CC)$ \cite{culler_and_shalen, flr2017, lawsik2017}.

The study of $G$-character varieties received a lot of attention. For instance, when $X$ is the underlying topological space of a smooth complex projective variety, the moduli space of $G$-representations, $\mathcal{M}_G(X)$ is one of the moduli spaces studied in non-abelian Hodge theory. In the case of a smooth complex projective curve $C$ and the algebraic group $G = \GL_n(\CC)$, the character variety $\mathcal{M}_G(C)$ parametrizes vector bundles over $C$ of rank $n$ and degree zero equipped with a flat connection. With this identification, the Riemann--Hilbert correspondence~\cite{simpson94_II} provides a real analytic isomorphism between the character variety and the moduli space of $G$-flat connections on the curve $C$. Moreover, the Hitchin--Kobayashi correspondence~\cite{simpson92} gives a real analytic isomorphism between $\mathcal{M}_G(C)$ and the moduli space of semistable $G$-Higgs bundles of rank $n$ and degree zero on $C$. These correspondences were used by Hitchin~\cite{hitchin} to compute the Poincaré polynomial of twisted character varieties for $G = \GL_2(\CC)$.

However, these correspondences are far from being algebraic. As a result, the mixed Hodge structures of the above-mentioned moduli spaces have been extensively studied, for instance, via the virtual Deligne--Hodge polynomial, or $E$-polynomial
\[ e(X) = \sum_{k,p,q} (-1)^k h_c^{k;p,q}(X)\,u^p v^q \in \ZZ[u,v] , \]
encoding the dimensions $h_c^{k;p,q}(X) = \dim_\CC H_c^{k;p,q}(X)$ of the associated graded components with respect to the weight and Hodge filtrations of the mixed Hodge structures on the compactly supported cohomology of a complex variety $X$. It follows from the excision long-exact sequence on cohomology and from the K\"unneth formula (see \cite{deligne1974theorie}) that the $E$-polynomial extends to a~motivic measure
\[ e \colon \ \K(\Var_\CC) \to \ZZ[u,v] , \]
from the Grothendieck ring of varieties to the polynomial ring in two variables $\ZZ[u,v]$, that maps the class of a complex variety $[X] \in \K(\Var_\CC)$ to its $E$-polynomial $e(X)$.

Inspired by the Weil conjectures, an arithmetic approach was introduced by Hausel and Rodríguez--Villegas \cite{hausel_villegas} to compute the $E$-polynomial of twisted $\GL_n(\CC)$-character varieties by counting its number of points over finite fields. In fact, using a theorem of Frobenius~\mbox{\cite{frobenius1896gruppencharaktere,witten1991quantum}} relating the number of points of a $G$-representation variety over a finite field $\FF_q$ to the dimensions of the irreducible representations of $G$ over $\FF_q$, Hausel and Rodríguez--Villegas show that the number of points of the twisted character variety over the finite fields $\mathbb{F}_q$ of $q$ elements is a~polynomial
 $p(q)\in \ZZ[q]$
in $q$, which in turn by a theorem of Katz, computes the $E$-polynomial of the complex twisted character variety by setting $q = uv$. This method was extended to the cases of $\SL_r(\CC)$-character varieties~\cite{mereb}, $\GL_r(\CC)$-character varieties with a generic parabolic structure~\cite{mellit2020poincare}, or non-orientable surfaces~\cite{letellier2020series}.

Recently, a geometric approach was introduced by Logares, Mu\~noz, and Newstead \cite{arXiv11066011} by dividing the representation variety into pieces and computing the $E$-polynomial piecewise. Using this technique, Martínez and Mu\~noz \cite{arXiv14076975} gave an explicit expression for the $E$-polynomial of the $\SL_2(\CC)$-representation variety. Moreover, by combining the arithmetic and geometric approaches, Baraglia and Hekmati \cite{baraglia_hekmati} gave explicit expressions for the cases $G = \GL_3(\CC), \SL_3(\CC)$.

The recursive patterns in these computations lead Gonz\'alez-Prieto, Logares, and Mu\~noz~\cite{arXiv170905724} to develop a new method using topological quantum field theories (TQFTs) to compute the virtual class of the representation varieties in the Grothendieck ring of varieties $\K(\Var_\CC)$. TQFTs, originated from physics, were first introduced by Witten \cite{witten} and axiomatized by Atiyah~\cite{atiyah}: a~TQFT is given by a monoidal functor
\[Z \colon \ \Bd_n \to R\text{-}\Mod\]
from the category of bordisms to the category of $R$-modules for some commutative ring $R$. In particular, any closed manifold $X$ can be seen as a bordism $X \colon \varnothing \to \varnothing$, so we obtain an $R$-module map $Z(X)\colon R \to R$, since $Z(\varnothing) = R$ by monoidality. As a consequence, any closed manifod $X$ has an associated invariant $Z(X)(1) \in R$. In~\cite{arXiv170905724}, Gonz\'alez-Prieto, Logares and Mu\~noz construct a lax monoidal TQFT with $R = \K(\Var_\CC)$ such that the associated invariant for a closed manifold $X$ is the virtual class of the representation variety $\mathfrak{X}_G(X)$. Then, a closed surface $X = \Sigma_g$ of genus $g$ can be considered as a composition of bordisms
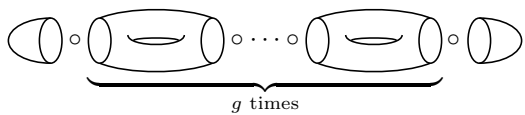
\begin{equation*}
 \renewcommand{\bdscale}{0.75}
 \Sigma_g = \bdcupright \circ \underbrace{\bdgenus \circ \cdots \circ \bdgenus}_{g \text{ times}} \circ \; \bdcupleft
\end{equation*}
so that computing the TQFT for these smaller bordisms will yield the virtual class of $\mathfrak{X}_G(\Sigma_g)$ for all $g$. This method was used in \cite{arXiv181009714} and \cite{arXiv190605222} to compute the virtual class of the (parabolic) $\SL_2(\CC)$-character variety in $\K(\Var_\CC)$. A clear advantage of this method is that it not only computes the $E$-polynomial of the representation variety, but more generally, its class in the Grothendieck ring of varieties.

\subsection{Main results}

{\bf Virtual classes of representation varieties.}
In this paper, we focus on the groups of complex upper triangular matrices $\U_n$ of rank $n = 2, 3, 4$. We compute the virtual classes of the corresponding representation varieties in a~suitably localized Grothendieck ring of varieties (see Theorems \ref{thm:result_u2}, \ref{thm:result_U3} and~\ref{thm:result_U4}).

\begin{Theorem} \label{thm:main1}
 Let $q = \big[\AA^1_\CC\big]$ be the class of the affine line in the Grothendieck ring of varieties over $\CC$. Then
 \begin{enumerate}\itemsep=0pt
 \item[$(i)$] the virtual class of the $\U_2$-representation variety $\mathfrak{X}_{\U_2}(\Sigma_g)$ is
 \[ [ \mathfrak{X}_{\U_2}(\Sigma_g) ] = q^{2g - 1} (q - 1)^{2g + 1} \big((q - 1)^{2g - 1} + 1\big) , \]
 \item[$(ii)$] the virtual class of the $\U_3$-representation variety $\mathfrak{X}_{\U_3}(\Sigma_g)$ is
 \begin{gather*}
 [ \mathfrak{X}_{\U_3}(\Sigma_g) ] = q^{3g - 3} (q - 1)^{2g} \big( q^2 (q - 1)^{2g + 1} + q^{3g} (q - 1)^2 \\
 \hphantom{[ \mathfrak{X}_{\U_3}(\Sigma_g) ] =}{}
 + q^{3g} (q - 1)^{4g} + 2 q^{3g} (q - 1)^{2g + 1} \big) ,
 \end{gather*}
 \item[$(iii)$] the virtual class of the $\U_4$-representation variety $\mathfrak{X}_{\U_4}(\Sigma_g)$ is
 \begin{gather*}
 [ \mathfrak{X}_{\U_4}(\Sigma_g) ] = q^{12 g - 6} (q - 1)^{8 g} + q^{12 g - 6} (q - 1)^{2 g + 3} + q^{10 g - 4} (q - 1)^{2 g + 3} \\
 \hphantom{[ \mathfrak{X}_{\U_4}(\Sigma_g) ] =}{}
 + q^{10 g - 3} (q - 1)^{4 g + 1} + q^{8 g - 2} (q - 1)^{6 g + 1} + q^{8 g - 2} (q - 1)^{4 g + 2} \\
 \hphantom{[ \mathfrak{X}_{\U_4}(\Sigma_g) ] =}{}
 + 2 q^{10 g - 4} (q - 1)^{6 g + 1} + 3 q^{12 g - 6} (q - 1)^{6 g + 1} + 3 q^{12 g - 6} (q - 1)^{4 g + 2}\\
 \hphantom{[ \mathfrak{X}_{\U_4}(\Sigma_g) ] =}{}
 + q^{10 g - 4} (q - 1)^{4 g + 1} (2 q - 3) .
 \end{gather*}
 \end{enumerate}
\end{Theorem}

By setting $q=uv$ in the above formulae, we obtain the $E$-polynomials of the representation varieties $\mathfrak{X}_{\U_2}(\Sigma_g)$, $\mathfrak{X}_{\U_3}(\Sigma_g)$ and $\mathfrak{X}_{\U_4}(\Sigma_g)$.

We remark that in an independent work \cite{arXiv200501841}, Gonz\'alez-Prieto, Logares, and Mu\~noz computed the virtual class of the $\AGL_1$-representation varieties, where $\AGL_1$ is the general affine group of the line. Their result can be deduced from our result on $\mathfrak{X}_{\U_2}(\Sigma_g)$ through the isomorphism of groups $\GG_m \times \AGL_1 \xrightarrow{\sim} \U_2$ (see Remark \ref{rem:result_genus_U2}).

We also provide a more general result with parabolic structures involved in the case of $\U_2$, see Theorem \ref{thm:result_U2_parabolic} of the paper.
\begin{Theorem}\label{thm:mainpar}
 Let $\Sigma_g$ be a compact oriented surface of genus $g$, with parabolic data
 \[Q = \{ (S_1, \mathcal{J}_{\lambda_1}), \dots, (S_k, \mathcal{J}_{\lambda_k}),(S_{k+1}, \mathcal{M}_{\mu_1, \sigma_1}), \dots, (S_{k+\ell}, \mathcal{M}_{\mu_\ell, \sigma_\ell}) \}.\]
 \begin{enumerate}\itemsep=0pt
 \item[$(i)$] If $\prod_{i = 1}^{k} \lambda_i \prod_{j = 1}^{\ell} \mu_j \ne 1$ or $\prod_{i = 1}^{k} \lambda_i \prod_{j = 1}^{\ell} \sigma_j \ne 1$, then
 \[ [ \mathfrak{X}_{\U_2}(\Sigma_g, Q) ] = 0 . \]

 \item[$(ii)$] Otherwise, and if $\ell = 0$, then
 \[ [ \mathfrak{X}_{\U_2}(\Sigma_g, Q) ] = q^{2g - 1} (q - 1)^{2g} \big( (-1)^k (q - 1) + (q - 1)^{2g + k} \big) , \]

 \item[$(iii)$] and if $\ell > 0$, then
 \[ [ \mathfrak{X}_{\U_2}(\Sigma_g, Q) ] = q^{2g + \ell - 1} (q - 1)^{4g + k} . \]
 \end{enumerate}
\end{Theorem}

In particular, we obtain that the representation varieties for $\U_2$, $\U_3$, and $\U_4$, and the representation varieties with parabolic structures for $\U_2$ are of balanced type, meaning that their virtual classes are generated by the virtual class of the affine line. This can be seen in relation to Higman's conjecture \cite{higman1960enumerating}, which states that the number $C_n(q)$ of conjugacy classes of $\U_n(\FF_q)$ over finite fields $\FF_q$ is polynomial in $q$, for all $n$. Indeed, using Burnside's lemma, one can relate the number of $\FF_q$-points of the $\U_n$-representation variety of a surface of genus 1 to the number of conjugacy classes:
\[ |\mathfrak{X}_{\U_n(\FF_q)}(\Sigma_1)| = |\U_n(\FF_q)| C_n(q).\]
Now, Higman's conjecture along with Katz's theorem imply that the $E$-polynomial of $\mathfrak{X}_{\U_n(\CC)}(\Sigma_1)$ is a polynomial in $uv$. We ask whether a motivic version of Higman's conjecture holds.

\begin{Conjecture}
 Let $\Sigma_g$ be a closed oriented surface of genus $g$. The class $[\mathfrak{X}_{\U_n(\CC)}(\Sigma_g)]$ of the $\U_n$-representation varieties in the Grothendieck ring of varieties is of balanced type for all~$n$ and~$g$.
\end{Conjecture}

\textbf{The moduli space of $\boldsymbol{G}$-representations and $\boldsymbol{G}$-character varieties.}
In this paper, we study the moduli space of $G$-representations of compact oriented surfaces~$\Sigma_g$ for the linear groups $\U_n$ with $n \ge 2$. As these groups are non-reductive, there is no guarantee, a priori, that the categorical quotient $\mathcal{M}_{\U_n}(\Sigma_g)$ is of finite type over $\CC$.
Nevertheless, we show the following result.
\begin{Theorem} \label{thm:moduli_space_Un_Sg}
 For all $n \ge 1$ and $g \ge 0$, there exists an isomorphism of varieties
 \[ \mathcal{M}_{\U_n}(\Sigma_g) \cong \big(\AA^1_\CC \setminus \{ 0 \}\big)^{2ng} . \]
\end{Theorem}

Furthermore, we show that the map $\overline{\chi}$ of~\eqref{eq:charmap} fails to be an isomorphism.
\begin{Theorem} \label{thm:main2}
 For $n \ge 2$ and $g \ge 1$, the natural morphism
 \[ \overline{\chi} \colon \ \mathcal{M}_{\U_n}(\Sigma_g) \to \chi_{\U_n}(\Sigma_g)\]
 is not an isomorphism.
\end{Theorem}

The paper is organized as follows. In Section \ref{sec:tqft}, we construct a topological quantum field theory (TQFT) computing the virtual classes of the representation varieties $\mathfrak{X}_G(\Sigma_g)$ in the Grothendieck ring of varieties. We mainly follow \cite{thesisangel, arXiv190605222, arXiv170905724, thesisvogel}.
The novelty of this section lies in Proposition~\ref{prop:conditions_for_reduction_TQFT}, which can be used to `reduce' the TQFT. While not strictly necessary to perform computations, it does allow for a simplification of the computations.
In Section~\ref{sec:app}, we will apply this theory to the groups of upper triangular matrices $G = \U_n$ of rank $n = 2, 3, 4$. We prove our main theorems, Theorems~\ref{thm:main1} and~\ref{thm:main2}. Explicitly, we compute the classes of the representation varieties $\mathfrak{X}_{G}(\Sigma_g)$ in a suitably localized Grothendieck ring of varieties for the groups of upper triangular matrices $\U_2$, $\U_3$, $\U_4$. Moreover, we show that the natural morphism \[\overline{\chi}\colon \ \mathcal{M}_{\U_n}(\Sigma_g)\to \chi_{\U_n}(\Sigma_g)\]
fails to be an isomorphism for $\U_n$ when $n\geq 2$ and $g \ge 1$.

\section{TQFTs and representation varieties}\label{sec:tqft}
In this section, we follow \cite{thesisangel, arXiv190605222, arXiv170905724, thesisvogel} to construct a topological quantum field theory (TQFT) that computes the virtual classes of the representation varieties $\mathfrak{X}_G(\Sigma_g)$ in the Grothendieck ring of varieties $\K(\Var_k)$.
More precisely, we construct a lax monoidal TQFT $Z$ over the ring $\K(\Var_k)$ such that the invariant associated to a closed manifold $X$ is $Z(X)(1) = [ \mathfrak{X}_G(X) ] \in \K(\Var_k)$. This construction allows to solve a more general problem: if $\Lambda$ is the set of conjugacy-closed subsets of $G$, one can put a parabolic structure $Q = \{ (S_1, \mathcal{E}_1), \dots, (S_s, \mathcal{E}_s) \}$ with data in~$\Lambda$ on~$\Sigma_g$, such that the invariant associated to $(\Sigma_g, Q)$ is the class of the variety
\begin{gather*}
 \mathfrak{X}_G(\Sigma_g, Q) = \bigg\{ (A_1, B_1, \dots, A_g, B_g, C_1, \dots, C_s) \in G^{2g + s} \,\bigg|\\
 \hphantom{\mathfrak{X}_G(\Sigma_g, Q) = \bigg\{}{}
 \prod_{i = 1}^{g} [A_i, B_i] \prod_{i = 1}^{s} C_i = 1 \text{ and } C_i \in \mathcal{E}_i \bigg\} .
\end{gather*}

The novelty of this section is Proposition \ref{prop:conditions_for_reduction_TQFT}, which allows us to modify the TQFT and simplify the computations. The rest is added for the sake of completeness. We begin by defining the categories involved in the construction of the TQFT.

\subsection{The 2-category of bordisms}

Let $i\colon M \to \partial W$ be an inclusion, where $W$ is an $n$-dimensional oriented manifold with boundary, and $M$ an $(n - 1)$-dimensional closed oriented manifold. (All manifolds we consider are assumed to be smooth.) Take a point $x \in i(M)$, let $\{ v_1, \dots, v_{n - 1} \}$ be a positively oriented basis for~$T_x i(M)$ with respect to the orientation induced by $M$, and pick some $w \in T_x i(M)$ that points inwards compared to~$W$. Then if $\{ v_1, \dots, v_{n - 1}, w \}$ is a positively oriented basis for $T_x W$, we say~$x$ is an \emph{in-boundary point}, and an \emph{out-boundary point} otherwise. Note that this is independent of the chosen vectors~$v_i$ and~$w$. If all $x \in i(M)$ are in-boundary (resp.\ out-boundary) points, we say $i$ is an \emph{in-boundary} (resp.\ \emph{out-boundary}).

\begin{Definition}
 Given two $(n - 1)$-dimensional closed oriented manifolds $M$ and $M'$, a \emph{bordism} from $M$ to $M'$ is an $n$-dimensional oriented manifold $W$ (with boundary) with maps
 \[ \begin{tikzcd} M' \arrow{r}{i'} & W & \arrow[swap]{l}{i} M, \end{tikzcd} \]
 where $i$ is an in-boundary, $i'$ an out-boundary and $\partial W = i(M) \sqcup i'(M')$. Two such bordisms~$W$,~$W'$ are said to be \emph{equivalent} if there exists an orientation-preserving diffeomorphism $W \xrightarrow{\sim} W'$ such that
 \[ \begin{tikzcd}[row sep=0.5em] & W \arrow{dd}{\wr} & \\ M' \arrow{ur} \arrow{dr} & & M \arrow{ul} \arrow{dl} \\ & W' & \end{tikzcd} \]
 commutes.
\end{Definition}
For a more precise definition of bordisms, see \cite{milnor} or \cite{kock}.

Suppose we have bordisms $W\colon M \to M'$ and $W'\colon M' \to M''$. One can glue~$W$ and~$W'$ as topological spaces by identifying the images of~$M'$, which we denote by $W \sqcup_{M'} W'$. By \cite[Theorem~1.4]{milnor}, there exists a smooth manifold structure on $W \sqcup_{M'} W'$ such that the inclusions $W \to W \sqcup_{M'} W'$ and $W' \to W \sqcup_{M'} W'$ are diffeomorphisms onto their images, which is unique up to (non-unique) diffeomorphism. Hence $W \sqcup_{M'} W'$ belongs to a well-defined equivalence class, and moreover this class only depends on the classes of~$W$ and~$W'$. Namely, if $\tilde{W}\colon M \to M'$ and $\tilde{W}'\colon M' \to M''$ are equivalent to $W$ and $W'$, respectively, then any such manifold structure on $W \sqcup_{M'} W'$ induces such a manifold structure on $\tilde{W} \sqcup_{M'} \tilde{W}'$ via the homeomorphism $W \sqcup_{M'} W' \to \tilde{W} \sqcup_{M'} \tilde{W}'$, showing $W \sqcup_{M'} W'$ and $\tilde{W} \sqcup_{M'} \tilde{W}'$ are equivalent. This implies that equivalence classes of bordisms can be composed to obtain an equivalence class of bordisms $M \to M''$.

 The discussion above gives rise to the following definition.

\begin{Definition}
 The \emph{category of $n$-bordisms}, denoted $\Bd_n$, is defined as follows. Its objects are $(n - 1)$-dimensional closed oriented manifolds, and morphisms $M \to M'$ are equivalence classes of bordisms from $M$ to $M'$. Composition is given by gluing along the common boundary: if $W\colon M \to M'$ and $W'\colon M' \to M''$, then $W' \circ W = W \sqcup_{M'} W'\colon M \to M''$.
\end{Definition}

\begin{Definition}
 The \emph{category of pointed $n$-bordisms}, denoted $\Bdp_n$, is the $2$-category consisting of:
 \begin{itemize}\itemsep=0pt
 \item Objects: pairs $(M, A)$ with $M$ an $(n - 1)$-dimensional closed oriented manifold, and $A \subset M$ a finite set of points intersecting each connected component of $M$.

 \item $1$-morphisms: a map $(M_1, A_1) \to (M_2, A_2)$ is given by a class of pairs $(W, A)$ with $W\colon M_1 \to M_2$ a bordism, and $A \subset W$ a finite set intersecting each connected component of $W$ such that $A \cap M_1 = A_1$ and $A \cap M_2 = A_2$. Two such pairs $(W, A)$ and $(W', A')$ are in the same class if there is a diffeomorphism $F\colon W \to W'$ such that $F(A) = A'$ and such that the diagram
 \begin{equation} \label{eq:diagram_from_definition_bdp} \begin{tikzcd}[row sep=0.5em] & W \arrow{dd}{\wr} & \\ M_2 \arrow{ur} \arrow{dr} & & M_1 \arrow{ul} \arrow{dl} \\ & W' & \end{tikzcd} \end{equation}
 commutes.

 The composition of $(W, A)\colon (M_1, A_1) \to (M_2, A_2)$ and $(W', A')\colon (M_2, A_2) \to (M_3, A_3)$ is $(W \sqcup_{M_2} W', A \cup A')\colon (M_1, A_1) \to (M_3, A_3)$.

 \item $2$-morphisms: a map $(W, A) \to (W', A')$ is given by a diffeomorphism $F\colon W \to W'$ such that $F(A) \subset A'$ and~\eqref{eq:diagram_from_definition_bdp} commutes.
 \end{itemize}
 Note that so far, no identity morphism exists for $(M, A)$, unless $M = A = \varnothing$. For this reason, we loosen the definition of a bordism a bit, and allow $M$ itself to be seen as a bordism $M \to M$, so that $(M, A)$ will be the identity morphism for $(M, A)$.
\end{Definition}

In this paper, we also consider manifolds which carry a so-called \textit{parabolic structure}. Fix a~set $\Lambda$ and call it the \emph{parabolic data}. We say a \emph{parabolic structure} on a manifold $M$ is a finite set $Q = \{ (S_1, \mathcal{E}_1), \dots, (S_s, \mathcal{E}_s) \}$ with $\mathcal{E}_i \in \Lambda$ and the $S_i$ are pairwise disjoint compact submanifolds of $M$ of codimension 2 with a co-orientation (i.e. an orientation of its normal bundle) such that $S_i \cap \partial M = \partial S_i$ transversally.
\begin{Definition}
 Let $\Lambda$ be a set. The \emph{$2$-category of pointed $n$-bordisms with parabolic structures over $\Lambda$}, denoted $\Bdp_n(\Lambda)$, is the $2$-category consisting of:
 \begin{itemize}\itemsep=0pt
 \item Objects: triples $(M, A, Q)$ with $M$ an $(n - 1)$-dimensional closed oriented manifold, $Q$~a~pa\-ra\-bolic structure on $M$, and $A \subset M$ a finite set of points not intersecting any of the~$S_i$ of~$Q$.

 \item $1$-morphisms: a map $(M_1, A_1, Q_1) \to (M_2, A_2, Q_2)$ is given by a class of triples $(W, A, Q)$ where $W\colon M_1 \to M_2$ is a bordism, $Q$ a parabolic structure on $W$, and $A \subset W$ a finite set intersecting each connected component of $W$ but not intersecting any of the $S_i$ of $Q$, such that $A \cap M_1 = A_1$, $A \cap M_2 = A_2$, $Q|_{M_1} = Q_1$ and $Q|_{M_2} = Q_2$. Here we use the notation $Q|_{M_i} = \{ (S_j \cap M_i, \mathcal{E}_j) \mid (S_j, \mathcal{E}_j) \in \Lambda \}$. Two such triples $(W, A, Q)$ and $(W', A', Q')$ are in the same class if there is a diffeomorphism $F\colon W \to W'$ such that $F(A) = A'$ and $(S, \mathcal{E}) \in Q$ if and only if $(F(S), \mathcal{E}) \in Q'$ and such that the diagram
 \begin{equation} \label{eq:diagram_from_definition_bdp_lambda} \begin{tikzcd}[row sep=0.5em] & W \arrow{dd}{\wr} & \\ M_2 \arrow{ur} \arrow{dr} & & M_1 \arrow{ul} \arrow{dl} \\ & W' & \end{tikzcd} \end{equation}
 commutes.

 The composition of bordisms $(W, A, Q)\colon (M_1, A_1, Q_1) \to (M_2, A_2, Q_2)$ and $(W', A', Q')\colon \allowbreak (M_2, A_2, Q_2) \to (M_3, A_3, Q_3)$ is given by $(W, A, Q) \circ (W', A', Q') = (W \sqcup_{M_2} W', A \cup A', Q \sqcup_{M_2} Q')$, where $Q \sqcup_{M_2} Q'$ denotes the union of $Q$ and $Q'$ but where we glue pairs $(S, \mathcal{E}) \in Q$ and $(S', \mathcal{E}) \in Q'$ that have a common boundary (in $M_2$).

 \item $2$-morphisms: a map $(W, A, Q) \to (W', A', Q')$ is given by a diffeomorphism $F\colon W \to W'$ such that $F(A) \subset A'$ and $(F(S), \mathcal{E}) \in Q'$ for each $(S, \mathcal{E}) \in Q$ and such that~\eqref{eq:diagram_from_definition_bdp_lambda} commutes.
 \end{itemize}
\end{Definition}

Actually, $\Bdp_n$ can be seen as a particular case of $\Bdp_n(\Lambda)$ for $\Lambda = \varnothing$. The category $\Bdp_n(\Lambda)$ (and thus $\Bdp_n$ as well) is a monoidal category. The tensor product is given by taking disjoint unions:
\[ (M, A, Q) \sqcup (M', A', Q') = (M \sqcup M', A \cup A', Q \cup Q') \]
for objects, and similarly for bordisms. The unital object is $(\varnothing, \varnothing, \varnothing)$, which we also denote simpy by $\varnothing$.

Note that (non-empty) parabolic structures can only exist on manifolds of dimension $\ge 2$. In particular for $\Bdp_2(\Lambda)$, its $1$-dimensional objects have $Q = \varnothing$ and the parabolic structures of its $2$-dimensional bordisms are of the form $\{ (p_1, \mathcal{E}_1), \dots, (p_s, \mathcal{E}_s) \}$ with $p_i$ points on the interior of the bordism that have a preferred orientation of small loops around them.

\subsection{The Grothendieck ring of varieties}

\begin{Definition}
 Let $S$ be a variety over a field $k$ (i.e., a reduced separated scheme of finite type over $k$). The \emph{Grothendieck ring of varieties} over $S$, denoted $\K(\Var/S)$, is defined as the quotient of the free abelian group on the set of isomorphism classes of varieties over $S$, by relations of the form
 \[ [X] = [X \backslash Z] + [Z], \]
 where $Z$ is a closed subvariety of $X$ and $X\backslash Z$ is its open complement. Multiplication is distributively induced by
 \[ [X] \cdot [Y] = [(X \times_S Y)_{\rm red} ] , \]
 which is indeed associative and commutative. It follows that $[ \varnothing ] = 0$ and $[ S ] = 1$ in $\K(\Var /S)$. When $S$ is the base field $k$ we denote the Grothendieck ring of varieties by $\K(\Var_k)$. To distinguish between the classes of different rings, we will write $[X]_S$ for the class of $X$ in $\K(\Var/S)$ and for the class of $X$ in $\K(\Var_k)$ we will simply write $[X]$.
\end{Definition}

Notice that $\K(\Var/S)$ is a monoid object in the category of $\K(\Var_k)$-modules. Indeed, $\K(\Var/S)$ has a natural $\K(\Var_k)$-module structure induced by
\[ [T] \cdot [X]_S = [T \times_k X]_S \]
for $T$ a variety over $k$ and $X$ a variety over $S$, such that the multiplication map $[X]_S \cdot [Y]_S = [(X\times_S Y)_\text{red}]_S$ is $\K(\Var_k)$-bilinear.

Now, we describe module maps $\K(\Var/X) \to \K(\Var/Y)$ which will be used in defining the TQFT. Let $f\colon X \to Y$ be a morphism of varieties over $k$. Composition with $f$ yields a functor
\begin{align*}
 f_! \colon \ \Var/X &\to \Var/Y \\
 \big(V \overset{g}{\longrightarrow} X\big) &\mapsto \big(V \overset{fg}{\longrightarrow} Y\big) .
\end{align*}
As $f_!$ sends isomorphisms to isomorphisms, and also
$T \times_k f_! V = f_!(T \times_k V)$ for any variety $T$ over $k$, we have that $f_!$ induces a $\K(\Var_k)$-module morphism
\[ f_! \colon \ \K(\Var/X) \to \K(\Var/Y) . \]
Note that this map will in general not be a ring morphism. For example, the unit $[X]_X \in \K(\Var/X)$ need not be sent to the unit $[Y]_Y \in \K(\Var/Y)$.

Similarly, pulling back along $f$ yields a functor
\[ f^* \colon \ \Var/Y \to \Var/X \]
sending $W \overset{h}{\longrightarrow} Y$ to $W \times_Y X \overset{f^*h}{\longrightarrow} X$. Also $f^*$ induces a map,
\[ f^* \colon \ \K(\Var/Y) \to \K(\Var/X) , \]
 which is a $\K(\Var_k)$-module morphism as $T \times_k f^*(V) = T \times_k (V \times_Y X) \cong (T \times_k V) \times_Y X = f^*(T \times_k V)$ for any variety $T$ over $k$. In contrast to $f_!$, the map $f^*$ is a ring morphism as $(V \times_Y W) \times_Y X \cong (V \times_Y X) \times_X (W \times_Y X)$ for any $V$, $W$ over $Y$.

\begin{Example} \label{ex:2morp}
 The rings $\K(\Var/X)$ are objects of the category of $\K(\Var_k)\text{-}\Mod$. Moreover, given a variety $Z$ over $k$ with morphisms $f \colon Z \to X$ and $g \colon Z \to Y$, we have an induced $\K(\Var_k)$-module morphism $g_!\circ f^* \colon \K(\Var/X) \to \K(\Var/Y)$.
 \[ \begin{tikzcd}[row sep=0.5em] \K(\Var/X) \arrow{r}{f^*} & \K(\Var/Z) \arrow{r}{g_!} & \K(\Var/Y). \end{tikzcd} \]
\end{Example}

\begin{Remark}
 The functors $f^*$ and $f_!$ are adjoint, as for any varieties $V \overset{v}{\longrightarrow} X$ and $W \overset{w}{\longrightarrow} Y$ there is a bijection
 \[ \Hom_{\Var/Y}(f_! V, W) \cong \Hom_{\Var/X}(V, W \times_Y X) \]
 natural in $V$ and $W$. Namely, by the universal property of the fiber product, to give a morphism $\varphi\colon V \to W \times_Y X$ is to give morphisms $V \overset{r}{\longrightarrow} W$ and $V \overset{s}{\longrightarrow} X$ such that $w \circ r = f \circ s$, and requiring $\varphi$ to be over $X$ means to have $s = v$. Hence, to give $\varphi$ over $X$ is to give $V \overset{r}{\longrightarrow} W$ such that $w \circ r = v$, i.e. a morphism $V \overset{r}{\longrightarrow} W$ over $Y$. The naturality of this bijection is easily seen.
\end{Remark}

In this paper, the target category of the TQFT is the 2-category of $\K(\Var_k)$-modules with twists.

\begin{Definition}
 Let $R$ be a commutative ring. Given $R$-module morphisms $f, g\colon M \to N$, we say $g$ is an \emph{immediate twist} of $f$ if there exists an $R$-module $P$ and $R$-module morphisms $f_1\colon M \to P$, $f_2\colon P \to N$ and $h\colon P \to P$ such that $f = f_2 \circ f_1$ and $g = f_2 \circ h \circ f_1$:
 \[ \begin{tikzcd} M \arrow{r}{f_1} \arrow[swap, bend right=20]{rr}{g} & \arrow[out=120,in=60,loop, "h"] P \arrow{r}{f_2} & N. \end{tikzcd} \]
 We say a \emph{twist} from $f$ to $g$ is a finite sequence $f = f_0, f_1, \dots, f_n = g\colon M \to N$ of $R$-module morphisms such that $f_{i + 1}$ is an immediate twist of $f_i$.

 Now the \emph{$2$-category of $R$-modules with twists}, denoted $R\text{-}\Mod_t$, is the category whose objects are $R$-modules, its 1-morphisms are $R$-module morphisms, and its 2-morphisms are twists.
\end{Definition}

\subsection{The 2-category of spans}
\begin{Definition}
 Given a category $\mathcal{C}$ with pullbacks, we define the 2-category $\Span(\mathcal{C})$ as follows:
 \begin{itemize}\itemsep=0pt
 \item Its objects are the objects of $\mathcal{C}$.
 \item An arrow from $C$ to $C'$ is given by a diagram $C \leftarrow D \rightarrow C'$ in $\mathcal{C}$, called a \textit{span}. Composition of the spans $C \leftarrow D \rightarrow C'$ and $C' \leftarrow D' \rightarrow C''$ is given by the span $C \leftarrow E \rightarrow C''$ such that the square in
 \[ \begin{tikzcd}[row sep=0em] & & E \arrow{ld} \arrow{rd} & & \\ & D \arrow{ld} \arrow{rd} & & D' \arrow{ld} \arrow{rd} & \\ C & & C' & & C'' \end{tikzcd} \]
 is a pullback square.
 \item A 2-morphism from $C \leftarrow D \rightarrow C'$ to $C \leftarrow D' \rightarrow C'$ is given by an arrow $D \to D'$ such that the following diagram commutes:
 \[ \begin{tikzcd}[row sep = 0em] & D \arrow{dd} \arrow{ld} \arrow{rd} & \\ C & & C'. \\ & D' \arrow{lu} \arrow{ru} & \end{tikzcd} \]
 \end{itemize}
\end{Definition}

The category $\Span^\op(\mathcal{C})$ is defined analogously on categories with pushouts, where we reverse all arrows.

If $\mathcal{C}$ is a monoidal category, then $\Span(\mathcal{C})$ naturally has the structure of a monoidal category as well. The tensor product and unital object naturally carry over, the associator will be given by the span
\[ \begin{tikzcd} A \otimes (B \otimes C) & \arrow[swap]{l}{\id} A \otimes (B \otimes C) \arrow{r}{\alpha_{A, B, C}} & (A \otimes B) \otimes C \end{tikzcd} \]
and the left and right unitor by
\[ \begin{tikzcd} I \otimes A & \arrow[swap]{l}{\id} I \otimes A \arrow{r}{\lambda_A} & A \end{tikzcd} \qquad \text{and} \qquad \begin{tikzcd} A \otimes I & \arrow[swap]{l}{\id} A \otimes I \arrow{r}{\rho_A} & A. \end{tikzcd} \]

We are ready to construct the TQFT computing the classes of the representation varieties $\mathfrak{X}_G(\Sigma_g)$ in $\K(\Var_k)$.

\subsection{Constructing the TQFT}

\begin{Definition}
 Let $(X, A)$ be a pair of topological spaces. The \emph{fundamental groupoid of $X$ w.r.t.~$A$} denoted $\Pi (X, A)$ is the groupoid category whose objects are elements of $A$, and an arrow $a \to b$ for each homotopy class of paths from $a$ to $b$. Composition of morphisms is given by concatenation of paths. Indeed this construction only depends on the homotopy type of~$(X, A)$. In particular, if $A = \{ x_0 \}$ is a single point, we obtain the fundamental group $\pi_1(X, x_0)$.
\end{Definition}

Note that if $f\colon (X, A) \to (X', A')$ is a map of pairs of topological spaces, there is an induced functor $\Pi(X, A) \to \Pi(X', A')$ between groupoids, mapping an object $a \in A$ to $f(a) \in A'$ and an arrow $\gamma\colon a \to b$ to $f \circ \gamma$. This allows us to construct the following functor.

\begin{Definition} \label{def:geometrization_functor}
 The \emph{geometrization functor} is a $2$-functor $\Pi\colon \Bdp_n \to \Span^\op(\Grpd)$, with $\Grpd$ the category of groupoids, defined as follows:
 \begin{itemize}\itemsep=0pt
 \item To each object $(X, A)$ we assign the fundamental groupoid $\Pi(X, A)$.
 \item For each $1$-morphism $(W, A)\colon (X_1, A_1) \to (X_2, A_2)$ we assign the cospan
 \[ \Pi(X_1, A_1) \overset{i_1}{\longrightarrow} \Pi(W, A) \overset{i_2}{\longleftarrow} \Pi(X_2, A_2) \]
 with $i_1$ and $i_2$ induced by inclusions.
 \item For each $2$-morphism $(W, A) \to (W, A')$ given by diffeomorphism $F\colon W \to W'$ with ${F(A) \subset A'}$, we obtain a groupoid morphism $\Pi F$ yielding the commutative diagram
 \[ \begin{tikzcd}[row sep=1em] & \Pi(W, A) \arrow{dd}{\Pi F} & \\ \Pi(X_1, A_1) \arrow{ur}{i_1} \arrow{dr}{i'_1} & & \Pi(X_2, A_2), \arrow[swap]{ul}{i_2} \arrow[swap]{dl}{i'_2} \\ & \Pi(W', A') & \end{tikzcd} \]
 which is a $2$-morphism in $\Span^\op(\Grpd)$.
 \end{itemize}
\end{Definition}

The Seifert--van~Kampen theorem for fundamental groupoids \cite{brown} provides that $\Pi$ defined above is indeed a functor (for more details see \cite{thesisangel}).

Suppose $X$ is a compact connected manifold (possibly with boundary), $A \subset X$ a finite set of points, and denote $\mathcal{G} = \Pi(X, A)$. We write $\mathcal{G}_a = \Hom_{\mathcal{G}}(a, a)$ for $a$ in $\mathcal{G}$. Since a compact connected manifold has the homotopy type of a~finite CW-complex, every $\mathcal{G}_a = \pi_1(X, a)$ is a~finitely generated group.

The groupoid $\mathcal{G}$ has finitely many connected components, where we say objects $a, b \in \mathcal{G}$ are \textit{connected} if $\Hom_{\mathcal{G}}(a, b)$ is non-empty.
Pick a subset $S = \{ a_1, \dots, a_s \} \subset A$ such that each connected component of $\mathcal{G}$ contains exactly one of the $a_i$. Also pick an arrow $f_a\colon a_i \to a$ for each $a \in A$ (with $a_i \in S$ in the connected component of $a$) such that $f_{a_i} = \id_{a_i}$ for each $a_i \in S$. Now if $G$ is a group, then a morphism of groupoids $\rho\colon \mathcal{G} \to G$ is uniquely determined by the group morphisms $\rho_i\colon \mathcal{G}_{a_i} \to G$ and a choice of $\rho(f_a) \in G$. Namely, any $\gamma\colon a \to b$ in $\mathcal{G}$ can be written as $\gamma = f_b \circ \gamma' \circ (f_a)^{-1}$ for some $\gamma' \in \mathcal{G}_{a_i}$ (with $a_i \in S$ in the connected component of $a$ and $b$). The elements $\rho(f_a)$ can take any value for $a \not\in A \backslash S$ (and $\rho(f_{a_i}) = 1 \in G$ for $a_i \in S$), so if $\mathcal{G}$ has $n$ objects and $s$ connected components, we have
\begin{equation}
 \label{eq:the_hom_equation}
 \Hom_{\Grpd}(\mathcal{G}, G) \cong \Hom(\mathcal{G}_{a_1}, G) \times \cdots \times \Hom(\mathcal{G}_{a_s}, G) \times G^{n - s} .
\end{equation}
If $G$ is an algebraic group, each of these factors naturally carries the structure of an algebraic variety. Namely, each $\mathcal{G}_{a_i}$ is finitely generated, so $\Hom(\mathcal{G}_{a_i}, G)$ can be identified with a subvariety of $G^m$ for some $m > 0$. This gives $\Hom(\mathcal{G}, G)$ the structure of an algebraic variety, and this structure can be shown not to depend on the choices. 

\begin{Definition} \label{def:representation_variety}
 Let $X$ be a compact connected manifold (possibly with boundary) and $A \subset X$ a finite set. Then we define the \emph{$G$-representation variety} of $(X, A)$ to be
 \[ \mathfrak{X}_G(X, A) = \Hom_{\Grpd}(\Pi(X, A), G) . \]
\end{Definition}

Note that the functor $\Hom_{\Grpd}(-, G)$ sends pushouts to pullbacks, so we obtain an induced 2-functor
\[ \mathcal{F}\colon \ \Bdp_n \to \Span(\Var_k), \]
which we refer to as the \emph{field theory}. This functor sends an object $(M, A)$ to $\mathfrak{X}_G(M, A)$, a~bordism $(W, A)\colon (M_1, A_1) \to (M_2, A_2)$ to the span
\[ \mathfrak{X}_G(M_1, A_1) \longleftarrow \mathfrak{X}_G(W, A) \longrightarrow \mathfrak{X}_G(M_2, A_2), \]
and a 2-morphism given by diffeomorphism $F\colon (W,A) \to (W',A')$ with $F(A) \subset A'$ to an inclusion of the corresponding varieties.

Recall that $\Var_k$ is a monoidal category, the tensor product being $\times_k$ and the unital object $\Spec k$, so the category $\Span(\Var_k)$ is monoidal as well. Following the above construction, one can see that $\mathcal{F}$ is a monoidal functor. Indeed, $\mathfrak{X}_G(\varnothing) = \Hom_{\Grpd}(\varnothing, G)$ is a point, and $\mathfrak{X}_G(X \sqcup X', A \cup A') \cong \mathfrak{X}_G(X, A) \times_k \mathfrak{X}_G(X', A')$ as can be easily shown from \eqref{eq:the_hom_equation}. Also, $\mathcal{F}$ is seen to be symmetric.

The last step in constructing the TQFT is the \emph{quantization functor}
\[ \mathcal{Q}\colon \ \Span(\Var_k) \to \K(\Var_k)\text{-}\Mod_t, \]
which assigns to an object $X$ the $\K(\Var_k)$-module $\K(\Var/X)$, and to a span $X \overset{f}{\longleftarrow} Z \overset{g}{\longrightarrow} Y$ the morphism $g_! \circ f^*\colon \K(\Var/X) \to \K(\Var/Y)$ (see Example~\ref{ex:2morp}). Given a 2-morphism
\[ \begin{tikzcd}[row sep=0.5em] & Z_1 \arrow{dd}{h} \arrow[swap]{ld}{f_1} \arrow{rd}{g_1} & \\ X & & Y \\ & Z_2 \arrow{lu}{f_2} \arrow[swap]{ru}{g_2} & \end{tikzcd} \]
we see that $(g_1)_! \circ (f_1)^* = (g_2)_! \circ h_! \circ h^* \circ f_2^*$, which defines an (immediate) twist from $(g_1)_! \circ (f_1)^*$ to $(g_2)_! \circ (f_2)^*$. For a proof that $\mathcal{Q}$ is a~lax (symmetric) monoidal 2-functor, see \cite[Theorem~4.13]{arXiv181009714}.

\begin{Remark}
 Contrary to $\mathcal{F}$, the quantization functor $\mathcal{Q}: \Span(\Var_k) \to \K(\Var_k)\text{-}\Mod_t$ is not a monoidal functor. Namely, even though there is a natural map \begin{gather*}
 \mathcal{Q}(X) \otimes \mathcal{Q}(Y) = \K(\Var/X) \otimes \K(\Var/Y) \to \K(\Var/(X \times Y)) = \mathcal{Q}(X \times Y), \\
 [V \to X] \otimes [W \to Y] \mapsto [V \times W \to X \times Y] ,
 \end{gather*}
 this map need not be an isomorphism. Relaxing the condition of this natural map to be an isomorphism, we obtain the notion of a \textit{lax} monoidal functor.
\end{Remark}

At last, we define the symmetric lax monoidal TQFT as the composition of the field theory and the quantization functor
\[ Z = \mathcal{Q} \circ \mathcal{F}\colon \ \Bdp_n \to \K(\Var_k)\text{-}\Mod_t . \]

Now, any closed connected oriented manifold $X$ of dimension $n$, with a point $\star$ on $X$, can be viewed as a bordism $(X, \star)\colon \varnothing \to \varnothing$. Then $\mathcal{F}(X, \star)$ is the span
\[ \star \overset{t}{\longleftarrow} \mathfrak{X}_G(X, \star) = \mathfrak{X}_G(X) \overset{t}{\longrightarrow} \star \]
and $Z(X, \star)(1) = t_! t^*(\star) = t_!\left([ \mathfrak{X}_G(X) ]_{\mathfrak{X}_G(X)}\right) = [ \mathfrak{X}_G(X) ]$ as desired.

\subsection{Parabolic structures} \label{sec:TQFT_parabolic_structures}
Let $\Lambda$ be a set of conjugacy-closed subsets of $G$. We will slightly modify the construction above to obtain a (lax symmetric) TQFT $Z_\Lambda\colon \Bdp_n(\Lambda) \to \K(\Var_k)\text{-}\Mod_t$. This will be an extension of $Z$ in the sense that it yields the same modules and morphisms as $Z$ in the absence of parabolic structures, i.e., $Z_\Lambda(X, A, \varnothing) = Z(X, A)$.

Let $X$ be a compact manifold (possibly with boundary) with a parabolic structure $Q$ given by
\[ Q = \{ (S_1, \mathcal{E}_1), \dots, (S_s, \mathcal{E}_s) \} , \]
and $A \subset X$ a finite set intersecting each connected component of $X$, but not intersecting $S = \cup_i S_i$. Then the representation variety of $(X, A, Q)$ is defined as
\begin{equation}
 \label{eq:representation_variety_parabolic}
 \mathfrak{X}_G(X, A, Q) = \left\{ \rho\colon \Pi(X - S, A) \to G \,\bigg|\, \begin{matrix} \rho(\gamma) \in \mathcal{E}_i \text{ for all loops } \gamma \text{ around } S_i \\ \text{ positive w.r.t. the co-orientation}, \\ \text{for all } (S_i, \mathcal{E}_i) \in Q\end{matrix} \right\} ,
\end{equation}
where `$\gamma$ around $S_i$' means a non-zero loop $\gamma$ in $\Pi(X - S, A)$ which is zero in $\Pi(X - (S - S_i), A)$.
Since the $\mathcal{E}_i$ are conjugacy-closed, the condition on the loops $\gamma$ around $S_i$ is independent on the chosen base point. Indeed, this definition of $\mathfrak{X}_G(X, A, Q)$ agrees with Definition \ref{def:representation_variety} for $Q = \varnothing$. When $X$ is connected, we write
\[ \mathfrak{X}_G(X, Q) = \mathfrak{X}_G(X, \star, Q) . \]
In the particular case of $X = \Sigma_g$ with parabolic structure $Q = \{ (S_1, \mathcal{E}_1), \dots, (S_s, \mathcal{E}_s) \}$ (co-orientation induced from orientation on $\Sigma_g$), we find
\begin{gather*}
 \mathfrak{X}_G(\Sigma_g, Q) = \bigg\{ (A_1, B_1, \dots, A_g, B_g, C_1, \dots, C_s) \in G^{2g + s} \,\bigg|\\
 \hphantom{\mathfrak{X}_G(\Sigma_g, Q) = \bigg\{}
 \prod_{i = 1}^{g} [A_i, B_i] \prod_{i = 1}^{s} C_i = 1 \text{ and } C_i \in \mathcal{E}_i \bigg\} .
\end{gather*}
Consider the modified field theory $\mathcal{F}_\Lambda\colon \Bdp_n(\Lambda) \to \Span(\Var_k)$ that maps an object $(M, A, Q)$ to $\mathfrak{X}_G(M, A, Q)$, a bordism $(W, A, Q)\colon (M_1, A_1, Q_1) \to (M_2, A_2, Q_2)$ to the span
\[ \mathfrak{X}_G(M_1, A_1, Q_1) \longleftarrow \mathfrak{X}_G(W, A, Q) \longrightarrow \mathfrak{X}_G(M_2, A_2, Q_2) \]
induced by the inclusions and a 2-morphism given by diffeomorphism $F\colon W \!\to\! W'$ with ${F(A) \!\subset\! A'}$ to an inclusion of the corresponding varieties (see~\cite{thesisangel}).
It is easy to see that this functor is still monoidal. We obtain the resulting (lax symmetric) TQFT
\[ Z_\Lambda = \mathcal{Q} \circ \mathcal{F}_\Lambda\colon \ \Bdp_n(\Lambda) \to \K(\Var_k)\text{-}\Mod_t . \]
To a closed connected oriented manifold $X$ with parabolic structure $Q$ is now associated the invariant
\[ Z_\Lambda(X, A, Q)(1) = [ \mathfrak{X}_G(X, A, Q) ] . \]
Since $Z_\Lambda$ is understood to be an extension of the earlier TQFT $Z\colon \Bdp_n \to \K(\Var_k)\text{-}\Mod_t$, and since it is clear what set $\Lambda$ we consider, we will just write $Z$ for $Z_\Lambda$.

\subsection{Field theory in dimension 2}\label{sec:field}
We focus on the case of dimension $n = 2$. Let $X = \Sigma_g$ be a closed oriented $2$-dimensional surface of genus $g$, possibly with a parabolic structure $Q$. Now $\Sigma_g$ can be considered as a~bordism $\varnothing \to \varnothing$, and after taking a~suitable finite set $A \subset \Sigma_g$, be written as a composition of the following bordisms:
\begin{equation}
 \label{eq:generators_bdp_2_Lambda}
 \begin{tikzcd}[row sep=0em, column sep=0.5em] \bdpcupright & \bdpgenus & \bdpparabolic & \bdpcupleft \\
 D^\dag\colon \big(S^1, \star \big) \to \varnothing & L\colon \big(S^1, \star \big) \to (S^1, \star) & L_{\mathcal{E}}\colon \big(S^1, \star \big) \to \big(S^1, \star \big) & D\colon \varnothing \to \big(S^1, \star \big) \end{tikzcd}
\end{equation}
Here $L_\mathcal{E}$ denotes the cylinder with parabolic structure $\{ (\star, \mathcal{E}) \}$ with $\mathcal{E} \in \Lambda$. Now indeed, if we write $Q = \{ (p_1, \mathcal{E}_1), \dots, (p_s, \mathcal{E}_s) \}$ for the parabolic structure on $\Sigma_g$, we have
\begin{equation}
 \label{eq:decomposition_of_Sigma_g}
 (\Sigma_g, A, Q) = D^\dag \circ L^g \circ L_{\mathcal{E}_1} \circ \cdots \circ L_{\mathcal{E}_s} \circ D .
\end{equation}

Of course, the category $\Bdp_2(\Lambda)$ consists of more objects and morphisms than just the ones mentioned in \eqref{eq:generators_bdp_2_Lambda}. However, as we are only interested in closed connected surfaces (possibly with a parabolic structure), we will restrict our attention to a subcategory of $\Bdp_2(\Lambda)$: we say a \emph{strict tube} is any composition of the bordisms in \eqref{eq:generators_bdp_2_Lambda}, and let $\Tb_2(\Lambda)$ be the subcategory of $\Bdp_2(\Lambda)$ whose objects are disjoint copies of $\big(S^1, \star\big)$ and bordisms are disjoint unions of strict tubes. Note that $\Tb_2(\Lambda)$ is still monoidal (with the same monoidal structure as $\Bdp_2(\Lambda)$). We refer to $\Tb_2(\Lambda)$ as the \emph{category of tubes}.

We restrict $Z$ to a functor $\Tb_2(\Lambda) \to \K(\Var_k)\text{-}\Mod_t$, and explicitly describe what the TQFT does to our objects and bordisms in \eqref{eq:generators_bdp_2_Lambda}.

The fundamental groups $\pi_1(D)$ and $\pi_1\big(D^\dag\big)$ are trivial, implying $\mathfrak{X}_G(D) = \mathfrak{X}_G\big(D^\dag\big) = \star$. Since $\pi_1\big(S^1, \star\big) = \ZZ$, we have $\mathfrak{X}_G\big(S^1, \star\big) = \Hom(\ZZ, G) = G$ and since $\Pi(\varnothing)$ is the empty groupoid, we have $\mathfrak{X}_G(\varnothing) = \star$. Hence the field theory for $D$ and $D^\dag$ is given by
\[ \begin{array}{@{}cccccc}
 \mathcal{F}(D)\colon & \star & \longleftarrow & \star & \longrightarrow & G \\
 & \star & \mapsfrom & \star & \mapsto & 1
\end{array} \qquad \text{and} \qquad \begin{array}{@{}cccccc}
 \mathcal{F}\big(D^\dag\big)\colon & G & \longleftarrow & \star & \longrightarrow & \star \\
 & 1 & \mapsfrom & \star & \mapsto & \star \rlap{.}
\end{array} \]
For the bordism $L$, call its two basepoints $a$ and $b$. The surface of $L$ is homotopic to a torus with two punctures, so its fundamental group (w.r.t.~$a$) is the free group $F_3$. We pick generators $\gamma$, $\gamma_1$, $\gamma_2$ as depicted in the following image, and a path $\alpha$ connecting $a$ and $b$:
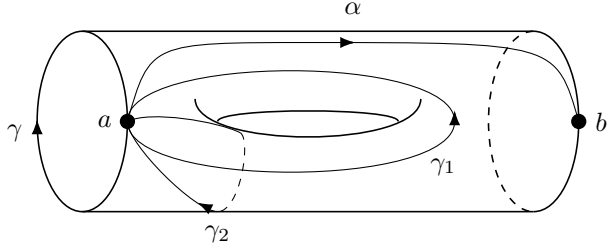
$$ \begin{tikzpicture}[semithick, scale=3.0]
 \begin{scope}
 \draw (-1,0) ellipse (0.2cm and 0.4cm);
 \draw (-1,0.4) -- (1,0.4);
 \draw (-1,-0.4) -- (1,-0.4);
 \draw (1,0.4) arc (90:-90:0.2cm and 0.4cm);
 \draw[dashed] (1,0.4) arc (90:270:0.2cm and 0.4cm);

 \draw (-0.5,0.1) .. controls (-0.5,-0.125) and (0.5,-0.125) .. (0.5,0.1);
 \draw (-0.4,0.0) .. controls (-0.4,0.0625) and (0.4,0.0625) .. (0.4,0.0);

 \draw[black,fill=black] (-0.8,0) circle (.2ex);
 \draw[black,fill=black] (1.2,0) circle (.2ex);

 \draw[thin] (-0.8,0) .. controls (-0.7,0.35) .. (0,0.35) .. controls (1.1,0.35) .. (1.2,0);
 \draw[thin] (-0.8,0) .. controls (-0.8,0.3) and (0.65,0.3) .. (0.65,0) .. controls (0.65,-0.3) and (-0.8,-0.3) .. (-0.8,0);
 \draw[thin] (-0.8,0) .. controls (-0.7,0.05) and (-0.4,0.0) .. (-0.3,-0.05);
 \draw[thin, dashed] (-0.3,-0.05) .. controls (-0.25,-0.08) and (-0.3,-0.4) .. (-0.4,-0.4);
 \draw[thin] (-0.4,-0.4) .. controls (-0.5,-0.4) and (-0.8,-0.1) .. (-0.8,0);

 \draw[-{Latex}] (-1.2,0) -- +(0,0.01);
 \draw[-{Latex}] (0.2,0.35) -- +(0.01,0);
 \draw[-{Latex}] (0.65,0) -- +(0,0.05);
 \draw[-{Latex}] (-0.47,-0.37) -- +(-0.02,0.01);

 \node at (-1.3,-0.05) {$\gamma$};
 \node at (0.6,-0.2) {$\gamma_1$};
 \node at (-0.4,-0.5) {$\gamma_2$};
 \node at (0.2,0.5) {$\alpha$};
 \node at (-0.9,0) {$a$};
 \node at (1.3,0) {$b$};
 \end{scope}
\end{tikzpicture} $$
According to \eqref{eq:the_hom_equation} we can now identify
\[ \mathfrak{X}_G(L) \cong \Hom(F_3, G) \times G \cong G^4, \qquad \rho \mapsto (\rho(\gamma), \rho(\gamma_1), \rho(\gamma_2), \rho(\alpha)) . \]
A generator for $\pi_1(S^1, b)$ is given by $\alpha \gamma [\gamma_1, \gamma_2] \alpha^{-1}$, and so the field theory for $L$ is found to be
\begin{equation}
 \label{eq:span_Z_L}
 \begin{array}{@{}cccccc}
 \mathcal{F}(L)\colon & G & \longleftarrow & G^4 & \longrightarrow & G \\
 & g & \mapsfrom & (g, g_1, g_2, h) & \mapsto & h g [g_1, g_2] h^{-1} .
 \end{array}
\end{equation}
Finally for the bordism $L_\mathcal{E}$, call its two basepoints $a$ and $b$. The fundamental group (w.r.t.~$a$) of the cylinder with a puncture is the free group $F_2$. We pick generators~$\gamma$,~$\gamma'$ as depicted in the following image, and a path~$\alpha$ connecting~$a$ and~$b$:
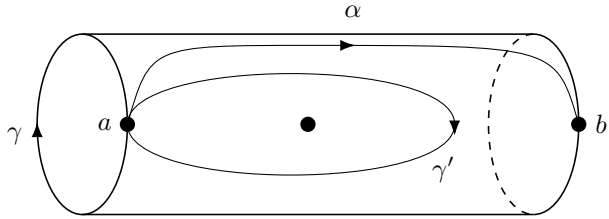
$$ \begin{tikzpicture}[semithick, scale=3.0]
 \begin{scope}
 \draw (-1,0) ellipse (0.2cm and 0.4cm);
 \draw (-1,0.4) -- (1,0.4);
 \draw (-1,-0.4) -- (1,-0.4);
 \draw (1,0.4) arc (90:-90:0.2cm and 0.4cm);
 \draw[dashed] (1,0.4) arc (90:270:0.2cm and 0.4cm);

 \draw[black,fill=black] (-0.8,0) circle (.2ex);
 \draw[black,fill=black] (1.2,0) circle (.2ex);
 \draw[black,fill=black] (0,0) circle (.2ex);

 \draw[thin] (-0.8,0) .. controls (-0.7,0.35) .. (0,0.35) .. controls (1.1,0.35) .. (1.2,0);
 \draw[thin] (-0.8,0) .. controls (-0.8,0.3) and (0.65,0.3) .. (0.65,0) .. controls (0.65,-0.3) and (-0.8,-0.3) .. (-0.8,0);

 \draw[-{Latex}] (-1.2,0) -- +(0,0.01);
 \draw[-{Latex}] (0.2,0.35) -- +(0.01,0);
 \draw[-{Latex}] (0.65,0) -- +(0,-0.05);

 \node at (-1.3,-0.05) {$\gamma$};
 \node at (0.6,-0.2) {$\gamma'$};
 \node at (0.2,0.5) {$\alpha$};
 \node at (-0.9,0) {$a$};
 \node at (1.3,0) {$b$};
 \end{scope}
\end{tikzpicture}$$
Using \eqref{eq:representation_variety_parabolic} we can now identify
\[ \mathfrak{X}_G(L_\mathcal{E}) \cong G^2 \times \mathcal{E}, \qquad \rho \mapsto (\rho(\gamma), \rho(\alpha), \rho(\gamma')) . \]
A generator for $\pi_1(S^1, b)$ is given by $\alpha \gamma \gamma' \alpha^{-1}$, and so the field theory for $L_\mathcal{E}$ is found to be
\[ \begin{array}{@{}cccccc}
 \mathcal{F}(L_\lambda)\colon & G & \longleftarrow & G^2 \times \mathcal{E} & \longrightarrow & G \\
 & g & \mapsfrom & (g, h, \xi) & \mapsto & h g \xi h^{-1} .
\end{array} \]
Finally, using \eqref{eq:the_hom_equation} and \eqref{eq:decomposition_of_Sigma_g} we can express the class of the character variety $\mathfrak{X}_G(\Sigma_g, Q)$ in terms of the TQFT (see also \cite[Theorem 4.11]{arXiv181009714}).

\begin{Theorem} Let $\Sigma_g$ be the closed oriented surface of genus $g$ with parabolic structure $Q = \{ (S_1, \mathcal{E}_1), \dots, (S_s, \mathcal{E}_s) \}$. Then,
\begin{align}
 [ \mathfrak{X}_G(\Sigma_g, Q) ]
 &= \frac{1}{[G]^{g + s}} [ \mathfrak{X}_G(\Sigma_g, \{ g + s + 1 \text{ points} \}, Q) ] \nonumber \\
 &= \frac{1}{[G]^{g + s}} Z\big(D^\dag \circ L^g \circ L_{\mathcal{E}_1} \circ \dots \circ L_{\mathcal{E}_s} \circ D\big)(1) . \label{eq:representation_variety_from_TQFT}
\end{align}
\end{Theorem}

Note that $[G]$ might be a zero-divisor. Indeed, in this paper, we focus on the groups of upper triangular matrices $\U_n$ of rank $n = 2, 3, 4$ whose classes are products of the class of the punctured affine line $\big[\AA^1_k \setminus \{ 0 \}\big]$ and the class of the affine line $\big[\AA^1_k\big]$, where the latter is a zero divisor in $\K(\Var_k)$ \cite{Borisov2014, Martin2016}. However, one can consider a suitable localization of $\K(\Var_k)$, the localization with the classes of $\big[\AA^1_k\big]$ and $\big[\AA^1_k \setminus \{ 0 \}\big]$, in which the computation described by~\eqref{eq:representation_variety_from_TQFT} holds.

\subsection{Reduction of the TQFT}\label{sec:reduction_TQFT}

Let $Z\colon \mathcal{B} \to R\text{-}\Mod_t$ be a (lax monoidal) TQFT, where $\mathcal{B}$ is some kind of bordism category (e.g., $\Bdp_n$, $\Bdp_n(\Lambda)$ or $\Tb_2(\Lambda)$) and $R$ is a commutative ring. For some TQFTs there is a symmetry, such as a group action, which can be used to `reduce' the TQFT, allowing for a~simplification of the computations. In this section, we will show how such a reduction can be obtained.

For each object $M \in \mathcal{B}$, let
\[ \begin{tikzcd} Z(M) \arrow[shift left=0.25em]{r}{\alpha_M} & \arrow[shift left=0.25em]{l}{\beta_M} N_M \end{tikzcd} \]
be $R$-module morphisms, with $N_M$ an $R$-module. Assume that $N_\varnothing = Z(\varnothing)$ with $\alpha_\varnothing$ and $\beta_\varnothing$ the identity maps. Let $V_M$ be a submodule of $N_M$ such that $(\alpha_{M'} \circ Z(W) \circ \beta_M)(V_M) \subset V_{M'}$ for all bordisms $W\colon M \to M'$ in $\mathcal{B}$. In particular, $(\alpha_M \circ \beta_M)(V_M) \subset V_M$ for any $M$.

\begin{Lemma} \label{lemma:TQFT_reduction}
 Suppose that the map $\alpha_M \circ \beta_M\colon V_M \to V_M$ is invertible for all $M \in \mathcal{B}$. Then for every $W\colon M \to M'$ there exists a unique $R$-linear map $\tilde{Z}(W)\colon V_M \to V_{M'}$ such that the following diagram commutes:
 \[ \begin{tikzcd}[row sep=3em, column sep=5em]
 \beta_M(V_M) \arrow{r}{Z(W)} \arrow[swap]{d}{\alpha_M} & Z(W)(\beta_M(V_M)) \arrow{d}{\alpha_{M'}} \\ V_M 
 \arrow[dashed,swap]{r}{\tilde{Z}(W)} & V_{M'}. 
 \end{tikzcd} \]
\end{Lemma}

\begin{proof}
 Indeed, the above diagram is well-defined by the assumption that $(\alpha_{M'} \circ Z(W) \circ \beta_M)(V_M) \subset V_{M'}$. We are looking for an $R$-linear map $\tilde{Z}(W)\colon V_M \to V_{M'}$ such that $\tilde{Z}(W) \circ \alpha_M = \alpha_{M'} \circ Z(W)$. Precomposing this equality with $\beta_M \circ (\alpha_M \circ \beta_M)^{-1}$ gives
 \begin{equation} \label{eq:expression_Z_tilde} \tilde{Z}(W) = \alpha_{M'} \circ Z(W) \circ \beta_M \circ (\alpha_M \circ \beta_M)^{-1} . \end{equation}
 This shows there is a unique choice of $\tilde{Z}(W)$, and it is easy to see that this choice makes the diagram commute: any $x \in \beta_M(V_M)$ can be written as $x = \beta_M(y)$ for some $y \in V_M$, so
 \[ \tilde{Z}(W) \circ \alpha_M(x) = \tilde{Z}(W) \circ (\alpha_M \circ \beta_M) (y) = \alpha_{M'} \circ Z(W) \circ \beta_M(y) = \alpha_{M'} \circ Z(W) (x) . \tag*{\qed}\]
 \renewcommand{\qed}{}
\end{proof}

Now our goal is to construct a $2$-functor $\tilde{Z}\colon \mathcal{B} \to R\text{-}\Mod_t$ with $\tilde{Z}(M) = V_M$ and $\tilde{Z}(W)\colon V_M \allowbreak \to V_{M'}$ for any $W\colon M \to M'$ as above. From \eqref{eq:expression_Z_tilde} we see that a twist from $Z(W_1)$ to $Z(W_2)$ induces a twist from $\tilde{Z}(W_1)$ to $\tilde{Z}(W_2)$. Note that $\tilde{Z}$ preserves compositions of 1-morphisms if the following additional assumption holds:
\begin{equation*}
 Z(W)(\beta_M(V_M)) \subset \beta_{M'}(V_{M'}) \qquad \text{for any bordism} \quad W\colon M \to M' .
\end{equation*}
Namely if so, let $W \colon M \to M'$ and $W' \colon M' \to M''$ be bordisms. Then
\[ \begin{tikzcd}[row sep=3em, column sep=5em] \beta_M(V_M) \arrow{r}{Z(W)} \arrow[swap]{d}{\alpha_M} & Z(W)(\beta_M(V_M)) \arrow{r}{Z(W')} \arrow{d}{\alpha_{M'}} & Z(W' \circ W)(\beta_M(V_M)) \arrow{d}{\alpha_{M''}} \\ V_M \arrow{r}{\tilde{Z}(W)} 
& V_{M'} \arrow{r}{\tilde{Z}(W')} 
& V_{M''} 
\end{tikzcd} \]
is a commutative diagram by the previous lemma. We have $\tilde{Z}(W') \circ \tilde{Z}(W) \circ \alpha_M = \alpha_{M''} \circ Z(W') \circ Z(W) = \alpha_{M''} \circ Z(W' \circ W)$, so uniqueness implies that $\tilde{Z}(W' \circ W) = \tilde{Z}(W') \circ \tilde{Z}(W)$, and hence~$\tilde{Z}$ is a functor.

Summarizing, we obtain the following definition.
\begin{Definition} \label{def:TQFT_reduction}
 For each object $M$ in $\mathcal{B}$, let $\begin{tikzcd} Z(M) \arrow[shift left=0.25em]{r}{\alpha_M} & \arrow[shift left=0.25em]{l}{\beta_M} N_M \end{tikzcd}$ be $R$-module morphisms with~$N_M$ an $R$-module, and $V_M \subset N_M$ a~submodule. If
 \begin{enumerate}\itemsep=0pt
 \item[(i)] $N_\varnothing = Z(\varnothing)$ and $\alpha_\varnothing, \beta_\varnothing$ are identity maps,
 \item[(ii)] $(\alpha_{M'} \circ Z(W) \circ \beta_M)(V_M) \subset V_{M'}$ for all bordisms $W\colon M \to M'$,
 \item[(iii)] the restriction $\alpha_M \circ \beta_M\colon V_M \to V_M$ is invertible for all $M$,
 \item[(iv)] $Z(W)(\beta_M(V_M)) \subset \beta_{M'}(V_{M'})$ for all bordisms $W\colon M \to M'$,
 \end{enumerate}
 then we speak of a \emph{reduction} of the TQFT, and call the functor $\tilde{Z}$ of Lemma \ref{lemma:TQFT_reduction} the \emph{reduced TQFT}.
\end{Definition}

The whole point of the reduced TQFT $\tilde{Z}$ is that it computes the same invariants as $Z$ for closed manifolds, while allowing for easier computations. Indeed, if $W\colon \varnothing \to \varnothing$ is a bordism, then
\[ \tilde{Z}(W)(1) = \tilde{Z}(W) \circ \alpha_\varnothing (1) = \alpha_\varnothing \circ Z(W)(1) = Z(W)(1) . \]
In Section \ref{sec:app}, we will apply this to the category $\mathcal{B} = \Tb_2(\Lambda)$ and TQFT of Section \ref{sec:TQFT_parabolic_structures}, as follows. We have $Z\big(S^1, \star\big) = \K(\Var/G)$, and there is an action of $G$ on itself by conjugation. Suppose there are conjugacy-closed strata $\mathcal{C}_1, \dots, \mathcal{C}_n$ for $G$, with maps $\pi_i\colon \mathcal{C}_i \to C_i$ whose fibers are precisely the orbits of $G$. Then we have an induced decomposition $K(\Var/G) = \K(\Var/\mathcal{C}_1) \oplus \cdots \oplus \K(\Var/\mathcal{C}_n)$, and the maps $(\pi_i)_!$ and $(\pi_i)^*$ induce maps
\begin{equation}
 \label{eq:usual_reduction}
 \begin{tikzcd} \K(\Var/\mathcal{C}_1) \oplus \cdots \oplus \K(\Var/\mathcal{C}_n) \arrow[shift left=0.25em]{r}{\pi_!} & \arrow[shift left=0.25em]{l}{\pi^*} \K(\Var/C_1) \oplus \cdots \oplus \K(\Var/C_n), \end{tikzcd}
\end{equation}
which by slight abuse of notation we denote by $\pi^*$ and $\pi_!$. For bordisms $W\colon \big(S^1, \star\big) \to (S^1, \star)$, we write $Z_\pi(W) = \pi_! \circ Z(W) \circ \pi^*$ as a shorthand. Let $V \subset \bigoplus_{i = 1}^{n} \K(\Var/C_i)$ be a submodule on which $\eta = \pi_! \pi^*$ is invertible. Now, not any such stratification will satisfy the conditions for a~reduction. To make this precise, consider the variety $\Delta = \big\{ (g_1, g_2) \in G^2 \mid g_1 \sim{} g_2 \big\}$ of the pairs of conjugate elements of $G$, and the conjugation map
\[ c\colon \ G^2 \to \Delta, \qquad (g, h) \mapsto \big(g, h g h^{-1}\big) . \]
The stratification of $G$ naturally induces a stratification of $\Delta$, whose strata we denote by $\Delta_i = \big\{ (g_1, g_2) \in \mathcal{C}_i^2\colon g_1 \sim{} g_2 \big\}$, and the map $c_! c^*$ naturally restricts to a map $\K(\Var/\Delta_i) \to \K(\Var/\Delta_i)$.

\begin{Proposition}
 \label{prop:conditions_for_reduction_TQFT}
 Let $G$ be an algebraic group, stratified by conjugacy-closed strata $\mathcal{C}_i$, with maps $\pi_i\colon \mathcal{C}_i \to C_i$ whose fibers are precisely the orbits of $G$, and let $V \subset \bigoplus_{i = 1}^{n} \K(\Var/C_i)$ be a submodule. Assume that
 \begin{enumerate}\itemsep=0pt
 \item[$(i)$] $Z_\pi(W)(V) \subset V$ for all bordisms $W\colon \big(S^1, \star\big) \to \big(S^1, \star\big)$,
 \item[$(ii)$] the map $\eta = \pi_! \pi^*\colon V \to V$ is invertible,
 \item[$(iii)$] $[\{ 1 \}] \in \pi^* V$ and whenever $X \in V$ then $X_i = X|_{C_i} \in V$ as well,
 \item[$(iv)$] for each stratum $\mathcal{C}_i$, the stabilizers $\Stab_i$ of the points are all isomorphic and special.\footnote{We say that a linear algebraic group $G$ is special if every $G$-torsor (locally trivial in the \'etale topology) is locally trivial in the Zariski topology.}
 \end{enumerate}
 Then the maps in \eqref{eq:usual_reduction} and the submodule $V$ yield a reduction of the TQFT.
\end{Proposition}

\begin{proof}
 The only remaining condition to show is (iv) of Definition \ref{def:TQFT_reduction}, and it suffices to show this holds for the bordisms $D$, $D^\dag$, $L$ and $L_\mathcal{E}$. This holds for~$D$ by assumption~$(iii)$, and for $D^\dag$ trivially because $V_\varnothing = \K(\Var_k)$.
 For~$L$, we consider the associated span $\mathcal{F}(L)$ given by
 \[ \begin{array}{@{}ccccc} G & \overset{p}{\longleftarrow} & G^4 & \overset{q}{\longrightarrow} & G, \\ g & \mapsfrom & (g, g_1, g_2, h) & \mapsto & h g [g_1, g_2] h^{-1}, \end{array} \]
 and also the modified span
 \[
 \begin{array}{@{}ccccc} G & \overset{\tilde{p}}{\longleftarrow} & G^3 & \overset{\tilde{q}}{\longrightarrow} & G, \\ g & \mapsfrom & (g, g_1, g_2) & \mapsto & g [g_1, g_2] . \end{array}
 \]
 It is not hard to see that
 \begin{equation}
 \label{eq:that_one_from_diagram_chasing}
 q_! p^* = (p_2)_! c_! c^* (p_1)^* \tilde{q}_! \tilde{p}^* ,
 \end{equation}
 with $p_1, p_2\colon \Delta \to G$ the projections, as both sides of the equality send $X \overset{f}{\longrightarrow} G$ to
 \[ \begin{array}{@{}ccc} \big\{ (x, g, g_1, g_2, h) \in X \times G^4 \mid g = f(x) \big\} & \longrightarrow & G, \\ (x, g, g_1, g_2, h) & \mapsto & h g [g_1, g_2] h^{-1} . \end{array} \]
 Let us describe the restrictions $(c_! c^*)_i\colon \K(\Var/\Delta_i) \to \K(\Var/\Delta_i)$. For any $X \xrightarrow{f} \Delta_i$, the pullback $c^*X \to X$ is a $\Stab_i$-torsor, so by assumption (iv), it is Zariski-locally trivial, yielding $[c^* X]_X = [\Stab_i] [X]_X$. In particular, the restrictions $(c_! c^*)_i$ are given by scalar multiplication with $[\Stab_i]$.

 Now, take any $X \in V$, let $Y = \tilde{q}_! \tilde{p}^* \pi^* X \in \K(\Var/G)$, and decompose $Y = \sum_{i = 1}^{n} Y_i$ with each $Y_i \in \K(\Var/\mathcal{C}_i)$ according to the stratification of $G$. Then $\pi_! Y = \pi_! \tilde{q}_! \tilde{p}^* \pi^* X = \frac{1}{[G]} Z_\pi(L)(X)$ lies in $V$ by (i), so $\pi_! Y_i = (\pi_! Y)_i$ lies in $V$ by (iii). Furthermore,
 $(p_2)_! c_! c^* (p_1)^* Y_i = [\Stab_i] (p_2)_! (p_1)^* Y_i = [\Stab_i] \pi^* \pi_! Y_i \in \pi^* V$,
 so it follows using \eqref{eq:that_one_from_diagram_chasing} that
 \[ Z(L)(\pi^* X) = \sum_i \pi^* \pi_! Y_i \in \pi^* V . \]
 A completely similar argument shows that $Z(L_\mathcal{E})(\pi^* V) \subset \pi^* V$ as well.
\end{proof}

In this paper, we consider the groups of upper triangular matrices $\U_2$, $\U_3$ and $\U_4$ over $\CC$. For these groups, there exists a natural stratification given by the types of conjugacy classes. In particular, for each stratum, the stabilizers of the points are isomorphic and special, hence satisfying condition (iv) of the proposition above.

Now we see that the reduced TQFT leads to a simplification of the computation, as it allows us to work over $\K(\Var/C_i)$ instead of $\K(\Var/\mathcal{C}_i)$. This is reflected in the computations by the fact that we can get rid of the conjugation by $h$ as in the span of~$\mathcal{F}(L)$~\eqref{eq:span_Z_L}.

\begin{Remark}For the groups $G$ and strata $\mathcal{C}_i$ we consider in the next sections, it might be that the map $\eta = \pi_! \pi^*$ is not invertible as a $\K(\Var_k)$-module morphism. However, one can replace the ring $\K(\Var_k)$ by a suitable localization (often it suffices to invert $[G]$), to make $\eta$ is invertible. As a consequence, the resulting classes $[ \mathfrak{X}_G(X, Q) ]$ will only be defined in that localization. This is not unreasonable, since $[G]$ needs to be invertible anyway in order to apply \eqref{eq:representation_variety_from_TQFT}.
Also in many cases we can still extract algebraic data from the localized class: given a multiplicative system $S \subset \K(\Var_k)$ and an element $\overline{x} \in S^{-1} \K(\Var_k)$ that admits a lift $x \in \K(\Var_k)$, this lift is defined up to a sum of annihilators of elements of $S$. If $\varphi\colon \K(\Var_k) \to R$ is a ring morphism with $R$ a domain such that $\varphi(s) \ne 0$ for all $s \in S$, then $\varphi(a) = 0$ for any annihilator $a$ of any $s \in S$. Hence $\varphi(x)$ is independent on the choice of lift.
The example to have in mind here is the $E$-polynomial $e\colon \K(\Var_\CC) \to \ZZ[u, v]$. Since $e(q) = uv \ne 0$, where $q = \big[ \AA^1_\CC \big] \in \K(\Var_\CC)$, to compute the $E$-polynomial of some variety $X$ over $\CC$ it is sufficient to know its class in the localized ring $S^{-1} \K(\Var_\CC)$ for $S = \big\{ 1, q, q^2, \dots \big\}$. (Similarly we could invert $q - 1$ or $q + 1$.)
\end{Remark}

\section{Applications}\label{sec:app}

In this section, we apply the technique developed in Sections \ref{sec:field} and \ref{sec:reduction_TQFT} to compute the class of the $G$-representation varieties $\mathfrak{X}_G(\Sigma_g)$ in the Grothendieck ring $\K(\Var_\CC)$, with $G$ being the groups of complex upper triangular $n \times n$ matrices for $n = 2, 3, 4$. We prove our main theorems, Theorems \ref{thm:main1} and~\ref{thm:main2}. We also discuss generalizations of Theorem \ref{thm:main1} to representation varieties with parabolic structures.

\subsection[Upper triangular 2 times 2]{Upper triangular $\boldsymbol{2 \times 2}$}

We denote the group of $2 \times 2$ upper triangular matrices over $\CC$ by
\[ \U_2 = \left\{ \begin{pmatrix} a & b \\ 0 & c \end{pmatrix} \bigg|\, a, c \ne 0 \right\} . \]
It is easily seen that the class of $\U_2$ in the Grothendieck ring of varieties is $\big[\AA^1_\CC\big]\big[\AA^1_\CC \setminus \{ 0 \}\big]^2 = q (q - 1)^2$. Moreover, this group contains the following three types of conjugacy classes according to their orbit with respect to conjugation.
\begin{enumerate}\itemsep=0pt
 \item All scalar matrices $\left(\begin{smallmatrix} \lambda & 0 \\ 0 & \lambda \end{smallmatrix}\right)$ have a singleton orbit.
 \item All matrices of the form $\left(\begin{smallmatrix} \lambda & b \\ 0 & \lambda \end{smallmatrix}\right)$ with $b \ne 0$ are conjugate to the Jordan block $\left(\begin{smallmatrix} \lambda & 1 \\ 0 & \lambda \end{smallmatrix}\right)$, and thus have an orbit isomorphic to $\AA^1_\CC \setminus \{ 0 \}$.
 \item All remaining matrices in $\U_2$ are of the form $\left(\begin{smallmatrix} \lambda & b \\ 0 & \mu \end{smallmatrix}\right)$, with $\lambda, \mu \ne 0$, which are conjugate if and only if they have the same diagonal. Hence, these matrices have an orbit isomorphic to $\AA^1_\CC$.
\end{enumerate} We denote these three types of conjugacy classes by
\begin{gather*}
 \mathcal{S} = \left\{ \begin{pmatrix} \lambda & 0 \\ 0 & \lambda \end{pmatrix} \bigg|\, \lambda \ne 0 \right\}, \qquad
 \mathcal{J} = \left\{ \begin{pmatrix} \lambda & b \\ 0 & \lambda \end{pmatrix} \bigg|\, \lambda, b \ne 0 \right\}, \\
 \mathcal{M} = \left\{ \begin{pmatrix} \lambda & b \\ 0 & \mu \end{pmatrix} \bigg|\, \lambda, \mu \ne 0,\, \lambda \ne \mu,\, b \in \CC \right\}
 \end{gather*}
and the orbits of their elements by
\[ \mathcal{S}_{\lambda} = \left\{ \begin{pmatrix} \lambda & 0 \\ 0 & \lambda \end{pmatrix} \right\}, \qquad
 \mathcal{J}_{\lambda} = \left\{ \begin{pmatrix} \lambda & b \\ 0 & \lambda \end{pmatrix} \bigg|\, b \ne 0 \right\}, \qquad
 \mathcal{M}_{\lambda, \mu} = \left\{ \begin{pmatrix} \lambda & b \\ 0 & \mu \end{pmatrix} \bigg|\, b \in \CC \right\} \]
for any $\lambda, \mu \ne 0$ with $\lambda \ne \mu$. It is easy to see that the classes of these varieties in $\K(\Var_\CC)$ are given by
 \begin{alignat*}{4}
& [ \mathcal{S} ] = q - 1, \qquad && [ \mathcal{J} ] = (q - 1)^2, \qquad && [ \mathcal{M} ] = q (q - 1) (q - 2), & \\
& [ \mathcal{S}_{\lambda} ] = 1, \qquad && [ \mathcal{J}_{\lambda} ] = q - 1, \qquad && [ \mathcal{M}_{\lambda, \mu} ] = q.&
\end{alignat*}
We denote the orbit spaces by $S$, $J$ and $M = \{ (\lambda, \mu) \in \CC^* \times \CC^* \mid \lambda \ne \mu \}$ respectively. We remark that both $S$ and $J$ are isomorphic to $\AA^1_\CC \setminus \{ 0 \}$ as varieties. The quotient maps that identify the orbits are
\begin{gather*}
 \pi_{\mathcal{S}}\colon \ \mathcal{S} \to S\colon \ \begin{pmatrix} \lambda & 0 \\ 0 & \lambda \end{pmatrix} \mapsto \lambda, \qquad \pi_{\mathcal{J}}\colon \ \mathcal{J} \to J\colon \ \begin{pmatrix} \lambda & b \\ 0 & \lambda \end{pmatrix} \mapsto \lambda, \\ \pi_{\mathcal{M}}\colon \ \mathcal{M} \to M\colon \ \begin{pmatrix} \lambda & b \\ 0 & \mu \end{pmatrix} \mapsto (\lambda, \mu) .
 \end{gather*}
These maps induce the morphisms
\[ \begin{tikzcd} \K(\Var/\mathcal{S}) \oplus \K(\Var/\mathcal{J}) \oplus \K(\Var/\mathcal{M}) \arrow[shift left=0.25em]{r}{\pi_!} & \arrow[shift left=0.25em]{l}{\pi^*} \K(\Var/S) \oplus \K(\Var/J) \oplus \K(\Var/M) \end{tikzcd} . \]

To obtain a reduction of the TQFT for $G = \U_2$ with stratification given by $\mathcal{S}$, $\mathcal{J}$ and $\mathcal{M}$, let us show that condition~(iv) of Proposition~\ref{prop:conditions_for_reduction_TQFT} holds.

\begin{Lemma} The stabilizer subgroup of any $g \in G$ is isomorphic to
 \[ \Stab_\mathcal{S} = \U_2, \qquad \Stab_\mathcal{J} = \GG_m \times \GG_a, \qquad \textup{or} \qquad \Stab_\mathcal{M} = \GG_m \times \GG_m , \]
 depending whether $g$ lies in $\mathcal{S}$, $\mathcal{J}$ or $\mathcal{M}$, respectively. Moreover, all of these subgroups are special.
\end{Lemma}
\begin{proof}
 The proof is straightforward, we leave it to the reader.
\end{proof}

As can be seen in the proof of Proposition \ref{prop:conditions_for_reduction_TQFT}, the restriction of $c_! c^*$ to $\Delta_{\mathcal{E}}$ is now just multiplication by $[\Stab_\mathcal{E}]$ for each stratum $\mathcal{E} = \mathcal{S}, \mathcal{J}$ or $\mathcal{M}$.

\begin{Remark}
 Alternatively, one could show even more explicitly that the restriction of $c_! c^*$ to~$\Delta_\mathcal{E}$ is multiplication by~$[\Stab_\mathcal{E}]$. Consider the representatives of the conjugacy classes
 \[ \xi^\mathcal{S}_\lambda = \left( \begin{smallmatrix} \lambda & 0 \\ 0 & \lambda \end{smallmatrix} \right), \qquad \xi^\mathcal{J}_\lambda = \left( \begin{smallmatrix} \lambda & 1 \\ 0 & \lambda \end{smallmatrix} \right), \qquad \xi^\mathcal{M}_{\lambda, \mu} = \left( \begin{smallmatrix} \lambda & 0 \\ 0 & \mu \end{smallmatrix} \right) . \]
 Their stabilizers are explicitly given by
 \[ \Stab\big(\xi^{\mathcal{S}}_\lambda\big) = \U_2, \qquad \Stab\big(\xi^{\mathcal{J}}_\lambda\big) = \left\{ \left(\begin{smallmatrix} \alpha & \beta \\ 0 & \alpha \end{smallmatrix}\right) \mid \alpha \ne 0 \right\}, \qquad \Stab\big(\xi^{\mathcal{M}}_{\lambda, \mu}\big) = \left\{ \left(\begin{smallmatrix} \alpha & 0 \\ 0 & \beta \end{smallmatrix}\right) \mid \alpha, \beta \ne 0 \right\} . \]
 For all $\mathcal{E} = \mathcal{S}, \mathcal{J}, \mathcal{M}$ it is straightforward to come up with a map $\sigma\colon \mathcal{E} \to G$ such that $g = \sigma(g) \xi^{\mathcal{E}}_{\pi_\mathcal{E}(g)} \sigma(g)^{-1}$ for any $g \in \mathcal{E}$. For example, for $\mathcal{E} = \mathcal{J}$ we can take $\sigma\left(\begin{smallmatrix} \lambda & b \\ 0 & \lambda \end{smallmatrix}\right) = \left(\begin{smallmatrix} 1 & 0 \\ 0 & 1/b \end{smallmatrix}\right)$, 
 and for $\mathcal{E} = \mathcal{M}$ one can take $\sigma\left(\begin{smallmatrix} \lambda & b \\ 0 & \mu \end{smallmatrix}\right) = \left(\begin{smallmatrix} 1 & b / (\mu - \lambda) \\ 0 & 1 \end{smallmatrix}\right)$.
 Now, for any variety $X \overset{(f_1, f_2)}{\longrightarrow} \Delta_{\mathcal{E}}$ we have an isomorphism
\begin{gather*}
 X \times \Stab_\mathcal{E} \xrightarrow{\sim} c_! c^* X = \big\{ (x, h) \in X \times G \mid f_2(x) = h f_1(x) h^{-1} \big\}, \\
 (x, s) \mapsto \big(x, \sigma(f_2(x)) s \sigma(f_1(x))^{-1} \big),
 \end{gather*}
 which shows that $c_! c^* [X]_{\Delta_\mathcal{E}} = [\Stab_\mathcal{E}] [X]_{\Delta_\mathcal{E}}$.
\end{Remark}

Write $T_{S_\lambda} \in \K(\Var/S), T_{J_\lambda} \in \K(\Var/J)$ and $T_{M_{\lambda, \mu}} \in \K(\Var/M)$ for the classes of the points $\{ \lambda \} \to S, \{ \lambda \} \to J$ and $\{ (\lambda, \mu) \} \to M$. We consider the submodule $V$ generated by these classes $V = \langle T_{S_\lambda}, T_{J_\lambda}, T_{M_{\lambda, \mu}} \rangle$. From the computations that follow, it will be clear that $V$ is invariant under $\eta = \pi_! \circ \pi^*$ and $Z_\pi$. Hence all conditions from Proposition \ref{prop:conditions_for_reduction_TQFT} are satisfied, so we have a reduction of the TQFT.

Since all fibrations $\mathcal{S} \to S$, $\mathcal{J} \to J$ and $\mathcal{M} \to M$ are trivial, we immediately find that
\begin{gather*}
 \eta(T_{S_\lambda}) = [ \mathcal{S}_\lambda ] T_{S_\lambda} = T_{S_\lambda}, \qquad \eta(T_{J_\lambda}) = [ \mathcal{J}_\lambda] T_{J_\lambda} = (q - 1) T_{J_\lambda}, \\ \eta(T_{M_{\lambda, \mu}}) = [ \mathcal{M}_{\lambda, \mu} ] T_{M_{\lambda, \mu}} = q T_{M_{\lambda, \mu}} ,
 \end{gather*}
that is,
\[ \eta = \left( \begin{array}{c|ccc}
& T_{S_\lambda} & T_{J_\lambda} & T_{M_{\lambda, \mu}} \\
\hline
T_{S_\lambda} & 1 & 0 & 0 \\
T_{J_\lambda} & 0 & q - 1 & 0 \\
T_{M_{\lambda, \mu}} & 0 & 0 & q
\end{array} \right). \]

For computing $Z_\pi(L)$, recall from \eqref{eq:span_Z_L} that $\mathcal{F}(L)$ is given by
\[ \begin{tikzcd}[row sep=0em] \U_2 & \arrow[swap]{l}{p} \U_2^4 \arrow{r}{q} & \U_2, \\ g & \arrow[maps to]{l} (g, g_1, g_2, h) \arrow[maps to]{r} & h g {[} g_1, g_2 {]} h^{-1} . \end{tikzcd} \]

 First we compute $Z_\pi(L)(T_{S_\lambda})$. We have $\pi^*(T_{S_\lambda}) = [ \mathcal{S}_\lambda ]_{\U_2}$. Note that for any group elements $g_1=\left( \begin{smallmatrix} a_1 & b_1 \\ 0 & c_1 \end{smallmatrix} \right)$ and $g_2=\left( \begin{smallmatrix} a_2 & b_2 \\ 0 & c_2 \end{smallmatrix} \right)$, the commutator $[g_1, g_2] = \left( \begin{smallmatrix} 1 & x \\ 0 & 1 \end{smallmatrix} \right)$ with
 \[x = \frac{a_1 b_2 - a_2 b_1 + b_1 c_2 - b_2 c_1}{c_1 c_2}.\]
 Hence $q(p^*(\mathcal{S}_\lambda)) \subset \mathcal{S}_\lambda \cup \mathcal{J}_\lambda$, and thus $Z_\pi(L)(T_{S_\lambda})$ are generated by $T_{S_\lambda}$ and $T_{J_\lambda}$. Then, we have that
 \begin{align*}
 Z_\pi(L)(T_{S_\lambda})|_{T_{S_\lambda}}
 &= \big[ \mathcal{S}_\lambda \times \U_2^3 \cap q^{-1}(\mathcal{S}_\lambda) \big] \\
 &= [ \{ g_1, g_2 \in \U_2 \mid [g_1, g_2] = 1 \} ] \cdot [ \U_2 ] \\
 &= \left[ \left\{ a_1, b_1, c_1, a_2, b_2, c_2 \in \CC \mid \substack{a_1 b_2 - a_2 b_1 + b_1 c_2 - b_2 c_1 = 0 \\ \text{ and } a_1 c_1 a_2 c_2 \ne 0} \right\} \right] \cdot [ \U_2 ] .
 \end{align*}
 We cut up the variety
 \[[ \left\{ a_1, b_1, c_1, a_2, b_2, c_2 \in \CC \mid a_1 b_2 - a_2 b_1 + b_1 c_2 - b_2 c_1 = 0 \text{ and } a_1 c_1 a_2 c_2 \ne 0 \right\} ]\]
 into three pieces given by extra conditions: $a_1\ne c_1$; $a_1=c_1$, $a_2\ne c_2$; and finally $a_1=c_1$, $a_2=c_2$. In the first case, the equation $a_1 b_2 - a_2 b_1 + b_1 c_2 - b_2 c_1 = 0$ can be solved for $b_2$ given the values of the other variables, hence
 \[ [ \left\{a_1 b_2 - a_2 b_1 + b_1 c_2 - b_2 c_1 = 0,\text{ } a_1 c_1 a_2 c_2 \ne 0 \text{ and } a_1\ne c_1\right\} ] = q(q - 1)^3 (q - 2) . \]
 In the second case, the equations $a_1=c_1$ and $a_2\ne c_2$ yield $b_1=0$, hence
 \[[ \left\{a_1 b_2 - a_2 b_1 + b_1 c_2 - b_2 c_1 = 0,\text{ } a_1 c_1 a_2 c_2 \ne 0,\text{ } a_1=c_1 \text{ and } a_2\ne c_2\right\} ]=q(q-1)^2(q-2).\]
 Finally, in the third case, the equations $a_1=c_1$ and $a_2=c_2$ imply that
 \[a_1 b_2 - a_2 b_1 + b_1 c_2 - b_2 c_1 = 0\]
 providing that
 \[[ \left\{a_1 b_2 - a_2 b_1 + b_1 c_2 - b_2 c_1 = 0,\text{ } a_1 c_1 a_2 c_2 \ne 0,\text{ } a_1=c_1 \text{ and } a_2= c_2\right\} ]=q^2(q-1)^2.\]
 As a consequence, we obtain
 \[[ \left\{a_1 b_2 - a_2 b_1 + b_1 c_2 - b_2 c_1 = 0 \text{ and } a_1 c_1 a_2 c_2\ne 0 \right\} ]=q^2(q-1)^3\]
 implying that
 \[Z_\pi(L)(T_{S_\lambda})|_{T_{S_\lambda}}=q^2 (q - 1)^3 \cdot [\U_2] = q^3 (q - 1)^5.\]
 Now it follows that $Z_\pi(L)(T_{S_\lambda})|_{T_{J_\lambda}}$ is simply $\big[\mathcal{S}_\lambda \times \U_2^3 \big] - q^3 (q - 1)^5 = q^3 (q - 1)^5 (q - 2)$.

Next we compute $Z_\pi(L)(T_{J_\lambda})$. We have $\pi^*(T_{J_\lambda}) = [ \mathcal{J}_\lambda ]_{\U_2}$. By the same observation as above about the commutator, we see that $Z_\pi(L)(T_{J_\lambda})$ is also generated by $T_{S_\lambda}$ and $T_{J_\lambda}$. Note that
 \[ \left(\begin{smallmatrix} \lambda & b \\ 0 & \lambda \end{smallmatrix}\right) [g_1, g_2] = \left(\begin{smallmatrix} \lambda & 0 \\ 0 & \lambda \end{smallmatrix}\right) \qquad \text{if and only if} \quad [g_1, g_2] = \left(\begin{smallmatrix} 1 & - b / \lambda \\ 0 & 1 \end{smallmatrix}\right) , \]
 which implies that
 \[Z_\pi(L)(T_{J_\lambda})|_{T_{S_\lambda}} = Z_\pi(L)(T_{S_\lambda})|_{T_{J_\lambda}} = q^3 (q - 1)^5 (q - 2).\]
 Now it follows that
 \[Z_\pi(L)(T_{J_\lambda})|_{T_{J_\lambda}} = \big[ \mathcal{J}_\lambda \times \U_2^3 \big] - q^3 (q - 1)^5 (q - 2) = q^3 (q - 1)^5 \big(q^2 - 3q + 3\big).\]

Lastly we compute $Z_\pi(L)(T_{M_{\lambda, \mu}})$. We have $\pi^*(T_{M_{\lambda, \mu}}) = [ \mathcal{M}_{\lambda, \mu} ]_{\U_2}$. By the observation about the commutator, we immediately see that $Z_\pi(L)(T_{M_{\lambda, \mu}})$ must be generated by $T_{M_{\lambda, \mu}}$. Therefore, $Z_\pi(L)(T_{M_{\lambda, \mu}}) = \big[ \mathcal{M}_{\lambda, \mu} \times \U_2^3 \big] T_{M_{\lambda, \mu}} = q^4 (q - 1)^6 T_{M_{\lambda, \mu}}$.

In summary,
\[ Z_\pi(L) = q^3 (q - 1)^5 \left( \begin{array}{c|ccc}
& T_{S_\lambda} & T_{J_\lambda} & T_{M_{\lambda, \mu}} \\
\hline
T_{S_\lambda} & 1 & q - 2 & 0 \\
T_{J_\lambda} & q - 2 & q^2 - 3q + 3 & 0 \\
T_{M_{\lambda, \mu}} & 0 & 0 & q(q - 1)
\end{array} \right). \]

Now, we turn our attention to the cylinder with parabolic structures and we compute $Z_\pi(L_{\mathcal{S}_\lambda})$, $Z_\pi(L_{\mathcal{J}_\lambda})$ and $Z_\pi(L_{\mathcal{M}_{\lambda, \mu}})$. Recall that $\mathcal{F}(L_\mathcal{E})$ is given by
\[ \begin{tikzcd}[row sep=0em] \U_2 & \arrow[swap]{l}{r} \U_2^2 \times \mathcal{E} \arrow{r}{s} & \U_2, \\ g & \arrow[maps to]{l} (g, h, \xi) \arrow[maps to]{r} & h g \xi h^{-1} . \end{tikzcd} \]

Let $g \in \mathcal{S}_\lambda$, and note that if $\xi \in \mathcal{S}_\sigma$ then $g \xi \in \mathcal{S}_{\lambda \sigma}$, if $\xi \in \mathcal{J}_\sigma$ then $g \xi \in \mathcal{J}_{\lambda \sigma}$, and if $\xi \in \mathcal{M}_{\sigma, \rho}$ then $g \xi \in \mathcal{M}_{\lambda \sigma, \lambda \rho}$. Hence we have that $Z_\pi(L_{\mathcal{S}_\lambda})$ is given by
 \begin{gather*} \left( \begin{array}{c|ccc} & T_{S_\sigma} & T_{J_\sigma} & T_{M_{\sigma, \rho}} \\ \hline
 T_{S_{\lambda \sigma}} & [ \mathcal{S}_\sigma \times \U_2] & 0 & 0 \\
 T_{J_{\lambda \sigma}} & 0 & [ \mathcal{J}_\sigma \times \U_2] & 0 \\
 T_{M_{\lambda \sigma, \lambda \rho}} & 0 & 0 & [ \mathcal{M}_{\sigma, \rho} \times \U_2]
 \end{array} \right) \\
 \qquad{} = q (q - 1)^2 \left( \begin{array}{c|ccc} & T_{S_\sigma} & T_{J_\sigma} & T_{M_{\sigma, \rho}} \\ \hline
 T_{S_{\lambda \sigma}} & 1 & 0 & 0 \\
 T_{J_{\lambda \sigma}} & 0 & q - 1 & 0 \\
 T_{M_{\lambda \sigma, \lambda \rho}} & 0 & 0 & q
 \end{array} \right) . \end{gather*}

Now let $g \in \mathcal{J}_\lambda$. We see that if $\xi \in \mathcal{S}_\sigma$ then $g \xi \in \mathcal{J}_{\lambda \sigma}$, and if $\xi \in \mathcal{M}_{\sigma, \rho}$ then $g \xi \in \mathcal{M}_{\lambda \sigma, \lambda \rho}$. If $\xi \in \mathcal{J}_\sigma$, then $g \xi \in \mathcal{S}_{\lambda \sigma}$ precisely if $g = \lambda \sigma \xi^{-1}$ and otherwise $g \xi \in \mathcal{J}_{\lambda \sigma}$. Hence we have
 \begin{align*}
 Z_\pi(L_{\mathcal{J}_\lambda}) & = \left( \begin{array}{c|ccc} & T_{S_\sigma} & T_{J_\sigma} & T_{M_{\sigma, \rho}} \\ \hline
 T_{S_{\lambda \sigma}} & 0 & [ \mathcal{J}_\sigma \times \U_2 ] & 0 \\
 T_{J_{\lambda \sigma}} & [ \mathcal{S}_\sigma \times \U_2] [ \mathcal{J}_\lambda ] & [ \mathcal{J}_\sigma \times \U_2 ] ([ \mathcal{J}_\lambda ] - 1) & 0 \\
 T_{M_{\lambda \sigma, \lambda \rho}} & 0 & 0 & [ \mathcal{M}_{\sigma, \rho} \times \U_2] [ \mathcal{J}_\lambda ]
 \end{array} \right) \\
 & = q (q - 1)^2 \left( \begin{array}{c|ccc} & T_{S_\sigma} & T_{J_\sigma} & T_{M_{\sigma, \rho}} \\ \hline
 T_{S_{\lambda \sigma}} & 0 & q - 1 & 0 \\
 T_{J_{\lambda \sigma}} & q - 1 & (q - 1)(q - 2) & 0 \\
 T_{M_{\lambda \sigma, \lambda \rho}} & 0 & 0 & q (q - 1)
 \end{array} \right) .
 \end{align*}

Lastly, let $g \in \mathcal{M}_{\lambda, \mu}$. If $\xi \in \mathcal{S}_\sigma$ then $g \xi \in \mathcal{M}_{\lambda \sigma, \mu \sigma}$, and if $\xi \in \mathcal{J}_\sigma$ then $g \xi \in \mathcal{M}_{\lambda \sigma, \mu \sigma}$ as well. If $\xi \in \mathcal{M}_{\sigma, \rho}$, then $g \xi \in \mathcal{M}_{\lambda \sigma, \mu \rho}$ if $\lambda \sigma \ne \mu \rho$ and otherwise $g \xi \in \mathcal{S}_{\lambda \sigma}$ precisely for $g = \lambda \sigma \xi^{-1}$ and else $g \xi \in \mathcal{J}_{\lambda \sigma}$. Hence we see that $Z_\pi(L_{\mathcal{M}_{\lambda, \mu}})$ is given by
 \begin{gather*}
\left( \begin{array}{@{\,}c@{\,}|@{\,}c@{\,}c@{\,}c@{\,}c@{\,}} & T_{S_\sigma} & T_{J_\sigma} & T_{M_{\sigma, \rho}} & T_{M_{\sigma', \rho'}} \\ \hline
 T_{S_{\lambda \sigma}} & 0 & 0 & 0 & {[} \mathcal{M}_{\sigma', \rho'} \times \U_2 {]} \\
 T_{J_{\lambda \sigma}} & 0 & 0 & 0 & {[} \mathcal{M}_{\sigma', \rho'} \!\times\! \U_2 {]} ([ \mathcal{M}_{\lambda, \mu} ] - 1 ) \\
 T_{M_{\lambda \sigma, \lambda \rho}} & [ \mathcal{S}_\sigma \!\times\! \U_2 ] [ \mathcal{M}_{\lambda, \mu} ] & [ \mathcal{J}_\sigma \!\times\! \U_2 ] [ \mathcal{M}_{\lambda, \mu} ] & {[} \mathcal{M}_{\sigma, \rho} \!\times\! \U_2 {]} {[} \mathcal{M_{\lambda, \mu}} {]} & 0
 \end{array} \right)
 \\
 \qquad{} = q (q - 1)^2 \left( \begin{array}{c|cccc} & T_{S_\sigma} & T_{J_\sigma} & T_{M_{\sigma, \rho}} & T_{M_{\sigma', \rho'}} \\ \hline
 T_{S_{\lambda \sigma}} & 0 & 0 & 0 & q \\
 T_{J_{\lambda \sigma}} & 0 & 0 & 0 & q (q - 1) \\
 T_{M_{\lambda \sigma, \lambda \rho}} & q & q (q - 1) & q^2 & 0
 \end{array} \right) ,
 \end{gather*}
 with $\lambda \sigma \ne \mu \rho$ and $\lambda \sigma' = \mu \rho'$.

Since the conditions of Proposition~\ref{prop:conditions_for_reduction_TQFT} are satisfied, we can consider the reduced TQFT, $\tilde{Z}$. We have $\tilde{Z}(L) = Z_\pi(L) \circ \eta^{-1}$, so
\[\tilde{Z}(L)= q^3 (q - 1)^4 \left( \begin{array}{c|ccc}
& T_{S_\lambda} & T_{J_\lambda} & T_{M_{\lambda, \mu}} \\
\hline
T_{S_\lambda} & q-1 & q - 2 & 0 \\
T_{J_\lambda} & (q - 2)(q - 1) & q^2 - 3q + 3 & 0 \\
T_{M_{\lambda, \mu}} & 0 & 0 & (q - 1)^2
\end{array} \right).\]
We can diagonalize this matrix as
\[
\tilde{Z}(L) = q^3 (q - 1)^4 A \begin{pmatrix} 1 & 0 & 0 \\ 0 & (q - 1)^2 & 0 \\ 0 & 0 & (q - 1)^2 \end{pmatrix} A^{-1} \qquad \text{with} \quad A = \begin{pmatrix} 1 & 1 & 0 \\ -1 & q - 1 & 0 \\ 0 & 0 & 1 \end{pmatrix} ,
\]
which yields
\begin{gather}
 \tilde{Z}(L^g) = q^{3g - 1} (q - 1)^{4g} \nonumber\\
 \hphantom{\tilde{Z}(L^g) =}{}
 \times \left( \begin{array}{c|ccc} & T_{S_\lambda} & T_{J_\lambda} & T_{M_{\lambda, \mu}} \\ \hline T_{S_\lambda} & (q - 1)((q - 1)^{2g - 1} + 1) & (q - 1)^{2g} - 1 & 0 \\ T_{J_\lambda} & (q - 1)((q - 1)^{2g} - 1) & (q - 1)^{2g + 1} + 1 & 0 \\ T_{M_{\lambda, \mu}} & 0 & 0 & q (q - 1)^{2g} \end{array} \right) . \label{eq:TQFT_U2_genus_power}
\end{gather}

In particular, we proved Theorem \ref{thm:main1} for the group $G=\U_2$.

\begin{Theorem}\label{thm:result_u2} The class of the representation variety $\mathfrak{X}_{\U_2}(\Sigma_g)$ in the localized Grothendieck ring of varieties is
 \begin{equation} \label{eq:class_representation_variety_U2_genus}
 [ \mathfrak{X}_{\U_2}(\Sigma_g) ] = \frac{1}{[ \U_2 ]^{g}} \tilde{Z}(L^g)(T_{S_1})|_{T_{S_1}} = q^{2g - 1} (q - 1)^{2g + 1} \big((q - 1)^{2g - 1} + 1\big) .
\end{equation}
\end{Theorem}

\begin{Remark} \label{rem:result_genus_U2}
 For small values of $g$, we find
 \begin{gather*}
 [ \mathfrak{X}_{\U_2}(\Sigma_1) ] = q^{2} (q - 1)^3 , \\
 [ \mathfrak{X}_{\U_2}(\Sigma_2) ] = q^{4} (q - 1)^5 \big(q^2 - 3q + 3\big) , \\
 [ \mathfrak{X}_{\U_2}(\Sigma_3) ] = q^{6} (q - 1)^7 \big(q^4 - 5q^3 + 10q^2 - 10q + 5\big) , \\
 [ \mathfrak{X}_{\U_2}(\Sigma_4) ] = q^{8} (q - 1)^{9} \big(q^{6} - 7 q^{5} + 21 q^{4} - 35 q^{3} + 35 q^{2} - 21 q + 7\big) .
 \end{gather*}
 Note that $[ \mathfrak{X}_{\U_2}(\Sigma_g) ]$ has a factor $(q - 1)^{2g + 1}$, which can be explained as follows. There is a~free action of $\GG_m^{2g}$ on $\mathfrak{X}_{\U_2}(\Sigma_g)$ given by scaling the $A_i, B_i$ (notation as in \eqref{eq:explicit_expression_representation_variety}). This yields a $\GG_m^{2g}$-torsor $\mathfrak{X}_{\U_2}(\Sigma_g) \to \mathfrak{X}_{\U_2}(\Sigma_g) \sslash \GG_m^{2g}$, which by is trivial in the Zariski topology as $\GG_m^{2g}$ is a special group. Hence, $[\mathfrak{X}_{\U_2}(\Sigma_g)]$ is divisible by $[\GG_m^{2g}] = (q - 1)^{2g}$. For the remaining factor $(q - 1)$, let $D \subset \mathfrak{X}_{\U_2}(\Sigma_g)$ be the subvariety where all $A_i, B_i$ are diagonal. Then $[D] = (q - 1)^{4g}$ and there is a free action of $\GG_m$ on $\mathfrak{X}_{\U_2}(\Sigma_g) \backslash D$ given by conjugation with $\left( \begin{smallmatrix} 1 & 0 \\ 0 & x \end{smallmatrix} \right)$ for $x \in \CC^*$.

 In fact, the affine GIT quotient $\mathfrak{X}_{\U_2}(\Sigma_g) \sslash \GG_m^{2g}$ can be identified with the representation variety $[\mathfrak{X}_{\AGL_1}(\Sigma_g)]$, where
 \[ \AGL_1 = \left\{ \begin{pmatrix} a & b \\ 0 & 1 \end{pmatrix}\colon a \ne 0 \right\} \]
 is the general affine group of the line. Therefore, we obtain that
 \[ [\mathfrak{X}_{\AGL_1}(\Sigma_g)] = \frac{[ \mathfrak{X}_{\U_2}(\Sigma_g) ]}{(q - 1)^{2g}} = q^{2g - 1} (q - 1)\big((q - 1)^{2g - 1} + 1\big) , \]
 recovering a result of Gonz\'alez-Prieto, Logares, and Mu\~noz \cite{arXiv200501841}. Alternatively, this result can also be obtained from the isomorphism of algebraic groups $\GG_m \times \AGL_1 \xrightarrow{\sim} \U_2$ which maps $\left(t, \left( \begin{smallmatrix} a & b \\ 0 & 1 \end{smallmatrix} \right)\right)$ to $\left(\begin{smallmatrix} t a & t b \\ 0 & t \end{smallmatrix}\right)$, as it implies
 \[
 [\mathfrak{X}_{\U_2}(\Sigma_g)] = [\mathfrak{X}_{\GG_m}(\Sigma_g)] [\mathfrak{X}_{\AGL_1}(\Sigma_g)],
 \]
  where $[\mathfrak{X}_{\GG_m}(\Sigma_g)] = (q - 1)^{2g}$.
\end{Remark}

Now, we compute the classes of the twisted representation varieties $\mathfrak{X}_{\U_2}(\Sigma_g, Q)$. Similar calculations as before shows that
\begin{gather} \tilde{Z}(L_{\mathcal{S}_\lambda}) = Z_\pi(L_{\mathcal{S}_\lambda})\circ \eta^{-1}= q (q - 1)^2 \left( \begin{array}{c|ccc} & T_{S_\sigma} & T_{J_\sigma} & T_{M_{\sigma, \rho}} \\ \hline T_{S_{\lambda \sigma}} & 1 & 0 & 0 \\ T_{J_{\lambda \sigma}} & 0 & 1 & 0 \\ T_{M_{\lambda \sigma, \mu \sigma}} & 0 & 0 & 1 \end{array} \right) ,
\nonumber\\
 \label{eq:TQFT_U2_parabolic_jordan}
 \tilde{Z}(L_{\mathcal{J}_\lambda}) = Z_\pi(L_{\mathcal{J}_\lambda})\circ \eta^{-1}= q (q - 1)^2 \left( \begin{array}{c|ccc} & T_{S_\sigma} & T_{J_\sigma} & T_{M_{\sigma, \rho}} \\ \hline T_{S_{\lambda \sigma}} & 0 & 1 & 0 \\ T_{J_{\lambda \sigma}} & q - 1 & q - 2 & 0 \\ T_{M_{\lambda \sigma, \mu \sigma}} & 0 & 0 & q - 1 \end{array} \right) ,
\\
 \label{eq:TQFT_U2_parabolic_two_eigenvalues}
 \tilde{Z}(L_{\mathcal{M}_{\lambda, \mu}}) =Z_\pi(L_{\mathcal{M}_{\lambda,\mu}})\circ \eta^{-1}= q (q - 1)^2 \left( \begin{array}{c|cccc} & T_{S_\sigma} & T_{J_\sigma} & T_{M_{\sigma, \rho}} & T_{M_{\sigma', \rho'}} \\ \hline T_{S_{\lambda \sigma}} & 0 & 0 & 0 & 1 \\ T_{J_{\lambda \sigma}} & 0 & 0 & 0 & q - 1 \\ T_{M_{\lambda \sigma, \mu \rho}} & q & q & q & 0 \end{array} \right) ,
\end{gather}
with $\lambda \sigma \ne \mu \rho$ but $\lambda \sigma' = \mu \rho'$.

As a consequence, we obtain the classes of the twisted representation varieties $\mathfrak{X}_{\U_2}(\Sigma_g, Q)$.

\begin{Theorem} \label{thm:result_U2_parabolic}
 Let $\Sigma_g$ be a surface of genus $g$, with parabolic data $Q = \{ (S_1, \mathcal{J}_{\lambda_1}), \dots, (S_k, \mathcal{J}_{\lambda_k}),$ \\ $(S_{k+1}, \mathcal{M}_{\mu_1, \sigma_1}), \dots, (S_{k+\ell}, \mathcal{M}_{\mu_\ell, \sigma_\ell}) \}$.
 \begin{enumerate}\itemsep=0pt
 \item[$(i)$] If $\prod_{i = 1}^{k} \lambda_i \prod_{j = 1}^{\ell} \mu_j \ne 1$ or $\prod_{i = 1}^{k} \lambda_i \prod_{j = 1}^{\ell} \sigma_j \ne 1$, then
 \[ [ \mathfrak{X}_{\U_2}(\Sigma_g, Q) ] = 0 . \]

 \item[$(ii)$] Otherwise, and if $\ell = 0$, then
 \[ [ \mathfrak{X}_{\U_2}(\Sigma_g, Q) ] = q^{2g - 1} (q - 1)^{2g} \big( (-1)^k (q - 1) + (q - 1)^{2g + k} \big) , \]

 \item[$(iii)$] and if $\ell > 0$, then
 \[ [ \mathfrak{X}_{\U_2}(\Sigma_g, Q) ] = q^{2g + \ell - 1} (q - 1)^{4g + k} . \]
 \end{enumerate}
\end{Theorem}

\begin{proof} First note that $(\Sigma_g, Q)$ can be seen as the composition
 \[ D^\dag \circ L^g \circ L_{\mathcal{J}_{\lambda_1}} \circ \cdots \circ L_{\mathcal{J}_{\lambda_k}} \circ L_{\mathcal{M}_{\mu_1, \sigma_1}} \circ \cdots L_{\mathcal{M}_{\mu_\ell, \sigma_\ell}} \circ D . \]

(i) From expressions \eqref{eq:TQFT_U2_genus_power}, \eqref{eq:TQFT_U2_parabolic_jordan} and \eqref{eq:TQFT_U2_parabolic_two_eigenvalues}, we can see that
\[
Z\big(L^g \circ L_{\mathcal{J}_{\lambda_1}} \circ \cdots \circ L_{\mathcal{J}_{\lambda_k}} \circ L_{\mathcal{M}_{\mu_1, \sigma_1}} \circ \cdots \circ L_{\mathcal{M}_{\mu_\ell, \sigma_\ell}}\big)(T_{S_1})|_{T_{S_1}} = 0,
\]
 and hence $[ \mathfrak{X}_{\U_2}(\Sigma_g, Q) ] = 0$.

(ii) Using \eqref{eq:TQFT_U2_parabolic_jordan} and the diagonalization
 \[ \begin{pmatrix} 0 & 1 & 0 \\ q - 1 & q - 2 & 0 \\ 0 & 0 & q - 1 \end{pmatrix} = A \begin{pmatrix} -1 & 0 & 0 \\ 0 & q - 1 & 0 \\ 0 & 0 & q - 1 \end{pmatrix} A^{-1}\qquad \text{with} \quad A = \begin{pmatrix} -1 & \frac{1}{q - 1} & 0 \\ 1 & 1 & 0 \\ 0 & 0 & 1 \end{pmatrix}, \]
 we find that
 \begin{gather*}
 \tilde{Z}(L_{\mathcal{J}_{\lambda_1}} \circ \cdots \circ L_{\mathcal{J}_{\lambda_k}})(T_{S_1}) \\
 \qquad{} = q^{k - 1} (q - 1)^{2k} \big( (-1)^k (q - 1) + (q - 1)^k \big) T_{S_\lambda} + \big( (-1)^{k + 1} (q - 1) + (q - 1)^{k + 1} \big) T_{J_\lambda} ,
 \end{gather*}
 where $\lambda = \prod_{i = 0}^{k} \lambda_i$.
 Then, using \eqref{eq:TQFT_U2_genus_power} and that $\lambda = 1$, we have
 \begin{gather*}
 \tilde{Z}\big(L^g \circ L_{\mathcal{J}_{\lambda_1}} \circ \cdots \circ L_{\mathcal{J}_{\lambda_k}}\big)(T_{S_1})|_{T_{S_1}}
 = q^{3g + k - 1} (q - 1)^{4g + 2k} \big( (-1)^k (q - 1) + (q - 1)^{2g + k} \big) .
 \end{gather*}
 So finally
 \begin{align*}
 [ \mathfrak{X}_{\U_2}(\Sigma_g, Q) ]
 &= \frac{1}{[\U_2]^{g + k}} \tilde{Z}\big(L^g \circ L_{\mathcal{J}_{\lambda_1}} \circ \cdots \circ L_{\mathcal{J}_{\lambda_k}}\big)(T_{S_1})|_{T_{S_1}} \\
 &= q^{2g - 1} (q - 1)^{2g} \big( (-1)^k (q - 1) + (q - 1)^{2g + k}\big) .
 \end{align*}
 Note that this is in accordance with \eqref{eq:class_representation_variety_U2_genus} for $k = 0$.

(iii) Note that $\prod_{i = 0}^{\ell} \mu_i = \prod_{i = 0}^{\ell} \sigma_i$. In combination with \eqref{eq:TQFT_U2_parabolic_two_eigenvalues} it follows that
 \[ \tilde{Z}\big(L_{\mathcal{M}_{\mu_1, \sigma_1}} \circ \cdots \circ L_{\mathcal{M}_{\mu_\ell, \sigma_\ell}}\big)(T_{S_1}) = q^{2\ell - 1} (q - 1)^{2\ell} (T_{S_{\mu}} + (q - 1) T_{J_{\mu}}) , \]
 where $\mu = \prod_{i = 0}^{\ell} \mu_i$. Similar as before, we use \eqref{eq:TQFT_U2_parabolic_jordan} to obtain
 \begin{align*}
 [ \mathfrak{X}_{\U_2}&(\Sigma_g, Q) ] \\
 &= \frac{1}{[\U_2]^{g + k + \ell}} \tilde{Z}(L^g \circ L_{\mathcal{J}_{\lambda_1}} \circ \cdots \circ L_{\mathcal{J}_{\lambda_k}} \circ L_{\mathcal{M}_{\mu_1, \sigma_1}} \circ \cdots \circ L_{\mathcal{M}_{\mu_\ell, \sigma_\ell}})(T_{S_1})|_{T_{S_1}} \\
 &= q^{2g + \ell - 1} (q - 1)^{4g + k} .\tag*{\qed}
 \end{align*}\renewcommand{\qed}{}
\end{proof}

\subsection[Upper triangular 3 times 3 matrices]{Upper triangular $\boldsymbol{3 \times 3}$ matrices}\label{sec:application_U3}

Now consider the case where $G = \U_3$, the group of upper triangular $3 \times 3$ matrices, that is,
\[ \U_3 = \left\{ \begin{pmatrix} a & b & c \\ 0 & d & e \\ 0 & 0 & f \end{pmatrix} \bigg|\, a, d, f \ne 0 \right\} . \]
It is easy to see that the class of $\U_3$ in the Grothendieck ring of varieties is $\big[\AA^1_\CC\big]^3\big[\AA^1_\CC \setminus \{ 0 \}\big]^3 = q^3(q-1)^3$.

For simplicity we will just consider the representation varieties without parabolic data. As any commutator $[g_1, g_2]$ in $\U_3$ has ones on the diagonal, we only need to consider the conjugacy classes of such elements in order to compute the classes of representation varieties $\mathfrak{X}_{\U_3}(\Sigma_g)$.
There are five such conjugacy classes are given by
\begin{gather*}
\mathcal{C}_1 = \left\{ \begin{pmatrix} 1 & 0 & 0 \\ 0 & 1 & 0 \\ 0 & 0 & 1 \end{pmatrix} \right\} , \qquad
\mathcal{C}_2 = \left\{ \begin{pmatrix} 1 & \alpha & \beta \\ 0 & 1 & \gamma \\ 0 & 0 & 1 \end{pmatrix} \bigg|\, \alpha, \gamma \ne 0 \right\} , \\
\mathcal{C}_3 = \left\{ \begin{pmatrix} 1 & \alpha & \beta \\ 0 & 1 & 0 \\ 0 & 0 & 1 \end{pmatrix} \bigg|\, \alpha \ne 0 \right\} , \qquad
\mathcal{C}_4 = \left\{ \begin{pmatrix} 1 & 0 & \beta \\ 0 & 1 & \alpha \\ 0 & 0 & 1 \end{pmatrix} \bigg|\, \alpha \ne 0 \right\} , \\
\mathcal{C}_5 = \left\{ \begin{pmatrix} 1 & 0 & \alpha \\ 0 & 1 & 0 \\ 0 & 0 & 1 \end{pmatrix}\bigg|\, \alpha \ne 0 \right\}
\end{gather*}
with representatives given by
\begin{gather*}
 \xi_1 = \begin{pmatrix} 1 & 0 & 0 \\ 0 & 1 & 0 \\ 0 & 0 & 1 \end{pmatrix}, \qquad \xi_2 = \begin{pmatrix} 1 & 1 & 0 \\ 0 & 1 & 1 \\ 0 & 0 & 1 \end{pmatrix}, \hspace{0.5em} \xi_3 = \begin{pmatrix} 1 & 1 & 0 \\ 0 & 1 & 0 \\ 0 & 0 & 1 \end{pmatrix}, \\
 \xi_4 = \begin{pmatrix} 1 & 0 & 0 \\ 0 & 1 & 1 \\ 0 & 0 & 1 \end{pmatrix}, \qquad \xi_5 = \begin{pmatrix} 1 & 0 & 1 \\ 0 & 1 & 0 \\ 0 & 0 & 1 \end{pmatrix} .
\end{gather*}
Hence, we obtain maps $\pi_i\colon \mathcal{C}_i \to C_i$ mapping each conjugacy class to the space of orbits under conjugaction. In this case, there is only one orbit for each conjugacy class, hence all the $C_i$ are points. Technically one should stratify $G \backslash \cup_i \mathcal{C}_i$ as well, but as everything happens over $\cup_i \mathcal{C}_i$ we omit this. To obtain a reduction of the TQFT for $G = \U_3$ with the above stratification, let us show that condition~(iv) of Proposition~\ref{prop:conditions_for_reduction_TQFT} holds.

\begin{Lemma}
The stabilizer subgroup of any $g \in G$ is isomorphic to
\begin{gather*}
 \Stab_1 = \U_3, \qquad
\Stab_2 = \left\{ \begin{pmatrix} x & y & z \\ 0 & x & y \\ 0 & 0 & x \end{pmatrix} \bigg|\, x \ne 0 \right\}, \\
\Stab_3 = \left\{ \begin{pmatrix} x & y & z \\ 0 & x & 0 \\ 0 & 0 & w \end{pmatrix} \bigg|\, x, w \ne 0 \right\}, \qquad \Stab_4 = \left\{ \begin{pmatrix} x & 0 & z \\ 0 & w & y \\ 0 & 0 & w \end{pmatrix} \bigg|\, x, w \ne 0 \right\} \qquad \textup{or} \\
\Stab_5 = \left\{ \begin{pmatrix} x & y & z \\ 0 & w & v \\ 0 & 0 & x \end{pmatrix} \bigg|\, x, w \ne 0 \right\} ,
\end{gather*}
depending on whether $g$ lies in $\mathcal{C}_1, \mathcal{C}_2, \dots$, or $\mathcal{C}_5$, respectively. Moreover, all of these subgroups are special.
\end{Lemma}
\begin{proof}
 The proof is straightforward. Note that all stabilizers are extensions of copies of $\GG_m$ and $\GG_a$, which are special.
\end{proof}

We write $T_i = [ C_i ]_{C_i} \in \K(\Var/C_i)$, and consider $V = \langle T_1, \dots, T_5 \rangle$. In what follows, all matrices and vectors will be written with respect to the basis $\{ T_1, \dots, T_5 \}$.

Since, all the fibrations $\pi:\mathcal{C}_i\to C_i$ are trivial, the map $\eta = \pi_! \pi^*$ is simply given by
\[ \eta = \begin{pmatrix} [ \mathcal{C}_1 ] & & & & \\ & [ \mathcal{C}_2 ] & & & \\ & & [ \mathcal{C}_3 ] & & \\ & & & [ \mathcal{C}_4 ] & \\ & & & & [ \mathcal{C}_5 ] \end{pmatrix} = \begin{pmatrix} 1 & & & & \\ & q (q - 1)^2 & & & \\ & & q (q - 1) & & \\ & & & q (q - 1) & \\ & & & & q - 1 \end{pmatrix} . \]

Now we compute $Z_\pi(L) = \pi_! \circ Z(L) \circ \pi^*$, starting with $Z_\pi(L)(T_1)$. Since the commutator $[g_1, g_2]$ has ones on the diagonal for all $g_1, g_2 \in \U_3$, indeed we have that $Z_\pi(L)(T_1) \in \langle T_1, \dots, T_5 \rangle$. We write $g_i = \left( \begin{smallmatrix} a_i & b_i & c_i \\ 0 & d_i & e_i \\ 0 & 0 & f_i \end{smallmatrix} \right) $.

We have that $Z_\pi(L)(T_1)|_{T_1}$ is the class of $\big\{ (g_1, g_2, h) \in \U_3^3 \mid g_1 g_2 = g_2 g_1 \big\}$. Since $g_1g_2$ and $g_2g_1$ have the same elements on the diagonal for every pair of group elements $g_1$ and $g_2$, we only need to check on three entries whether $g_1g_2=g_2g_1$. This gives us three equations in the entries of~$g_1$ and~$g_2$. Explicitly, we obtain
 \begin{gather*}
 a_{1} b_{2} - a_{2} b_{1} + b_{1} d_{2} - b_{2} d_{1} = 0, \\
 a_{1} c_{2} - a_{2} c_{1} + b_{1} e_{2} - b_{2} e_{1} + c_{1} f_{2} - c_{2} f_{1} = 0, \\
 d_{1} e_{2} - d_{2} e_{1} + e_{1} f_{2} - e_{2} f_{1} = 0 .
 \end{gather*}

 \begin{Lemma}\label{lem:disgust}
 The variety in the affine space $\AA_\CC^{12}$ $($with coordinates being the $a_i$, $b_i$, $c_i$, $d_i$, $e_i$, $f_i$ for $i=1$ and $i=2)$ cut out by the three equations above has class $q^3 (q - 1)^4 \big(q^2 + q - 1\big)$ in the Grothendieck ring of varieties $\K(\Var_\CC)$.
 \end{Lemma}

 \begin{proof}
 We cut the variety in pieces as follows.
 \begin{itemize}\itemsep=0pt
 \item $a_1\ne d_1$, $a_1\ne f_1$, $d_1\ne f_1$: In this case, the first equation can be solved for $b_2$, the second for $c_2$ and the third for $e_2$ yielding that the class of this piece is $(q-1)(q-2)(q-3)(q-1)^3q^3$.
 \item $a_1=d_1$, $a_1\ne f_1$, $a_2\ne d_2$: In this case, the first equation can be solved for $b_1$, the second for $c_2$ and the third for $e_2$ yielding that the class of this piece is $(q-1)(q-2)(q-1)^2(q-2)q^3$.
 \item $a_1=d_1$, $a_1\ne f_1$, $a_2= d_2$: In this case, the first equation is always satisfied, moreover, the second can be solved for $c_2$ and the third for $e_2$ yielding that the class of this piece is $(q-1)(q-2)(q-1)^2q^4$.
 \item $a_1\ne d_1$, $d_1=f_1$, $d_2\ne f_2$: In this case, the first equation can be solved for $b_2$, the second for $c_2$ and the third for $e_1$ yielding that the class of this piece is $(q-1)(q-2)(q-1)^2(q-2)q^3$.
 \item $a_1\ne d_1$, $d_1=f_1$, $d_2=f_2$: In this case, the first equation can be solved for $b_2$, the second for $c_2$ and the third equation is always satisfied yielding that the class of this piece is $(q-1)(q-2)(q-1)^2q^4$.
 \item $a_1\ne d_1$, $a_1=f_1$, $a_2\ne f_2$: In this case, the first equation can be solved for $b_2$, the second for $c_1$ and the third for $e_2$ yielding that the class of this piece is $(q-1)(q-2)(q-1)^2(q-2)q^3$.
 \item $a_1\ne d_1$, $a_1=f_1$, $a_2= f_2$: In this case, the first equation can be solved for $b_2$, the third for $e_2$, and the second is then satisfied yielding that the class of this piece is $(q-1)(q-2)(q-1)^2q^4$.
 \item $a_1=d_1=f_1$: In this case, we separate again into cases:
 \begin{itemize}\itemsep=0pt
 \item $b_1=e_1=0$, then $c_1(a_2-f_2)=0$ yielding that the class is $2(q-1)^4q^3$
 \item $b_1=0$, $e_1\ne 0$, then $d_2=f_2$ and the second equation can be solved for $b_2$ yielding that the class is $(q-1)^3(q-1)q^3$,
 \item $b_1\ne 0$, $e_1= 0$, then $a_2=d_2$ and the second equation can be solved for $e_2$ yielding that the class is $(q-1)^3(q-1)q^3$,
 \item $b_1\ne 0$, $e_1\ne 0$, then $a_2=d_2=f_2$ and the second equation can be solved for $e_2$ yielding that the class is
 $(q-1)^2(q-1)^2q^3$.
 \end{itemize}
 So in total, we have that the class is $(q-1)^3q^3(5q-5)$.
 \end{itemize}
 Adding the pieces together we obtain that the class of this variety is $q^3 (q - 1)^4 (q^2 + q - 1)$.
 \end{proof}
 As a consequence, we obtain that class of $\big\{ (g_1, g_2, h) \in \U_3^3 \mid g_1 g_2 = g_2 g_1 \big\}$ is $q^3 (q - 1)^4 (q^2 + q - 1) [\U_3]$ in $\K(\Var_\CC)$.

 \begin{Remark}
 In general, in the computations of $Z_\pi(L)(T_1)$, we follow the same strategy as we explained in the proof above. We cut the variety into pieces given by some variables being~0 or some variables being equal to each other. This idea can be made into an Algorithm~\ref{alg:computing_classes} (see Appendix~\ref{app:A}) which we use to compute the other classes. We added the proof of the above lemma for sake of completeness.
 \end{Remark}

 We have that $Z_\pi(L)(T_1)|_{T_2}$ is the class of $\big\{ (g_1, g_2, h) \in \U_3^3 \mid [ g_1, g_2 ] \in \mathcal{C}_2 \big\}$, which is given by the equations
 \begin{gather*}
 a_{1} b_{2} - a_{2} b_{1} + b_{1} d_{2} - b_{2} d_{1} \ne 0, \\
 d_{1} e_{2} - d_{2} e_{1} + e_{1} f_{2} - e_{2} f_{1} \ne 0 .
 \end{gather*}
 This evaluates to $q^6 (q - 2)^2 (q - 1)^4 [\U_3]$ using Algorithm~\ref{alg:computing_classes} or a similar calculation as in Lemma~\ref{lem:disgust}.

We have that $Z_\pi(L)(T_1)|_{T_3}$ is the class of $\big\{ (g_1, g_2, h) \in \U_3^3 \mid [ g_1, g_2 ] \in \mathcal{C}_3 \big\}$, which is given by the equations
 \begin{gather*}
 a_{1} b_{2} - a_{2} b_{1} + b_{1} d_{2} - b_{2} d_{1} \ne 0, \\
 d_{1} e_{2} - d_{2} e_{1} + e_{1} f_{2} - e_{2} f_{1} = 0.
 \end{gather*}
 This evaluates to $q^6 (q - 2) (q - 1)^4 [\U_3]$ using Algorithm~\ref{alg:computing_classes} or a~similar calculation as in Lemma~\ref{lem:disgust}.

We have that $Z_\pi(L)(T_1)|_{T_4}$ is the class of $\big\{ (g_1, g_2, h) \in \U_3^3 \mid [ g_1, g_2 ] \in \mathcal{C}_4 \big\}$, which is given by the equations
 \begin{gather*}
 a_{1} b_{2} - a_{2} b_{1} + b_{1} d_{2} - b_{2} d_{1} = 0, \\
 d_{1} e_{2} - d_{2} e_{1} + e_{1} f_{2} - e_{2} f_{1} \ne 0.
 \end{gather*}
 This is symmetric to the previous case, so it also evaluates to $q^6 (q - 2) (q - 1)^4 [\U_3]$.

 We have that $Z_\pi(L)(T_1)|_{T_5}$ is the class of $\big\{ (g_1, g_2, h) \in \U_3^3 \mid [ g_1, g_2 ] \in \mathcal{C}_5 \big\}$, which is given by the equations
 \begin{gather*}
 a_{1} b_{2} - a_{2} b_{1} + b_{1} d_{2} - b_{2} d_{1} = 0, \\
 d_{1} e_{2} - d_{2} e_{1} + e_{1} f_{2} - e_{2} f_{1} = 0, \\
 - a_{1} b_{2} d_{2} e_{1} - a_{1} b_{2} e_{2} f_{1} + a_{1} c_{2} d_{1} d_{2} + a_{2} b_{1} d_{2} e_{1} + a_{2} b_{1} e_{2} f_{1} - a_{2} c_{1} d_{1} d_{2} \\
 \qquad{} + b_{1} d_{1} d_{2} e_{2} - b_{1} d_{2}^{2} e_{1} - b_{1} d_{2} e_{2} f_{1} + b_{2} d_{1} e_{2} f_{1} + c_{1} d_{1} d_{2} f_{2} - c_{2} d_{1} d_{2} f_{1} \ne 0.
 \end{gather*}
 The latter inequality can be simplified to
 \[a_{1} c_{2} - a_{2} c_{1} + b_{1} e_{2} - b_{2} e_{1} + c_{1} f_{2} - c_{2} f_{1} \ne 0.\]
 This evaluates to $q^3 (q - 1)^6 (q + 1) [\U_3]$ using Algorithm~\ref{alg:computing_classes} or a similar calculation as in Lemma~\ref{lem:disgust}.

So far we have computed the first column of the matrix of $Z_\pi(L)$. As a check, indeed we have that the sum of the entries of this column equals $[\U_3]$:
\begin{gather*}
 q^3 (q - 1)^4 \big(q^2 + q - 1\big) [\U_3] + q^6 (q - 2)^2 (q - 1)^4 [\U_3] + q^6 (q - 2) (q - 1)^4 [\U_3] \\
 \qquad{} + q^6 (q - 2) (q - 1)^4 [\U_3] + q^3 (q - 1)^6 (q + 1) [\U_3] = [\U_3]^3 .
\end{gather*}

We use a similar strategy as in the previous section to determine $Z_\pi(L)(T_i)$ for $i = 2, 3, 4, 5$ from the case $i = 1$.
We have
\[ Z_\pi(L)(T_j)|_{T_i} = [X_{ij}] \cdot [ \U_3 ] \qquad \text{with} \quad X_{ij} = \big\{ (g, g_1, g_2) \in \mathcal{C}_j \times \U_3^2 \mid g [g_1, g_2] \in \mathcal{C}_i \big\} . \]
We can stratify $X_{ij}$ by
\[ X_{ijk} = \big\{ (g, g_1, g_2) \in \mathcal{C}_j \times \U_3^2 \mid g [g_1, g_2] \in \mathcal{C}_i \text{ and } [ g_1, g_2 ] \in \mathcal{C}_k \big\} \qquad \text{for } k = 1, \dots, 5 . \]
Note that for each conjugacy class $\mathcal{C}_k$ we have an algebraic normal subgroup $N$ so that every element $g\in \mathcal{C}_k$ is given by $n\xi_k n^{-1}$. In other words, there exists a (non-unique) morphism of varieties $\sigma_k\colon \mathcal{C}_k \to \U_3$ such that $\sigma_k(g) \xi_k \sigma_k(g)^{-1} = g$ for all $g\in \mathcal{C}_k$. Moreover, if $[g_1,g_2]=n\xi_k n^{-1}$, then $g [g_1, g_2] \in \mathcal{C}_i$ if and only if $n^{-1}gn\xi\in \mathcal{C}_i$. Thus, for each $i,j,k$ we have an isomorphism of varieties
 \begin{gather*}
 X_{ijk} \xrightarrow{\sim} \{ g \in \mathcal{C}_j \mid g \xi_k \in \mathcal{C}_i \} \times \big\{ (g_1, g_2) \in G^2 \mid [g_1, g_2] \in \mathcal{C}_k \big\}, \\
 (g, g_1, g_2) \mapsto (\sigma_k([g_1, g_2])^{-1} g \sigma_k([g_1, g_2]), g_1, g_2),
\end{gather*}
so we find that
\begin{equation}
 \label{eq:Z_pi_L_from_first_row}
 Z_\pi(L)(T_j)|_{T_i} = \sum_{k = 1}^{5} F_{ijk} \cdot Z_\pi(L)(T_1)|_{T_k} \qquad \text{with} \quad F_{ijk} = [ \{ g \in \mathcal{C}_j \mid g \xi_k \in \mathcal{C}_i \} ] .
\end{equation}
Although there are about $5^3 = 125$ computations to be done to determine the coefficients $F_{ijk}$, all of them are quite simple. For instance, it is clear that $F_{i,1,k} = \delta_{ik}$, the Kronecker delta. For $j = 2$, take any $g = \left(\begin{smallmatrix} 1 & \alpha & \beta \\ 0 & 1 & \gamma \\ 0 & 0 & 1 \end{smallmatrix} \right) \in \mathcal{C}_2$. Then $g \xi_1, g \xi_5 \in \mathcal{C}_2$. We have $g \xi_3 \in \mathcal{C}_4$ if $\alpha = -1$ and $g \xi_3 \in \mathcal{C}_2$ otherwise. Similarly, $g \xi_4 \in \mathcal{C}_3$ if $\gamma = -1$ and $g \xi_4 \in \mathcal{C}_2$ otherwise. Finally,
\[ g \xi_2 \in \begin{cases}
 \mathcal{C}_1 & \text{if } \alpha, \gamma = -1 \text{ and } \beta = 0 , \\
 \mathcal{C}_2 & \text{if } \alpha, \gamma \ne -1, \\
 \mathcal{C}_3 & \text{if } \alpha \ne -1 \text{ and } \gamma = -1 , \\
 \mathcal{C}_4 & \text{if } \alpha = -1 \text{ and } \gamma \ne -1 , \\
 \mathcal{C}_5 & \text{if } \alpha, \gamma = -1 \text{ and } \beta \ne 0 .
 \end{cases} \]
This gives
\[ F_{i,2,k} = \begin{pmatrix} 0 & 1 & 0 & 0 & 0\\q (q - 1)^{2} & q (q - 2)^{2} & q (q - 2) (q - 1) & q (q - 2) (q - 1) & q (q - 1)^{2}\\0 & q (q - 2) & 0 & q (q - 1) & 0\\0 & q (q - 2) & q (q - 1) & 0 & 0\\0 & q - 1 & 0 & 0 & 0 \end{pmatrix} , \]
with $i$ the row index and $k$ the column index. By completely similar arguments, one can show that
\begin{gather*}
 F_{i,3,k} = \begin{pmatrix} 0 & 0 & 1 & 0 & 0\\0 & q (q - 2) & 0 & q (q - 1) & 0\\q (q - 1) & 0 & q (q - 2) & 0 & q (q - 1)\\0 & q & 0 & 0 & 0\\0 & 0 & q - 1 & 0 & 0 \end{pmatrix} , \\ F_{i,4,k} = \begin{pmatrix} 0 & 0 & 0 & 1 & 0\\0 & q (q - 2) & q (q - 1) & 0 & 0\\0 & q & 0 & 0 & 0\\q (q - 1) & 0 & 0 & q (q - 2) & q (q - 1)\\0 & 0 & 0 & q - 1 & 0 \end{pmatrix} ,\\
 F_{i,5,k} = \begin{pmatrix} 0 & 0 & 0 & 0 & 1\\0 & q - 1 & 0 & 0 & 0\\0 & 0 & q - 1 & 0 & 0\\0 & 0 & 0 & q - 1 & 0\\q - 1 & 0 & 0 & 0 & q - 2 \end{pmatrix} .
\end{gather*}
Using \eqref{eq:Z_pi_L_from_first_row} we obtain the matrix representing $Z_\pi(L)$. Precomposing this map with $\eta^{-1}$ gives $\tilde{Z}(L)$, which can be diagonalized as follows:
\begin{gather*}
 \tilde{Z}(L) = q^6 (q - 1)^5 A \begin{pmatrix} q^3 &&&& \\ & q (q - 1)^2 &&& \\ && q^3 (q - 1)^2 && \\ &&& q^3 (q - 1)^2 & \\ &&&& q^3 (q - 1)^4 \end{pmatrix} A^{-1} \\
\hphantom{\tilde{Z}(L) =}{} \text{with} \quad A = \begin{pmatrix} 1 & 1 & 0 & 1 & 1 \\ q & 0 & 0 & -q(q - 1) & q(q - 1)^2 \\ -q & 0 & 1 & 0 & q (q - 1) \\ -q & 0 & -1 & q (q - 2) & q (q - 1) \\ q - 1 & -1 & 0 & q - 1 & q - 1 \end{pmatrix} .
\end{gather*}

\begin{Remark}
 We see that $Z_\pi(L)$ is symmetric. This can be explained by the fact that for each~$i$ and~$j$ we have an isomorphism:
 \begin{gather*}
 \big\{ (g, g_1, g_2) \in \mathcal{C}_i \times G^2 \mid g [g_1, g_2] \in \mathcal{C}_j \big\} \overset{\sim}{\longleftrightarrow} \big\{ (g, g_1, g_2) \in \mathcal{C}_j \times G^2 \mid g [g_1, g_2] \in \mathcal{C}_i \big\}, \\
 (g, g_1, g_2) \mapsto (g [g_1, g_2], g_2, g_1), \\
 (g [g_1, g_2], g_2, g_1) \mapsfrom (g, g_1, g_2).
 \end{gather*}
 Hence the classes of both sides are equal in $\K(\Var_\CC)$, and so $Z_\pi(L)(T_i)|_{T_j} = Z_\pi(L)(T_j)|_{T_i}$. This does not hold for any TQFT, but it relies on the fact that the $C_i$ are points.
\end{Remark}

\begin{Theorem} \label{thm:result_U3}
 For any $g \ge 0$, the virtual class of the $\U_3$-representation variety $\mathfrak{X}_{\U_3}(\Sigma_g)$ is
 \begin{gather*}
 [ \mathfrak{X}_{\U_3}(\Sigma_g) ] = q^{3g - 3} (q - 1)^{2g} \big( q^2 (q - 1)^{2g + 1} + q^{3g} (q - 1)^2 + q^{3g} (q - 1)^{4g} + 2 q^{3g} (q - 1)^{2g + 1} \big) .
 \end{gather*}
\end{Theorem}

\begin{proof}
 One can check that
 \[ A^{-1} = \frac{1}{q^3} \begin{pmatrix} (q - 1)^2 & 1 & 1 - q & 1 - q & (q - 1)^2 \\ q^2 (q - 1) & 0 & 0 & 0 & -q^2 \\ q(q - 2)(q - 1) & -q (q - 2) & q^3 - 2q^2 + 2q & -2q(q - 1) & q^3 - 3q^2 + 2q \\ 2q - 2 & -2 & q - 2 & q - 2 & 2q - 2 \\ 1 & 1 & 1 & 1 & 1 \end{pmatrix} . \]
 By matrix multiplication, we find that
 \begin{gather*}
 \mathfrak{X}_{\U_3}(\Sigma_g)
 = \frac{1}{[\U_3]^g} \tilde{Z}(L)^g(T_1)|_{T_1} = q^{3g - 3} (q - 1)^{2g} \\
 \hphantom{\mathfrak{X}_{\U_3}(\Sigma_g)=}{}
 \times \big( q^2 (q - 1)^{2g + 1} + q^{3g} (q - 1)^2 + q^{3g} (q - 1)^{4g} + 2 q^{3g} (q - 1)^{2g + 1} \big) .\tag*{\qed}
 \end{gather*}\renewcommand{\qed}{}
\end{proof}

\begin{Remark}
 In particular, for small values of $g$, we find
 \begin{gather*}
 [ \mathfrak{X}_{\U_3}(\Sigma_1) ] = q^3 (q - 1)^4 \big(q^2 + q - 1\big) , \\
 [ \mathfrak{X}_{\U_3}(\Sigma_2) ] = q^{7} (q - 1)^{6} \big(q^{8} - 6 q^{7} + 15 q^{6} - 18 q^{5} + 9 q^{4} + q^{3} - 3 q^{2} + 3 q - 1\big), \\
 [ \mathfrak{X}_{\U_3}(\Sigma_3) ] = q^{11} (q - 1)^{8} \big(q^{14} - 10 q^{13} + 45 q^{12} - 120 q^{11} + 210 q^{10} - 250 q^{9} + 200 q^{8} \\
 \hphantom{[ \mathfrak{X}_{\U_3}(\Sigma_3) ] =}{}
 - 100 q^{7} + 25 q^{6} + q^{5} - 5 q^{4} + 10 q^{3} - 10 q^{2} + 5 q - 1\big) .
 \end{gather*}
 As in Remark \ref{rem:result_genus_U2}, the factor $(q - 1)^{2g + 2}$ can be explained from the actions of $\GG_m^{2g}$ (given by scaling the $A_i$, $B_i$) and $\GG_m^2$ (given by conjugating with $\left(\begin{smallmatrix} 1 & & \\ & x & \\ & & y \end{smallmatrix} \right)$, $x, y \in \CC^*$).
\end{Remark}

\subsection[Upper triangular 4 times 4 matrices]{Upper triangular $\boldsymbol{4 \times 4}$ matrices}
\label{sec:application_U4}
The last case we will treat is the group $\U_4$ of upper triangular $4 \times 4$ matrices. We can use the same strategies as in the previous case of $\U_3$, but all computations are done using Algorithm~\ref{alg:computing_classes}. Source code for these computations is given in~\cite{github}. The group $\U_4$ contains sixteen unipotent conjugacy classes \cite{bhunia}. We consider the following representatives of these classes
\begin{gather*}
 \begin{pmatrix} 1 & & & \\ & 1 & & \\ & & 1 & \\ & & & 1 \end{pmatrix} ,\quad
 \begin{pmatrix} 1 & 1 & & \\ & 1 & & \\ & & 1 & \\ & & & 1 \end{pmatrix} ,\quad
 \begin{pmatrix} 1 & & 1 & \\ & 1 & & \\ & & 1 & \\ & & & 1 \end{pmatrix} ,\quad
 \begin{pmatrix} 1 & & & 1 \\ & 1 & & \\ & & 1 & \\ & & & 1 \end{pmatrix} ,\\
 \begin{pmatrix} 1 & & & \\ & 1 & 1 & \\ & & 1 & \\ & & & 1 \end{pmatrix} ,\quad
 \begin{pmatrix} 1 & & & \\ & 1 & & 1 \\ & & 1 & \\ & & & 1 \end{pmatrix} ,
\quad
 \begin{pmatrix} 1 & & & \\ & 1 & & \\ & & 1 & 1 \\ & & & 1 \end{pmatrix} ,\quad
 \begin{pmatrix} 1 & 1 & & \\ & 1 & 1 & \\ & & 1 & \\ & & & 1 \end{pmatrix} ,\\
 \begin{pmatrix} 1 & 1 & & \\ & 1 & & 1 \\ & & 1 & \\ & & & 1 \end{pmatrix} ,\quad
 \begin{pmatrix} 1 & 1 & & \\ & 1 & & \\ & & 1 & 1 \\ & & & 1 \end{pmatrix} ,\quad
 \begin{pmatrix} 1 & & 1 & \\ & 1 & & 1 \\ & & 1 & \\ & & & 1 \end{pmatrix} ,\quad
 \begin{pmatrix} 1 & & 1 & \\ & 1 & & \\ & & 1 & 1 \\ & & & 1 \end{pmatrix} ,\\
 \begin{pmatrix} 1 & & & \\ & 1 & 1 & \\ & & 1 & 1 \\ & & & 1 \end{pmatrix} ,\quad
 \begin{pmatrix} 1 & & & 1 \\ & 1 & 1 & \\ & & 1 & \\ & & & 1 \end{pmatrix} ,\quad
 \begin{pmatrix} 1 & 1 & & \\ & 1 & 1 & \\ & & 1 & 1 \\ & & & 1 \end{pmatrix} ,\quad
 \begin{pmatrix} 1 & 1 & 1 & \\ & 1 & & \\ & & 1 & 1 \\ & & & 1 \end{pmatrix} ,
\end{gather*}
which we denote in order by $\xi_1, \dots, \xi_{16}$. Explicitly, the conjugacy classes are given by
\begin{gather*}
 \mathcal{C}_{1} = \{ a_{0,1} = a_{0,2} = a_{0,3} = a_{1,2} = a_{1,3} = a_{2,3} = 0 \} , \\
 \mathcal{C}_{2} = \{ a_{1,2} = a_{1,3} = a_{2,3} = 0, a_{0,1} \ne 0 \} , \\
 \mathcal{C}_{3}= \{ a_{0,1} = a_{1,2} = a_{1,3} = a_{2,3} = 0, a_{0,2} \ne 0 \} , \\
 \mathcal{C}_{4} = \{ a_{0,1} = a_{0,2} = a_{1,2} = a_{1,3} = a_{2,3} = 0, a_{0,3} \ne 0 \} , \\
 \mathcal{C}_{5} = \{ a_{0,1} = a_{2,3} = a_{0,3} a_{1,2} - a_{0,2} a_{1,3} = 0, a_{1,2} \ne 0 \} , \\
 \mathcal{C}_{6} = \{ a_{0,1} = a_{0,2} = a_{1,2} = a_{2,3} = 0, a_{1,3} \ne 0 \} , \\
 \mathcal{C}_{7} = \{ a_{0,1} = a_{0,2} = a_{1,2} = 0, a_{2,3} \ne 0 \} , \\
 \mathcal{C}_{8}= \{ a_{2,3} = 0, a_{0,1} \ne 0, a_{1,2} \ne 0 \} , \\
 \mathcal{C}_{9} = \{ a_{1,2} = a_{2,3} = 0, a_{0,1} \ne 0, a_{1,3} \ne 0 \} , \\
 \mathcal{C}_{10}= \{ a_{1,2} = a_{0,2} a_{2,3} + a_{0,1} a_{1,3} = 0, a_{0,1} \ne 0, a_{2,3} \ne 0 \} , \\
 \mathcal{C}_{11} = \{ a_{0,1} = a_{1,2} = a_{2,3} = 0, a_{0,2} \ne 0, a_{1,3} \ne 0 \} , \\
 \mathcal{C}_{12} = \{ a_{0,1} = a_{1,2} = 0, a_{0,2} \ne 0, a_{2,3} \ne 0 \} , \\
 \mathcal{C}_{13} = \{ a_{0,1} = 0, a_{1,2} \ne 0, a_{2,3} \ne 0 \} , \\
 \mathcal{C}_{14} = \{ a_{0,1} = a_{2,3} = 0, a_{1,2} \ne 0, a_{0,3} a_{1,2} - a_{0,2} a_{1,3} \ne 0 \} , \\
 \mathcal{C}_{15} = \{ a_{0,1} \ne 0, a_{1,2} \ne 0, a_{2,3} \ne 0 \} , \\
 \mathcal{C}_{16} = \{ a_{1,2} = 0, a_{0,1} \ne 0, a_{2,3} \ne 0, a_{0,2} a_{2,3} + a_{0,1} a_{1,3} \ne 0 \} .
\end{gather*}
As before, we write $T_i = [C_i]_{C_i} \in \K(\Var/C_i)$ and consider $V = \langle T_1, \dots, T_{16} \rangle$. We wish to compute $Z_\pi(L)(T_1)$ first and deduce the other columns from this column as we did before. We have
\[ Z_\pi(L)(T_1)|_{T_i} = [ \{ (g_1, g_2) \in \U_4 \mid [g_1, g_2] \in \mathcal{C}_i \} ]. \]
In terms of coordinates, this will yield systems of equations in 20 variables.
Using Algorithm~\ref{alg:computing_classes} we obtain
\[ Z_\pi(L)(T_1) = q^6 (q - 1)^2 \left(\begin{matrix}q (q - 1)^{3} \big(q^{2} + 3 q - 2\big)\\q (q - 2) (q - 1)^{2}\big(q^{2} + q - 1\big)\\
q (q - 1)^{2} \big(q^{3} - 3 q + 1\big)\\
q (q - 1)^{3} \big(q^{2} + 2 q - 2\big)\\
(q - 2) (q - 1)^{2} \big(q^{3} + q - 1\big)\\
q (q - 1)^{2} \big(q^{3} - 3 q + 1\big)\\
q (q - 2) (q - 1)^{2}\big(q^{2} + q - 1\big)\\
q^{3} (q - 2)^{2} (q - 1)\\q (q - 2) (q - 1)^{3} (q + 1)\\
q (q - 2) (q - 1)^{2} \big(q^{2} - q - 1\big)\\
q (q - 1) \big(q^{4} - 2 q^{3} - q^{2} + 4 q - 1\big)\\
q (q - 2) (q - 1)^{3} (q + 1 )\\q^{3} (q - 2)^{2} (q - 1)\\
(q - 2) (q - 1)^{3} \big(q^{2} + q + 1\big)\\
q^{3} (q - 2)^{3}\\
q (q - 2) \big(q^{4} - 3 q^{3} + 2 q^{2} - 1\big)\end{matrix}\right). \]
Completely analogous to the previous section, the other columns can be computed from this result via
\[ Z_\pi(L)(T_j)|_{T_i} = \sum_{k = 1}^{16} F_{ijk} \cdot Z_\pi(L)(T_1|_{T_k}), \]
where
\[ F_{ijk} = [ \{ g \in \mathcal{C}_j \mid g \xi_k \in \mathcal{C}_i \} ] . \]
Again, see~\cite{github} for the actual computations. As usual, the map $\eta = \pi_! \pi^*$ is diagonal, with $\eta(T_i) = [ \mathcal{C}_i ] T_i$. Composing $Z_\pi(L)$ with $\eta^{-1}$ gives us $\tilde{Z}(L)$, which can be found in Appendix~\ref{app:B}. After diagonalizing the reduced TQFT $\tilde{Z}(L) = Z_\pi(L) \circ \eta^{-1}$, we obtain the following result.

\begin{Theorem} \label{thm:result_U4}
 For any $g \ge 0$, the virtual class of the $\U_4$-representation variety $\mathfrak{X}_{\U_4}(\Sigma_g)$ is
 \begin{gather*}
 [ \mathfrak{X}_{\U_4}(\Sigma_g) ]
 =  q^{12 g - 6} (q - 1)^{8 g} + q^{12 g - 6} (q - 1)^{2 g + 3} + q^{10 g - 4} (q - 1)^{2 g + 3} + q^{10 g - 3} (q - 1)^{4 g + 1} \\
\hphantom{[ \mathfrak{X}_{\U_4}(\Sigma_g) ]=}{} + q^{8 g - 2} (q - 1)^{6 g + 1} + q^{8 g - 2} (q - 1)^{4 g + 2} + 2 q^{10 g - 4} (q - 1)^{6 g + 1} \\
\hphantom{[ \mathfrak{X}_{\U_4}(\Sigma_g) ]=}{} + 3 q^{12 g - 6} (q - 1)^{6 g + 1} + 3 q^{12 g - 6} (q - 1)^{4 g + 2} + q^{10 g - 4} (q - 1)^{4 g + 1} (2 q - 3 ) .
 \end{gather*}
\end{Theorem}

\begin{Remark}
 For small values of $g$, we have
 \begin{gather*}
 [ \mathfrak{X}_{\U_4}(\Sigma_1) ] = q^{7} (q - 1)^{5} \big(q^{2} + 3 q - 2\big), \\
 [ \mathfrak{X}_{\U_4}(\Sigma_2) ] = q^{15} (q - 1)^{7} \big(q^{2} - 3 q + 3\big) \big(q^{10} - 6 q^{9} + 15 q^{8} - 18 q^{7} + 9 q^{6} + 2 q^{5} - 6 q^{4} + 7 q^{3} \\
 \hphantom{[ \mathfrak{X}_{\U_4}(\Sigma_2) ] =}{} - 4 q^{2} + 3 q - 1\big), \\
 [ \mathfrak{X}_{\U_4}(\Sigma_3) ] = q^{23} (q - 1)^{9} \big(q^{4} - 5 q^{3} + 10 q^{2} - 10 q + 5\big) \big(q^{18} - 10 q^{17} + 45 q^{16} - 120 q^{15} \\
 \hphantom{[ \mathfrak{X}_{\U_4}(\Sigma_3) ] =}{}
 + 210 q^{14} - 250 q^{13} + 200 q^{12} - 100 q^{11} + 25 q^{10} + 2 q^{9} - 10 q^{8} + 20 q^{7} - 20 q^{6}
 \\
 \hphantom{[ \mathfrak{X}_{\U_4}(\Sigma_3) ] =}{}
 + 11 q^{5} - 6 q^{4} + 10 q^{3} - 10 q^{2} + 5 q - 1 \big).
 \end{gather*}
\end{Remark}

\subsection{Moduli space of representations and character variety}\label{sec:character_varieties}
Let $X$ be a path-connected topological space with finitely generated fundamental group and $G$ a linear algebraic group over an algebraically closed field $k$. There is a natural action of $G$ on the affine representation variety $\mathfrak{X}_G(X)$ given by conjugation. If $G$ is reductive, one can look at the affine \textit{geometric invariant theory} (GIT) quotient
\[ \mathcal{M}_G(X) = \mathfrak{X}_G(X) \sslash G , \]
which is defined as the spectrum of the ring of invariants $\Spec (\mathcal{O}_{\mathfrak{X}_G(X)})^G$. This scheme is known as the \emph{moduli space of $G$-representations}.

A theorem of Nagata \cite{nagata} shows that the affine GIT quotient is finitely generated over $k$, using that $G$ is reductive. However, specializing to the non-reductive groups $G = \U_n$, the ring of invariants is no longer guaranteed to be finitely generated over $k$. Hence, we will instead focus on the categorical quotient of $\mathfrak{X}_G(X)$ by $G$, which coincides with the affine GIT quotient for $G$ reductive \cite{newstead_intro}. We begin with the following lemma.

\begin{Lemma} \label{lemma:about_categorical_quotient}
 Let $X$ be a variety over an algebraically closed field $k$ equipped with a $G$-action. Let $\pi\colon X \to Y$ be a $G$-equivariant morphism of varieties over $k$, such that the action of $G$ on $Y$ is trivial. Assume that there exists a $G$-equivariant morphism $\sigma\colon Y \to X$ such that $\pi \circ \sigma = \id_Y$. If for any $x \in X$ the Zariski-closure of the $G$-orbit of $x$ contains $\sigma(\pi(x))$, then $\pi$ is a categorical quotient.
\end{Lemma}

\begin{proof}
 Let $f\colon X \to Z$ be any $G$-equivariant morphism, where $Z$ has trivial $G$-action. We need to show there exists a unique $G$-invariant morphism $g\colon Y \to Z$ such that $f = g \circ \pi$. If such a~morphism exists, it must be given by $g = f \circ \sigma$ since $f\circ \sigma=(g\circ \pi)\circ \sigma=g$. This already shows uniqueness.
 Now we show that for $g = f \circ \sigma$ we have $f = g \circ \pi$. Take $x \in X$ and note that as $f$ is $G$-equivariant, so $f(\tilde{x}) = f(x)$ for any $\tilde{x}$ in the orbit of $x$. By continuity, we find that $f(\sigma(\pi(x))) = f(x)$ finishing the proof.
\end{proof}

Let us apply this lemma to the case of $G = \U_2$.

\begin{Lemma} There exists an isomorphism of varieties $\mathcal{M}_{\U_2}(\Sigma_g) \cong \big(\AA^1_\CC \setminus \{ 0 \}\big)^{4g}$ over $\CC$.
\end{Lemma}

\begin{proof}
Consider
\[ M = \big\{ (A_1, B_1, \dots, A_g, B_g) \in \mathfrak{X}_{\U_2}(\Sigma_g) \mid \text{all $A_i, B_i$ are diagonal} \big\} , \]
which is isomorphic to $\big(\AA^1_\CC \setminus \{ 0 \}\big)^{4g}$. Let $\pi\colon \mathfrak{X}_{\U_2}(\Sigma_g) \to M$ be the morphism that sends every matrix $\left(\begin{smallmatrix} a_i & b_i \\ 0 & c_i \end{smallmatrix}\right)$ to the matrix $\left(\begin{smallmatrix} a_i & 0 \\ 0 & c_i \end{smallmatrix}\right)$, and take $\sigma\colon M \to \mathfrak{X}_{\U_2}(\Sigma_g)$ to be the inclusion. Indeed $\pi$ is $\U_2$-invariant, and for any
\[ A = \left( \begin{pmatrix} a_1 & b_1 \\ 0 & c_1 \end{pmatrix}, \dots, \begin{pmatrix} a_{2g} & b_{2g} \\ 0 & c_{2g} \end{pmatrix} \right) \in \mathfrak{X}_{\U_2}(\Sigma_g) \]
we find that
\[ \lim_{x \to 0} \begin{pmatrix} x & 0 \\ 0 & 1 \end{pmatrix} A \begin{pmatrix} x & 0 \\ 0 & 1 \end{pmatrix}^{-1}
= \left( \begin{pmatrix} a_1 & 0 \\ 0 & c_1 \end{pmatrix}, \dots, \begin{pmatrix} a_{2g} & 0 \\ 0 & c_{2g} \end{pmatrix} \right)
= \sigma(\pi(A)) , \]
so $\sigma(\pi(A))$ lies in the analytic-closure (and thus Zariski-closure) of the orbit of $A$. By Lem\-ma~\ref{lemma:about_categorical_quotient}, we conclude that $\pi\colon \mathfrak{X}_{\U_2}(\Sigma_g) \to M$ is the categorical quotient of $\mathfrak{X}_{\U_2}(\Sigma_g)$ by the action of~$\U_2$ providing the required isomorphism $\mathcal{M}_{\U_2}(\Sigma_g) \cong \big(\AA^1_\CC \setminus \{ 0 \}\big)^{4g}$.
\end{proof}

We remark that the moduli space of $G$-representations $\mathcal{M}_{\U_2}(\Sigma_g)$ is a variety even though $\U_2$ is a non-reductive linear group.

Now, we turn our attention to the $G$-character variety of $X$. We write $\Gamma = \pi_1(X)$, which is assumed to be finitely generated. The \emph{character} of a representation $\rho \in \mathfrak{X}_G(X)$ is defined as the map
\[ \chi_\rho\colon \ \Gamma \to k\colon \ \gamma \mapsto \operatorname{tr}(\rho(\gamma)) , \]
and the \emph{character map} as
\[ \chi\colon \ \mathfrak{X}_G(X) \to k^\Gamma\colon \ \rho \mapsto \chi_\rho . \]
The image of $\chi$ is called the \emph{$G$-character variety}, denoted $\chi_G(X)$. By the results from~\cite{culler_and_shalen}, there exists a finite set of elements $\gamma_1, \dots, \gamma_a \in \pi_1(X)$ such that $\chi_\rho$ is determined by the characters $\chi_\rho(\gamma_1), \dots, \chi_\rho(\gamma_a)$ for any $\rho$. This way $\chi_G(X)$ can be identified with the image of the map $\mathfrak{X}_G(X) \to k^a\colon \rho \mapsto (\chi_\rho(\gamma_1), \dots, \chi_\rho(\gamma_a))$, which gives the $G$-character variety the structure of a~variety. This structure is independent of the chosen $\gamma_i$.

Note that the character map $\chi$ is a $G$-invariant morphism: indeed the trace map is invariant under conjugation. By the universal property of the categorical quotient, we obtain an induced map
\[ \overline{\chi}\colon \ \mathcal{M}_G(X) \to \chi_G(X) . \]
In the case of $G = \SL_n(\CC)$, $\Sp_{2n}
(\CC)$ or $\SO_{2n+1}(\CC)$ this map is an isomorphism \cite{culler_and_shalen, flr2017, lawsik2017}. However, $\overline{\chi}$ fails to be an isomorphism for the surfaces $\Sigma_g$ and the non-reductive groups $G = \U_2$ with $g \ge 1$.

\begin{Theorem}The map \[ \overline{\chi}\colon \mathcal{M}_{\U_2}(\Sigma_g) \to \chi_{\U_2}(\Sigma_g) \]
is not an isomorphism.
\end{Theorem}

\begin{proof} Arguing as in Lemma~\ref{lemma:about_categorical_quotient}, the map $\overline{\chi}$
 must be given by $\overline{\chi}(A) = \chi_A$. But this cannot be an isomorphism: for general $A \in \mathcal{M}_{\U_2}(\Sigma_g)$ one can consider $B = \left(\begin{smallmatrix} 0 & 1 \\ 1 & 0 \end{smallmatrix}\right) A \left(\begin{smallmatrix} 0 & 1 \\ 1 & 0 \end{smallmatrix}\right)^{-1}$ (where the diagonal entries of $A$ are interchanged), and we have $\chi_A = \chi_B$, even though in general $A \ne B$. Therefore, the moduli space $\mathcal{M}_{\U_2}(\Sigma_g)$ is not isomorphic to the character variety $\chi_{\U_2}(\Sigma_g)$ through the natural map.
\end{proof}

All of the above can easily be generalized to the case $G = \U_n$ for any $n \ge 2$. Namely, similar to before, let
\[ M = \big\{ ( A_1, B_1, \dots, A_{g}, B_{g} ) \in \mathfrak{X}_{\U_n}(\Sigma_g) \mid \text{all $A_i, B_i$ are diagonal} \big\} , \]
and $\pi\colon \mathfrak{X}_{\U_n}(\Sigma_g) \to M$ the map that sets all off-diagonal entries to zero, and $\sigma\colon M \to \mathfrak{X}_{\U_n}(\Sigma_g)$ the inclusion. Then $\sigma \circ \pi = \id$, and one easily checks that $\sigma(\pi(A))$ lies in the closure of the orbit of $A$ for any $A \in \mathfrak{X}_{\U_n}(\Sigma_g)$, e.g., by conjugating $A$ with $\left(\begin{smallmatrix} x^{n - 1} & & & \\ & x^{n - 2} & & \\ & & \ddots & \\ & & & 1 \end{smallmatrix}\right)$ and taking the limit $x \to 0$. This shows that $\mathcal{M}_{\U_n}(\Sigma_g) \cong \big(\AA^1_\CC \setminus \{ 0 \}\big)^{2gn}$, proving Theorem \ref{thm:moduli_space_Un_Sg}. Again, note that the natural map $\overline{\chi}\colon \mathcal{M}_{\U_n}(\Sigma_g) \to \chi_{\U_n}(\Sigma_g)$ cannot be an isomorphism as there are symmetries (permuting diagonal entries) that are invariant under the character map, proving Theorem \ref{thm:main2}.

\begin{Remark}
In fact, as pointed out by the reviewers, the discussion above implies that the $\U_n$-character variety $\chi_{\U_n}(\Sigma_g)$ is the GIT quotient of the moduli space of $\U_n$-representations $\mathcal{M}_{\U_n}(\Sigma_g)$ under the symmetric group $S_n$,
\[ \chi_{\U_n}(\Sigma_g) \cong \mathcal{M}_{\U_n}(\Sigma_g) \sslash S_n = \big(\AA^1_\CC \setminus \{ 0 \}\big)^{2gn} \sslash S_n, \]
where $S_n$ permutes the eigenvalues of the $2g$ generators simultaneously. Indeed, the traces of the products of the powers of the $2g$ generators determine the eigenvalues of the $2g$ generators up to a simultaneous $S_n$-action, which in turn, determine the character of any representation $\chi_\rho\colon \Gamma \to k$.

In this way, the $\U_n$-character variety $\chi_{\U_n}(\Sigma_g)$ can be understood geometrically. The variety $\big(\AA^1_{\CC}\setminus \{0\}\big)^{n}$ is identified with variety of diagonal matrices of $\GL_n(\CC)$, and thus, with a maximal torus $T$. The symmetric group $S_n$ can be identified with the Weyl group $W$ acting on the maximal torus, and thus, one gets isomorphisms
\[\chi_{\U_n}(\Sigma_g) \cong \big(\AA^1_{\CC}\setminus \{ 0 \}\big)^{2gn} \sslash S_n \cong T^{2g} \sslash S_n \cong \Hom\big(\ZZ^{2g}
, T\big) \sslash W.\]
Now the map
\[T^{2g} \to \Hom\big(\ZZ^{2g}, T\big) \to \Hom\big(\ZZ^{2g}, \GL_n(\CC)\big) \to \Hom\big(\ZZ^{2g}, \GL_n(\CC)\big) \sslash \GL_n(\CC)\]
factors through $T^{2g} \sslash W$, so we obtain a map from the $\U_n$-character variety to the moduli space of $\GL_n(\CC)$-representations of the free abelian group $\ZZ^{2g}$ (which is the abelianization of $\pi_1(\Sigma_g)$)
\[ \chi_{\U_n}(\Sigma_g) \cong \Hom\big(\ZZ^{2g}
, T\big) \sslash W \to \Hom\big(\ZZ^{2g}, \GL_n(\CC)\big) \sslash \GL_n(\CC) , \]
which is, in fact, an isomorphism, see~\cite{florentino2021hodge,sikora2014character} for more details.
\end{Remark}

\appendix

\section{Algorithmically computing virtual classes}
\label{app:A}
In this appendix, we describe our algorthim for computing classes of affine varieties over $\CC$ in the Grothendieck ring of varieties $\K(\Var_\CC)$ in terms of $q = \big[ \AA^1_\CC \big]$.

Let $S = \{ x_1, \dots, x_n \}$ be a finite set (of variables), and $F$, $G$ be finite subsets of $\CC[S]$. Then we write $X(S, F, G)$ for the (reduced) subvariety of $\AA_\CC^n$ given by $f = 0$ for all $f \in F$ and $g \ne 0$ for all $g \in G$. For example,
\[ \AA_\CC^n = X(\{ x_1, \dots, x_n \}, \varnothing, \varnothing) \qquad \text{and} \qquad \GL_2(\CC) = X(\{ a, b, c, d \}, \varnothing, \{ ad - bc \}) . \]
For convenience we will write $\text{ev}_x(f, u)$ with $f \in \CC[S]$ for the polynomial where $x \in S$ in $f$ is substituted for $u \in \CC[S]$. Then we write
\[ \text{ev}_x(F, u) = \{ \text{ev}_x(f, u) \mid f \in F \} \]
and
\[ \text{ev}_x(F, u, v) = \big\{ v^{\deg_x(f)} \cdot \text{ev}_x(f, u/v) \mid f \in F \big\}\]
for a set of polynomials $F$ and $x \in S$ and $u, v \in \CC[S]$. Note that for $\text{ev}_x(F, u, v)$ the substituted polynomials are multiplied by a suitable number of factors $v$, in order to clear denominators.

We describe our recursive algorithm which computes the class $[X(S,F,G)]$ in terms of ${q=\big[\AA^1_\CC\big]}$.
\begin{Algorithm} \label{alg:computing_classes}
 Let $X = X(S, F, G)$ for some $S$, $F$ and $G$ as above.
 \begin{enumerate}\itemsep=0pt
 \item If $F$ contains a non-zero constant or if $0 \in G$, then $X = \varnothing$, so $[X] = 0$.

 \item If $F = \varnothing$ and $G = \varnothing$, then $X = \AA_\CC^{\# S}$, so $[X] = q^{\# S}$.

 \item If some $x \in S$ `does not appear' in any $f \in F$ and any $g \in G$, then we can factor $X \simeq \AA^1_\CC \times X'$ with $X' = X(S \setminus \{ x \}, F, G)$. We have $[X] = q [X']$.

 \item If $f = u^n$ (with $n > 1$) for some $f \in F$ and $u \in \CC[S]$, then we can replace $f$ with $u$, not changing $X$. That is, $X = X(S, (F \setminus \{ f \}) \cup \{ u \}, G)$. Similarly, if $g = u^n$ (with $n > 1$) for some $g \in G$ and $u \in \CC[S]$, then $X = X(S, F, (G \setminus \{ g \}) \cup \{ u \})$.

 \item If some $f \in F$ is univariate in $x \in S$, we write $f = (x - \alpha_1) \cdots (x - \alpha_m)$, and we have $[X] = \sum_{i = 1}^{m} [X_i]$ with $X_i = X(S \setminus \{ x \}, \text{ev}_x(F \setminus \{ f \}, \alpha_i), \text{ev}_x(G, \alpha_i))$.

 \item If $f = u v$ for some $f \in F$ and $u, v \in \CC[S]$ (both not constant), then $X_1 = X(S, (F \setminus \{ f \}) \cup \{ u \}, G) = X \cap \{ u = 0 \}$ and $X_2 = X(S, (F \setminus \{ f \}) \cup \{ v \}, G \cup \{ u \}) = X \cap \{ u \ne 0, v = 0 \}$ define a stratification for $X$, and thus $[X] = [X_1] + [X_2]$.

 \item If $f = x u + v$ for some $f \in F$, $x \in S$ and $u, v \in \CC[S]$ with $x$ not appearing in $u$ and $v$, then we consider the following cases. For any point $p$ of $X$, either $u(p) = 0$, implying $v(p) = 0$ as well, or $u(p) \ne 0$, implying $x(p) = -v(p)/u(p)$. Therefore $[X] = [X_1] + [X_2]$ with $X_1 = X(S, (F \setminus \{ f \}) \cup \{ u, v \}, G)$ and $X_2 = X(S, \text{ev}_x(F \setminus \{ f \}, -v, u)), \text{ev}_x(G, -v, u) \cup \{ u \})$.

 \item Suppose $f = x^2 u + x v + w$ for some $f \in F$, $x \in S$ and $u, v, w \in \CC[S]$ with $x$ not appearing in $u, v$ and $w$. Moreover, suppose that the discriminant $D = v^2 - 4uw$ is a square, i.e., we can write $D = h^2$ for some $h \in \CC[S]$. Then for any point $p$ of $X$, we consider the following cases. Either $u(p) = 0$, in which case $(xv + w)(p) = 0$. If $u(p) \ne 0$, we distinguish between $D(p) = 0$ and $D(p) \ne 0$. In the first case we find that $x(p) = \left(\frac{-v}{2u}\right)(p)$, and in the latter case we have the two possibilities $x(p) = \left(\frac{-v \pm h}{2u}\right)(p)$. Hence $[X] = [X_1] + [X_2] + [X_3] + [X_4]$, with
 \begin{gather*}
 X_1= X(S, (F \setminus \{ f \}) \cup \{ u, x v + w \}, G) , \\
 X_2 = X(S, \text{ev}_x(F \setminus \{ f \}, -v, 2u) \cup \{ D \}, \text{ev}_x(G, -v, 2u) \cup \{ u \}) , \\
 X_3= X(S, \text{ev}_x(F \setminus \{ f \}, -v-h, 2u), \text{ev}_x(G, -v-h, 2u) \cup \{ u, D \}) , \\
 X_4= X(S, \text{ev}_x(F \setminus \{ f \}, -v+h, 2u), \text{ev}_x(G, -v+h, 2u) \cup \{ u, D \}) .
 \end{gather*}

 \item If $G \ne \varnothing$, pick any $g \in G$. We have $[X] = [X_1] - [X_2]$ where $X_1 = X(S, F, G \setminus \{ g \})$ and $X_2 = X(S, F \cup \{ g \}, G)$.
 \end{enumerate}
\end{Algorithm}

An implementation of this algorithm in Python can be found at~\cite{github}, together with the code for the computations done Sections~\ref{sec:application_U3} and~\ref{sec:application_U4}.

\section[Matrix Z(L) for U\_4]{Matrix $\boldsymbol{\tilde{Z}(L)}$ for $\boldsymbol{\U_4}$}\label{app:B}

The map $\tilde{Z}(L)$ for $G = \U_4$ is represented by the following matrix:

\scalebox{0.70}{\tiny
\[ \left[\begin{array}{c|cccc} & T_{1} & T_{2} & T_{3} & T_{4} \\ \hline
T_1 & q^{7} (q - 1)^{5} \left(q^{2} + 3 q - 2\right) & q^{9} (q - 2) (q - 1)^{5} \left(q^{2} + q - 1\right) & q^{8} (q - 1)^{5} \left(q^{3} - 3 q + 1\right) & q^{7} (q - 1)^{6} \left(q^{2} + 2 q - 2\right)\\ T_2 & q^{7} (q - 2) (q - 1)^{4} \left(q^{2} + q - 1\right) & q^{9} (q - 1)^{4} \left(q^{2} - 3 q + 3\right) \left(q^{2} + q - 1\right) & q^{8} (q - 2) (q - 1)^{5} \left(q^{2} + q - 1\right) & q^{7} (q - 2) (q - 1)^{5} \left(q^{2} + q - 1\right)\\ T_3 & q^{7} (q - 1)^{4} \left(q^{3} - 3 q + 1\right) & q^{9} (q - 2) (q - 1)^{5} \left(q^{2} + q - 1\right) & q^{8} (q - 1)^{4} \left(q^{4} - q^{3} - 2 q^{2} + 4 q - 1\right) & q^{7} (q - 1)^{5} \left(q^{3} - 3 q + 1\right)\\ T_4 & q^{7} (q - 1)^{5} \left(q^{2} + 2 q - 2\right) & q^{9} (q - 2) (q - 1)^{5} \left(q^{2} + q - 1\right) & q^{8} (q - 1)^{5} \left(q^{3} - 3 q + 1\right) & q^{7} (q - 1)^{5} \left(q^{3} + q^{2} - 3 q + 2\right)\\ T_5 & q^{6} (q - 2) (q - 1)^{4} \left(q^{3} + q - 1\right) & q^{11} (q - 2)^{2} (q - 1)^{4} & q^{10} (q - 2) (q - 1)^{5} & q^{6} (q - 2) (q - 1)^{6} \left(q^{2} + q + 1\right)\\ T_6 & q^{7} (q - 1)^{4} \left(q^{3} - 3 q + 1\right) & q^{9} (q - 2) (q - 1)^{6} \left(q + 1\right) & q^{8} (q - 1)^{4} \left(q^{4} - 2 q^{3} - q^{2} + 4 q - 1\right) & q^{7} (q - 1)^{5} \left(q^{3} - 3 q + 1\right) \\ T_7 & q^{7} (q - 2) (q - 1)^{4} \left(q^{2} + q - 1\right) & q^{11} (q - 2)^{2} (q - 1)^{4} & q^{8} (q - 2) (q - 1)^{6} \left(q + 1\right) & q^{7} (q - 2) (q - 1)^{5} \left(q^{2} + q - 1\right)\\ T_8 & q^{9} (q - 2)^{2} (q - 1)^{3} & q^{11} (q - 2) (q - 1)^{3} \left(q^{2} - 3 q + 3\right) & q^{10} (q - 2)^{2} (q - 1)^{4} & q^{9} (q - 2)^{2} (q - 1)^{4} \\ T_9 & q^{7} (q - 2) (q - 1)^{5} \left(q + 1\right) & q^{9} (q - 1)^{5} \left(q + 1\right) \left(q^{2} - 3 q + 3\right) & q^{8} (q - 2) (q - 1)^{6} \left(q + 1\right) & q^{7} (q - 2) (q - 1)^{6} \left(q + 1\right) \\ T_{10} & q^{7} (q - 2) (q - 1)^{4} \left(q^{2} - q - 1\right) & q^{11} (q - 2) (q - 1)^{3} \left(q^{2} - 3 q + 3\right) & q^{8} (q - 2) (q - 1)^{3} \left(q^{4} - 3 q^{3} + 2 q^{2} - 1\right) & q^{7} (q - 2) (q - 1)^{5} \left(q^{2} - q - 1\right) \\ T_{11} & q^{7} (q - 1)^{3} \left(q^{4} - 2 q^{3} - q^{2} + 4 q - 1\right) & q^{9} (q - 2) (q - 1)^{6} \left(q + 1\right) & q^{8} (q - 1)^{3} \left(q^{5} - 3 q^{4} + 2 q^{3} + 3 q^{2} - 5 q + 1\right) & q^{7} (q - 1)^{4} \left(q^{4} - 2 q^{3} - q^{2} + 4 q - 1\right) \\ T_{12} & q^{7} (q - 2) (q - 1)^{5} \left(q + 1\right) & q^{11} (q - 2)^{2} (q - 1)^{4} & q^{8} (q - 2) (q - 1)^{4} \left(q^{3} - q^{2} + 1\right) & q^{7} (q - 2) (q - 1)^{6} \left(q + 1\right) \\ T_{13} & q^{9} (q - 2)^{2} (q - 1)^{3} & q^{11} (q - 2)^{3} (q - 1)^{3} & q^{10} (q - 2)^{2} (q - 1)^{4} & q^{9} (q - 2)^{2} (q - 1)^{4} \\ T_{14} & q^{6} (q - 2) (q - 1)^{5} \left(q^{2} + q + 1\right) & q^{11} (q - 2)^{2} (q - 1)^{4} & q^{10} (q - 2) (q - 1)^{5} & q^{6} (q - 2) (q - 1)^{4} \left(q^{4} - q^{3} + 1\right) \\ T_{15} & q^{9} (q - 2)^{3} (q - 1)^{2} & q^{11} (q - 2)^{2} (q - 1)^{2} \left(q^{2} - 3 q + 3\right) & q^{10} (q - 2)^{3} (q - 1)^{3} & q^{9} (q - 2)^{3} (q - 1)^{3} \\ T_{16} & q^{7} (q - 2) (q - 1)^{2} \left(q^{4} - 3 q^{3} + 2 q^{2} - 1\right) & q^{11} (q - 2) (q - 1)^{3} \left(q^{2} - 3 q + 3\right) & q^{8} (q - 2) (q - 1)^{2} \left(q^{5} - 4 q^{4} + 5 q^{3} - 2 q^{2} + 1\right) & q^{7} (q - 2) (q - 1)^{3} \left(q^{4} - 3 q^{3} + 2 q^{2} - 1\right)
\end{array} \right.\]}

\scalebox{0.70}{\tiny
\[ \begin{array}{cccc} T_{5} & T_{6} & T_{7} & T_{8} \\ \hline
q^{8} (q - 2) (q - 1)^{5} \left(q^{3} + q - 1\right) & q^{8} (q - 1)^{5} \left(q^{3} - 3 q + 1\right) & q^{9} (q - 2) (q - 1)^{5} \left(q^{2} + q - 1\right) & q^{12} (q - 2)^{2} (q - 1)^{5}\\q^{11} (q - 2)^{2} (q - 1)^{4} & q^{8} (q - 2) (q - 1)^{6} \left(q + 1\right) & q^{11} (q - 2)^{2} (q - 1)^{4} & q^{12} (q - 2) (q - 1)^{4} \left(q^{2} - 3 q + 3\right)\\q^{11} (q - 2) (q - 1)^{5} & q^{8} (q - 1)^{4} \left(q^{4} - 2 q^{3} - q^{2} + 4 q - 1\right) & q^{9} (q - 2) (q - 1)^{6} \left(q + 1\right) & q^{12} (q - 2)^{2} (q - 1)^{5}\\q^{8} (q - 2) (q - 1)^{6} \left(q^{2} + q + 1\right) & q^{8} (q - 1)^{5} \left(q^{3} - 3 q + 1\right) & q^{9} (q - 2) (q - 1)^{5} \left(q^{2} + q - 1\right) & q^{12} (q - 2)^{2} (q - 1)^{5}\\q^{7} (q - 1)^{4} \left(q^{2} - 3 q + 3\right) \left(q^{4} + q - 1\right) & q^{10} (q - 2) (q - 1)^{5} & q^{11} (q - 2)^{2} (q - 1)^{4} & q^{12} (q - 2) (q - 1)^{4} \left(q^{2} - 3 q + 3\right)\\q^{11} (q - 2) (q - 1)^{5} & q^{8} (q - 1)^{4} \left(q^{4} - q^{3} - 2 q^{2} + 4 q - 1\right) & q^{9} (q - 2) (q - 1)^{5} \left(q^{2} + q - 1\right) & q^{12} (q - 2)^{2} (q - 1)^{5}\\q^{11} (q - 2)^{2} (q - 1)^{4} & q^{8} (q - 2) (q - 1)^{5} \left(q^{2} + q - 1\right) & q^{9} (q - 1)^{4} \left(q^{2} - 3 q + 3\right) \left(q^{2} + q - 1\right) & q^{12} (q - 2)^{3} (q - 1)^{4}\\q^{11} (q - 2) (q - 1)^{3} \left(q^{2} - 3 q + 3\right) & q^{10} (q - 2)^{2} (q - 1)^{4} & q^{11} (q - 2)^{3} (q - 1)^{3} & q^{12} (q - 1)^{3} \left(q^{2} - 3 q + 3\right)^{2}\\q^{11} (q - 2)^{2} (q - 1)^{4} & q^{8} (q - 2) (q - 1)^{4} \left(q^{3} - q^{2} + 1\right) & q^{11} (q - 2)^{2} (q - 1)^{4} & q^{12} (q - 2) (q - 1)^{4} \left(q^{2} - 3 q + 3\right)\\q^{11} (q - 2)^{3} (q - 1)^{3} & q^{8} (q - 2) (q - 1)^{3} \left(q^{4} - 3 q^{3} + 2 q^{2} - 1\right) & q^{11} (q - 2) (q - 1)^{3} \left(q^{2} - 3 q + 3\right) & q^{12} (q - 2)^{2} (q - 1)^{3} \left(q^{2} - 3 q + 3\right)\\q^{11} (q - 2) (q - 1)^{5} & q^{8} (q - 1)^{3} \left(q^{5} - 3 q^{4} + 2 q^{3} + 3 q^{2} - 5 q + 1\right) & q^{9} (q - 2) (q - 1)^{6} \left(q + 1\right) & q^{12} (q - 2)^{2} (q - 1)^{5}\\q^{11} (q - 2)^{2} (q - 1)^{4} & q^{8} (q - 2) (q - 1)^{6} \left(q + 1\right) & q^{9} (q - 1)^{5} \left(q + 1\right) \left(q^{2} - 3 q + 3\right) & q^{12} (q - 2)^{3} (q - 1)^{4}\\q^{11} (q - 2) (q - 1)^{3} \left(q^{2} - 3 q + 3\right) & q^{10} (q - 2)^{2} (q - 1)^{4} & q^{11} (q - 2) (q - 1)^{3} \left(q^{2} - 3 q + 3\right) & q^{12} (q - 2)^{2} (q - 1)^{3} \left(q^{2} - 3 q + 3\right)\\q^{7} (q - 1)^{5} \left(q + 1\right) \left(q^{2} + 1\right) \left(q^{2} - 3 q + 3\right) & q^{10} (q - 2) (q - 1)^{5} & q^{11} (q - 2)^{2} (q - 1)^{4} & q^{12} (q - 2) (q - 1)^{4} \left(q^{2} - 3 q + 3\right)\\q^{11} (q - 2)^{2} (q - 1)^{2} \left(q^{2} - 3 q + 3\right) & q^{10} (q - 2)^{3} (q - 1)^{3} & q^{11} (q - 2)^{2} (q - 1)^{2} \left(q^{2} - 3 q + 3\right) & q^{12} (q - 2) (q - 1)^{2} \left(q^{2} - 3 q + 3\right)^{2}\\q^{11} (q - 2)^{3} (q - 1)^{3} & q^{8} (q - 2) (q - 1)^{2} \left(q^{5} - 4 q^{4} + 5 q^{3} - 2 q^{2} + 1\right) & q^{11} (q - 2) (q - 1)^{3} \left(q^{2} - 3 q + 3\right) & q^{12} (q - 2)^{2} (q - 1)^{3} \left(q^{2} - 3 q + 3\right)
\end{array} \]}

\scalebox{0.70}{\tiny
\[ \begin{array}{cccc} T_{9} & T_{10} & T_{11} & T_{12}\\ \hline
q^{9} (q - 2) (q - 1)^{7} \left(q + 1\right) & q^{9} (q - 2) (q - 1)^{6} \left(q^{2} - q - 1\right) & q^{8} (q - 1)^{5} \left(q^{4} - 2 q^{3} - q^{2} + 4 q - 1\right) & q^{9} (q - 2) (q - 1)^{7} \left(q + 1\right)\\q^{9} (q - 1)^{6} \left(q + 1\right) \left(q^{2} - 3 q + 3\right) & q^{11} (q - 2) (q - 1)^{4} \left(q^{2} - 3 q + 3\right) & q^{8} (q - 2) (q - 1)^{7} \left(q + 1\right) & q^{11} (q - 2)^{2} (q - 1)^{5}\\q^{9} (q - 2) (q - 1)^{7} \left(q + 1\right) & q^{9} (q - 2) (q - 1)^{4} \left(q^{4} - 3 q^{3} + 2 q^{2} - 1\right) & q^{8} (q - 1)^{4} \left(q^{5} - 3 q^{4} + 2 q^{3} + 3 q^{2} - 5 q + 1\right) & q^{9} (q - 2) (q - 1)^{5} \left(q^{3} - q^{2} + 1\right)\\q^{9} (q - 2) (q - 1)^{7} \left(q + 1\right) & q^{9} (q - 2) (q - 1)^{6} \left(q^{2} - q - 1\right) & q^{8} (q - 1)^{5} \left(q^{4} - 2 q^{3} - q^{2} + 4 q - 1\right) & q^{9} (q - 2) (q - 1)^{7} \left(q + 1\right)\\q^{11} (q - 2)^{2} (q - 1)^{5} & q^{11} (q - 2)^{3} (q - 1)^{4} & q^{10} (q - 2) (q - 1)^{6} & q^{11} (q - 2)^{2} (q - 1)^{5}\\q^{9} (q - 2) (q - 1)^{5} \left(q^{3} - q^{2} + 1\right) & q^{9} (q - 2) (q - 1)^{4} \left(q^{4} - 3 q^{3} + 2 q^{2} - 1\right) & q^{8} (q - 1)^{4} \left(q^{5} - 3 q^{4} + 2 q^{3} + 3 q^{2} - 5 q + 1\right) & q^{9} (q - 2) (q - 1)^{7} \left(q + 1\right)\\q^{11} (q - 2)^{2} (q - 1)^{5} & q^{11} (q - 2) (q - 1)^{4} \left(q^{2} - 3 q + 3\right) & q^{8} (q - 2) (q - 1)^{7} \left(q + 1\right) & q^{9} (q - 1)^{6} \left(q + 1\right) \left(q^{2} - 3 q + 3\right)\\q^{11} (q - 2) (q - 1)^{4} \left(q^{2} - 3 q + 3\right) & q^{11} (q - 2)^{2} (q - 1)^{3} \left(q^{2} - 3 q + 3\right) & q^{10} (q - 2)^{2} (q - 1)^{5} & q^{11} (q - 2)^{3} (q - 1)^{4}\\q^{9} (q - 1)^{4} \left(q^{2} - 3 q + 3\right) \left(q^{3} - q^{2} + 1\right) & q^{11} (q - 2) (q - 1)^{4} \left(q^{2} - 3 q + 3\right) & q^{8} (q - 2) (q - 1)^{5} \left(q^{3} - q^{2} + 1\right) & q^{11} (q - 2)^{2} (q - 1)^{5}\\q^{11} (q - 2) (q - 1)^{4} \left(q^{2} - 3 q + 3\right) & q^{9} (q - 1)^{3} \left(q^{2} - 3 q + 3\right) \left(q^{4} - 3 q^{3} + 3 q^{2} + 1\right) & q^{8} (q - 2) (q - 1)^{3} \left(q^{5} - 4 q^{4} + 5 q^{3} - 2 q^{2} + 1\right) & q^{11} (q - 2) (q - 1)^{4} \left(q^{2} - 3 q + 3\right)\\q^{9} (q - 2) (q - 1)^{5} \left(q^{3} - q^{2} + 1\right) & q^{9} (q - 2) (q - 1)^{3} \left(q^{5} - 4 q^{4} + 5 q^{3} - 2 q^{2} + 1\right) & q^{8} (q - 1)^{3} \left(q^{6} - 4 q^{5} + 6 q^{4} - 2 q^{3} - 5 q^{2} + 6 q - 1\right) & q^{9} (q - 2) (q - 1)^{5} \left(q^{3} - q^{2} + 1\right)\\q^{11} (q - 2)^{2} (q - 1)^{5} & q^{11} (q - 2) (q - 1)^{4} \left(q^{2} - 3 q + 3\right) & q^{8} (q - 2) (q - 1)^{5} \left(q^{3} - q^{2} + 1\right) & q^{9} (q - 1)^{4} \left(q^{2} - 3 q + 3\right) \left(q^{3} - q^{2} + 1\right)\\q^{11} (q - 2)^{3} (q - 1)^{4} & q^{11} (q - 2)^{2} (q - 1)^{3} \left(q^{2} - 3 q + 3\right) & q^{10} (q - 2)^{2} (q - 1)^{5} & q^{11} (q - 2) (q - 1)^{4} \left(q^{2} - 3 q + 3\right)\\q^{11} (q - 2)^{2} (q - 1)^{5} & q^{11} (q - 2)^{3} (q - 1)^{4} & q^{10} (q - 2) (q - 1)^{6} & q^{11} (q - 2)^{2} (q - 1)^{5}\\q^{11} (q - 2)^{2} (q - 1)^{3} \left(q^{2} - 3 q + 3\right) & q^{11} (q - 2) (q - 1)^{2} \left(q^{2} - 3 q + 3\right)^{2} & q^{10} (q - 2)^{3} (q - 1)^{4} & q^{11} (q - 2)^{2} (q - 1)^{3} \left(q^{2} - 3 q + 3\right)\\q^{11} (q - 2) (q - 1)^{4} \left(q^{2} - 3 q + 3\right) & q^{9} (q - 1)^{2} \left(q^{2} - 3 q + 3\right) \left(q^{5} - 4 q^{4} + 6 q^{3} - 3 q^{2} - 1\right) & q^{8} (q - 2) (q - 1)^{2} \left(q^{6} - 5 q^{5} + 9 q^{4} - 7 q^{3} + 2 q^{2} - 1\right) & q^{11} (q - 2) (q - 1)^{4} \left(q^{2} - 3 q + 3\right)
\end{array} \]}

\scalebox{0.70}{\tiny
\[ \left.\begin{array}{cccc}T_{13} & T_{14} & T_{15} & T_{16}\\ \hline
q^{12} (q - 2)^{2} (q - 1)^{5} & q^{8} (q - 2) (q - 1)^{7} \left(q^{2} + q + 1\right) & q^{12} (q - 2)^{3} (q - 1)^{5} & q^{9} (q - 2) (q - 1)^{5} \left(q^{4} - 3 q^{3} + 2 q^{2} - 1\right)\\q^{12} (q - 2)^{3} (q - 1)^{4} & q^{11} (q - 2)^{2} (q - 1)^{5} & q^{12} (q - 2)^{2} (q - 1)^{4} \left(q^{2} - 3 q + 3\right) & q^{11} (q - 2) (q - 1)^{5} \left(q^{2} - 3 q + 3\right)\\q^{12} (q - 2)^{2} (q - 1)^{5} & q^{11} (q - 2) (q - 1)^{6} & q^{12} (q - 2)^{3} (q - 1)^{5} & q^{9} (q - 2) (q - 1)^{4} \left(q^{5} - 4 q^{4} + 5 q^{3} - 2 q^{2} + 1\right)\\q^{12} (q - 2)^{2} (q - 1)^{5} & q^{8} (q - 2) (q - 1)^{5} \left(q^{4} - q^{3} + 1\right) & q^{12} (q - 2)^{3} (q - 1)^{5} & q^{9} (q - 2) (q - 1)^{5} \left(q^{4} - 3 q^{3} + 2 q^{2} - 1\right)\\q^{12} (q - 2) (q - 1)^{4} \left(q^{2} - 3 q + 3\right) & q^{7} (q - 1)^{6} \left(q + 1\right) \left(q^{2} + 1\right) \left(q^{2} - 3 q + 3\right) & q^{12} (q - 2)^{2} (q - 1)^{4} \left(q^{2} - 3 q + 3\right) & q^{11} (q - 2)^{3} (q - 1)^{5}\\q^{12} (q - 2)^{2} (q - 1)^{5} & q^{11} (q - 2) (q - 1)^{6} & q^{12} (q - 2)^{3} (q - 1)^{5} & q^{9} (q - 2) (q - 1)^{4} \left(q^{5} - 4 q^{4} + 5 q^{3} - 2 q^{2} + 1\right)\\q^{12} (q - 2) (q - 1)^{4} \left(q^{2} - 3 q + 3\right) & q^{11} (q - 2)^{2} (q - 1)^{5} & q^{12} (q - 2)^{2} (q - 1)^{4} \left(q^{2} - 3 q + 3\right) & q^{11} (q - 2) (q - 1)^{5} \left(q^{2} - 3 q + 3\right)\\q^{12} (q - 2)^{2} (q - 1)^{3} \left(q^{2} - 3 q + 3\right) & q^{11} (q - 2) (q - 1)^{4} \left(q^{2} - 3 q + 3\right) & q^{12} (q - 2) (q - 1)^{3} \left(q^{2} - 3 q + 3\right)^{2} & q^{11} (q - 2)^{2} (q - 1)^{4} \left(q^{2} - 3 q + 3\right)\\q^{12} (q - 2)^{3} (q - 1)^{4} & q^{11} (q - 2)^{2} (q - 1)^{5} & q^{12} (q - 2)^{2} (q - 1)^{4} \left(q^{2} - 3 q + 3\right) & q^{11} (q - 2) (q - 1)^{5} \left(q^{2} - 3 q + 3\right)\\q^{12} (q - 2)^{2} (q - 1)^{3} \left(q^{2} - 3 q + 3\right) & q^{11} (q - 2)^{3} (q - 1)^{4} & q^{12} (q - 2) (q - 1)^{3} \left(q^{2} - 3 q + 3\right)^{2} & q^{9} (q - 1)^{3} \left(q^{2} - 3 q + 3\right) \left(q^{5} - 4 q^{4} + 6 q^{3} - 3 q^{2} - 1\right)\\q^{12} (q - 2)^{2} (q - 1)^{5} & q^{11} (q - 2) (q - 1)^{6} & q^{12} (q - 2)^{3} (q - 1)^{5} & q^{9} (q - 2) (q - 1)^{3} \left(q^{6} - 5 q^{5} + 9 q^{4} - 7 q^{3} + 2 q^{2} - 1\right)\\q^{12} (q - 2) (q - 1)^{4} \left(q^{2} - 3 q + 3\right) & q^{11} (q - 2)^{2} (q - 1)^{5} & q^{12} (q - 2)^{2} (q - 1)^{4} \left(q^{2} - 3 q + 3\right) & q^{11} (q - 2) (q - 1)^{5} \left(q^{2} - 3 q + 3\right)\\q^{12} (q - 1)^{3} \left(q^{2} - 3 q + 3\right)^{2} & q^{11} (q - 2) (q - 1)^{4} \left(q^{2} - 3 q + 3\right) & q^{12} (q - 2) (q - 1)^{3} \left(q^{2} - 3 q + 3\right)^{2} & q^{11} (q - 2)^{2} (q - 1)^{4} \left(q^{2} - 3 q + 3\right)\\q^{12} (q - 2) (q - 1)^{4} \left(q^{2} - 3 q + 3\right) & q^{7} (q - 1)^{4} \left(q^{2} - 3 q + 3\right) \left(q^{5} - q^{4} + 1\right) & q^{12} (q - 2)^{2} (q - 1)^{4} \left(q^{2} - 3 q + 3\right) & q^{11} (q - 2)^{3} (q - 1)^{5}\\q^{12} (q - 2) (q - 1)^{2} \left(q^{2} - 3 q + 3\right)^{2} & q^{11} (q - 2)^{2} (q - 1)^{3} \left(q^{2} - 3 q + 3\right) & q^{12} (q - 1)^{2} \left(q^{2} - 3 q + 3\right)^{3} & q^{11} (q - 2) (q - 1)^{3} \left(q^{2} - 3 q + 3\right)^{2}\\q^{12} (q - 2)^{2} (q - 1)^{3} \left(q^{2} - 3 q + 3\right) & q^{11} (q - 2)^{3} (q - 1)^{4} & q^{12} (q - 2) (q - 1)^{3} \left(q^{2} - 3 q + 3\right)^{2} & q^{9} (q - 1)^{2} \left(q^{2} - 3 q + 3\right) \left(q^{6} - 5 q^{5} + 10 q^{4} - 9 q^{3} + 3 q^{2} + 1\right)
\end{array}\right] \]
}

\subsection*{Acknowledgements}
The authors thank Bas Edixhoven and David Holmes for reading a previous version of this paper and giving valuable comments; Ángel Gonz\'alez-Prieto whose papers were the starting point of this research and who was kind enough to answer any questions; and Sean Lawton for pointing out an error regarding $\overline{\chi}$. The authors also thank the reviewers for their detailed feedback and invaluable comments. The paper is part of the master's thesis \cite{thesisvogel} of the second author.

\pdfbookmark[1]{References}{ref}
\LastPageEnding


\begin{thebibliography}{99}
\footnotesize\itemsep=0pt

\bibitem{atiyah}
Atiyah M., Topological quantum field theories, \href{https://doi.org/10.1016/0040-9383(94)00051-4}{\textit{Inst. Hautes \'Etudes
 Sci. Publ. Math.}} \textbf{68} (1988), 175--186.

\bibitem{baraglia_hekmati}
Baraglia D., Hekmati P., Arithmetic of singular character varieties and their
 $E$-polynomials, \href{https://doi.org/10.1112/plms.12008}{\textit{Proc. Lond. Math. Soc.}} \textbf{114} (2017),
 293--332, \href{https://arxiv.org/abs/1602.06996}{arXiv:1602.06996}.

\bibitem{bhunia}
Bhunia S., Conjugacy classes of centralizers in the group of upper triangular
 matrices, \href{https://doi.org/10.1142/S0219498820500085}{\textit{J.~Algebra Appl.}} \textbf{19} (2020), 2050008, 14~pages,
 \href{https://arxiv.org/abs/1901.07869}{arXiv:1901.07869}.

\bibitem{Borisov2014}
Borisov L.A., The class of the affine line is a zero divisor in the
 {G}rothendieck ring, \href{https://doi.org/10.7282/T33B62H9}{\textit{J.~Algebraic Geom.}} \textbf{27} (2014),
 203--209, \href{https://arxiv.org/abs/1412.6194}{arXiv:1412.6194}.

\bibitem{brown}
Brown R., Groupoids and van {K}ampen's theorem, \href{https://doi.org/10.1112/plms/s3-17.3.385}{\textit{Proc. London Math.
 Soc.}} \textbf{17} (1967), 385--401.

\bibitem{culler_and_shalen}
Culler M., Shalen P.B., Varieties of group representations and splittings of
 $3$-manifolds, \href{https://doi.org/10.2307/2006973}{\textit{Ann. of Math.}} \textbf{117} (1983), 109--146.

\bibitem{deligne1974theorie}
Deligne P., Th\'eorie de {H}odge.~{III}, \href{https://doi.org/10.1007/BF02685881}{\textit{Inst. Hautes \'Etudes Sci.
 Publ. Math.}} \textbf{44} (1974), 5--77.

\bibitem{flr2017}
Florentino C., Lawton S., Ramras D., Homotopy groups of free group character
 varieties, \href{https://doi.org/10.2422/2036-2145.201510_004}{\textit{Ann. Sc. Norm. Super. Pisa Cl. Sci.}} \textbf{17} (2017),
 143--185, \href{https://arxiv.org/abs/1412.0272}{arXiv:1412.0272}.

\bibitem{florentino2021hodge}
Florentino C., Silva J., Hodge-{D}eligne polynomials of character varieties of
 free abelian groups, \href{https://doi.org/10.1515/math-2021-0038}{\textit{Open Math.}} \textbf{19} (2021), 338--362,
 \href{https://arxiv.org/abs/1711.07909}{arXiv:1711.07909}.

\bibitem{frobenius1896gruppencharaktere}
Frobenius G., \"Uber {G}ruppencharaktere, Wissenschaften Berlin,
 Sitzungsberichte der Preu{\ss}ischen Akademie der Wissenschaften zu Berlin,
 Reichsdr., 1896.

\bibitem{arXiv181009714}
Gonz\'alez-Prieto A., Motivic theory of representation varieties via
 topological quantum field theories, \href{https://arxiv.org/abs/1810.09714}{arXiv:1810.09714}.

\bibitem{thesisangel}
Gonz\'alez-Prieto A., Topological quantum field theories for character
 varieties, Ph.D.~Thesis, {U}niversidad Complutense de Madrid, 2018,
 \href{https://arxiv.org/abs/1812.11575}{arXiv:1812.11575}.

\bibitem{arXiv190605222}
Gonz\'alez-Prieto A., Virtual classes of parabolic {${\rm SL}_2(\mathbb
 C)$}-character varieties, \href{https://doi.org/10.1016/j.aim.2020.107148}{\textit{Adv. Math.}} \textbf{368} (2020), 107148,
 41~pages, \href{https://arxiv.org/abs/1906.05222}{arXiv:1906.05222}.

\bibitem{arXiv170905724}
Gonz\'alez-Prieto A., Logares M., Mu\~noz V., A lax monoidal topological
 quantum field theory for representation varieties, \href{https://doi.org/10.1016/j.bulsci.2020.102871}{\textit{Bull. Sci. Math.}}
 \textbf{161} (2020), 102871, 34~pages, \href{https://arxiv.org/abs/1709.05724}{arXiv:1709.05724}.

\bibitem{arXiv200501841}
Gonz\'alez-Prieto A., Logares M., Mu\~noz V., Representation variety for the
 rank one affine group, in Mathematical analysis in interdisciplinary
 research, \textit{Springer Optim. Appl.}, Vol.~179, \href{https://doi.org/10.1007/978-3-030-84721-0_18}{Springer}, Cham, 2021,
 381--416, \href{https://arxiv.org/abs/2005.01841}{arXiv:2005.01841}.

\bibitem{hausel_villegas}
Hausel T., Rodriguez-Villegas F., Mixed {H}odge polynomials of character
 varieties (with an appendix by {N}icholas {M}.~{K}atz), \href{https://doi.org/10.1007/s00222-008-0142-x}{\textit{Invent.
 Math.}} \textbf{174} (2008), 555--624, \href{https://arxiv.org/abs/math.AG/0612668}{arXiv:math.AG/0612668}.

\bibitem{higman1960enumerating}
Higman G., Enumerating $p$-groups. {I}.~{I}nequalities, \href{https://doi.org/10.1112/plms/s3-10.1.24}{\textit{Proc. London
 Math. Soc.}} \textbf{3} (1960), 24--30.

\bibitem{hitchin}
Hitchin N.J., The self-duality equations on a {R}iemann surface, \href{https://doi.org/10.1112/plms/s3-55.1.59}{\textit{Proc.
 London Math. Soc.}} \textbf{55} (1987), 59--126.

\bibitem{kock}
Kock J., Frobenius algebras and 2{D} topological quantum field theories,
 \textit{London Math. Soc. Stud. Texts}, Vol.~59, \href{https://doi.org/10.1017/CBO9780511615443}{Cambridge University Press},
 Cambridge, 2003.

\bibitem{lawsik2017}
Lawton S., Sikora A.S., Varieties of characters, \href{https://doi.org/10.1007/s10468-017-9679-y}{\textit{Algebr. Represent.
 Theory}} \textbf{20} (2017), 1133--1141, \href{https://arxiv.org/abs/1604.02164}{arXiv:1604.02164}.

\bibitem{letellier2020series}
Letellier E., Rodriguez-Villegas F., $E$-series of character varieties of
 non-orientable surfaces, \href{https://arxiv.org/abs/2008.13435}{arXiv:2008.13435}.

\bibitem{arXiv11066011}
Logares M., Mu\~noz V., Newstead P.E., Hodge polynomials of {${\rm
 SL}(2,\mathbb{C})$}-character varieties for curves of small genus,
 \href{https://doi.org/10.1007/s13163-013-0115-5}{\textit{Rev. Mat. Complut.}} \textbf{26} (2013), 635--703.

\bibitem{Martin2016}
Martin N., The class of the affine line is a zero divisor in the {G}rothendieck
 ring: an improvement, \href{https://doi.org/10.1016/j.crma.2016.05.016}{\textit{C.~R.~Math. Acad. Sci. Paris}} \textbf{354}
 (2016), 936--939, \href{https://arxiv.org/abs/1604.06703}{arXiv:1604.06703}.

\bibitem{arXiv14076975}
Mart\'inez J., Mu\~noz V., $E$-polynomials of the {${\rm SL}(2,\mathbb
 C)$}-character varieties of surface groups, \href{https://doi.org/10.1093/imrn/rnv163}{\textit{Int. Math. Res. Not.}}
 \textbf{2016} (2016), 926--961, \href{https://arxiv.org/abs/1705.04649}{arXiv:1705.04649}.

\bibitem{mellit2020poincare}
Mellit A., Poincar\'e polynomials of character varieties, {M}acdonald
 polynomials and affine {S}pringer fibers, \href{https://doi.org/10.4007/annals.2020.192.1.3}{\textit{Ann. of Math.}} \textbf{192}
 (2020), 165--228, \href{https://arxiv.org/abs/1710.04513}{arXiv:1710.04513}.

\bibitem{mereb}
Mereb M., On the {$E$}-polynomials of a family of character varieties, Ph.D.~Thesis, {T}he University of Texas at Austin, 2010, \href{https://arxiv.org/abs/1006.1286}{arXiv:1006.1286}.

\bibitem{milnor}
Milnor J., Lectures on the $h$-cobordism theorem, Princeton University Press,
 Princeton, N.J., 1965.

\bibitem{nagata}
Nagata M., Invariants of a group in an affine ring, \href{https://doi.org/10.1215/kjm/1250524787}{\textit{J.~Math. Kyoto
 Univ.}} \textbf{3} (1964), 369--377.

\bibitem{newstead_intro}
Newstead P.E., Introduction to moduli problems and orbit spaces, \textit{Tata
 Institute of Fundamental Research Lectures on Mathematics and Physics},
 Vol.~51, Tata Institute of Fundamental Research, Bombay, Narosa Publishing
 House, New Delhi, 1978.

\bibitem{sikora2014character}
Sikora A.S., Character varieties of abelian groups, \href{https://doi.org/10.1007/s00209-013-1252-8}{\textit{Math.~Z.}}
 \textbf{277} (2014), 241--256, \href{https://arxiv.org/abs/1207.5284}{arXiv:1207.5284}.

\bibitem{simpson92}
Simpson C.T., Higgs bundles and local systems, \href{https://doi.org/10.1007/bf02699491}{\textit{Inst. Hautes \'Etudes
 Sci. Publ. Math.}} \textbf{75} (1992), 5--95.

\bibitem{simpson94_II}
Simpson C.T., Moduli of representations of the fundamental group of a smooth
 projective variety.~{II}, \href{https://doi.org/10.1007/bf02698895}{\textit{Inst. Hautes \'Etudes Sci. Publ. Math.}}
 \textbf{80} (1994), 5--79.

\bibitem{thesisvogel}
Vogel J., Computing virtual classes of representation varieties using {TQFTs},
 {M}aster's thesis, Leiden University, 2020.

\bibitem{github}
Vogel J., Grothendieck ring, available at
 \url{https://github.com/jessetvogel/grothendieck-ring}.

\bibitem{witten}
Witten E., Topological quantum field theory, \href{https://doi.org/10.1007/BF01223371}{\textit{Comm. Math. Phys.}}
 \textbf{117} (1988), 353--386.

\bibitem{witten1991quantum}
Witten E., On quantum gauge theories in two dimensions, \href{https://doi.org/10.1007/BF02100009}{\textit{Comm. Math.
 Phys.}} \textbf{141} (1991), 153--209.

\end{thebibliography}
\end{document}